\documentclass{article}

\usepackage[english]{babel}
\usepackage{lscape}

\usepackage[letterpaper,top=2cm,bottom=2cm,left=2cm,right=2cm,marginparwidth=1.75cm]{geometry}

\makeatletter
\renewcommand{\maketitle}{\bgroup\setlength{\parindent}{0pt}
\begin{flushleft}
  \textbf{\huge \@title}\vspace{10pt}  
    
  \@author
\end{flushleft}\egroup
}
\makeatother

\usepackage{amsmath}
\usepackage{amssymb}
\usepackage[demo]{graphicx}
\usepackage[colorlinks=true, allcolors=blue]{hyperref}
\usepackage{physics}
\usepackage{caption}
\usepackage{subcaption}
\usepackage{comment}
\usepackage{float}

\usepackage{makecell}
\usepackage{multicol}
\usepackage{multirow}

\title{Mathematical Modelling of Mechanotransduction via RhoA Signalling Pathway}
\author{
Sofie Verhees$^{1}$, Chandrasekhar Venkataraman$^{2}$, Mariya Ptashnyk$^{1}$\\
$^{1}$Department of Mathematics, Heriot-Watt University, The Maxwell Institute for Mathematical Sciences, Edinburgh, UK; $^{2}$Department of Mathematics, University of Sussex, Brighton, UK}
\date{}

\begin{document}
\maketitle

\section*{Abstract}
We derive and simulate a mathematical model for mechanotransduction related to the Rho GTPase signalling pathway. The model addresses the bidirectional coupling between signalling processes and cell mechanics. A numerical method based on bulk-surface finite elements is proposed for the approximation of the  coupled system of nonlinear reaction-diffusion equations, defined inside the cell and on the cell membrane, and the equations of elasticity. Our simulation results illustrate novel emergent features such as the strong dependence of the dynamics on cell shape, a threshold-like response to changes in substrate stiffness, and the fact that coupling mechanics and signalling can lead to the robustness of cell deformation to larger changes in substrate stiffness, ensuring mechanical homeostasis in agreement with experiments.

\section*{Introduction}
Intercellular signalling processes constitute the mechanisms through which cells communicate with and respond to their environment. Hence, signalling pathways are important in all physiological activities of the cell, such as cell division, cell movement, the immune response, and tissue development ~\cite{Tomar_2009}. Aberrant cell signalling  can often result in the development of diseases ~\cite{Valls_2022}. It is therefore important to understand signalling phenomena. Recent studies have found that alongside biochemical reactions, mechanics plays an important role in many signalling pathways ~\cite{romani_crosstalk_2021, Cai_2021}.  This phenomena is referred to as mechanotransduction which, broadly speaking, is  any process by which cells convert mechanical stimuli into chemical signals ~\cite{humphrey_mechanotransduction_2014, saraswathibhatla_cellextracellular_2023}. 

A large number of recent works study the role of Rho GTPases, primarily RhoA,   in mechanotransduction in relation to  different mechanical cues: extracellular matrix (ECM) stiffness and viscoelasticity, tensile stress (stretching), compressive stress (compression), and shear stress (fluid flow shear), see e.g.~\cite{burridge_mechanotransduction_2019,  xie_cell_2023} for a review. Moreover, the coupling between biochemistry and mechanics is bidirectional, i.e., chemical signals can also affect the mechanical properties of the cell, such as molecules like  focal adhesion kinases (FAKs) that influence F-actin dynamics and therefore the stiffness of the cell~\cite{martino_cellular_2018, saraswathibhatla_cellextracellular_2023, Sun_2016, Young_2023}. 

The formidable complexity of the phenomena involved in mechanotransduction means that much about how the mechanics and the chemical processes of the cell communicate is not yet understood and mathematical modelling is crucial in this regard. Whilst the mathematical modelling of biochemical cell signalling processes   is fairly well developed,  e.g., ~\cite{gilbert2006computational, garcia_2014, ptashnyk_multiscale_2020}, the study of mechanotransduction is comparatively more recent, see~\cite{cheng_cellular_2017} for a review.  Typically the modelling involves solving coupled systems of partial differential equations  (PDEs) with reaction-diffusion equations modelling the biochemistry coupled to equations based on (visco)elastic constitutive laws for the mechanics. The progress of such efforts has been rapid, ranging from early models employing simplifications such as one-dimensional geometries ~\cite{besser_coupling_2007, novev_spatiotemporal_2021} to full three-dimensional simulations ~\cite{scott_spatial_2021} using advanced computational techniques.
Alongside continuum models,  a number of recent works have employed discrete approaches such as spring-based models ~\cite{kang_structurally_2015},  or models that employ a  Potts formalism~\cite{bar-ziv_pearling_1999,vianay_single_2010,albert_dynamics_2014}.
Despite this rapid progress, the existing models typically make major simplifying assumptions such as assuming a constant stiffness of the ECM~\cite{scott_spatial_2021,sun_computational_2016, eroume_exploring_2021}, as well as neglecting the two-way coupling in which signalling pathways affect the mechanics alongside mechanical cues inducing signalling processes. 

In the present work, we seek to develop,  analyse and simulate a model for mechanotransduction through the Rho GTPase signalling pathway which allows for a two-way coupling between the mechanics and the biochemistry.  The dynamics of the signalling molecules FAK and RhoA are modelled using reaction-diffusion equations, where the ECM stiffness and elastic stresses of the cell activate FAK. Under simplifying assumptions, i.e., assuming no dependence on the cell elastic stresses, the biochemical component of the model is derived as a reduction of the model proposed in~\cite{scott_spatial_2021}.  For the cell's mechanical properties, we assume an elastic constitutive relationship~\cite{gould_introduction_2013} and allow the material properties to depend on the concentrations of the signalling molecules.  We propose a numerical method based on bulk and surface finite elements~\cite{dziuk2013finite} for the approximation of the model equations. 
 
The results presented here show that our model can reproduce the qualitative results of~\cite{scott_spatial_2021}, i.e., the mass of activated FAK and RhoA depend on ECM stiffness, with the dependence captured well by a Hill function. On the inclusion of the two-way coupling between signalling processes and cell mechanics, we observe novel dynamics, such as the conservation of cell deformation under different values of the ECM stiffness, which underlines the importance of including these more complex models of the mechanics.
The role of mechanotransduction in homeostasis in biological processes has been discussed in a number of biological works, e.g.,  \cite{humphrey_mechanotransduction_2014, Cai_2021, Gilbert_2017, martino_cellular_2018, Sun_2016} and our work presents a concrete example of how modelling can help elucidate potential mechanisms that underlay the mechanical homeostasis. Homeostasis of cell deformation,  as observed in simulations of our model,  has been observed experimentally~\cite{Grolleman_2023}. Our focus is on elastic constitutive assumptions for the mechanics of the cell to enhance clarity of exposition and to avoid unnecessary technical complexities. This can be extended to allow for other constitutive laws such as viscoelasticity of the cell and/or of the ECM as has been done elsewhere in the literature in simpler settings in 1D~\cite{besser_coupling_2007, mcnicol_theoretical_2025}. 
This work thus serves as a starting point  in  modelling and analysis of the two-way coupling between mechanics and chemistry. 

The paper is organised as follows. We first derive the reduced model for the Rho GTPase signalling pathway, based on the model proposed in~\cite{scott_spatial_2021}.  Next, the mathematical model for the mechanotransduction related to the Rho GTPase signalling pathway is derived. Simulations of the model are presented in the section after.
We conclude the paper with a discussion of the results. Details on the numerical method applied to simulate the model equations are given in \nameref{S1_Appendix}, Section A.3.

\section*{Methods}
\subsection*{A Mathematical model for the Rho GTPase signalling pathway}\label{sec:model}
One of the main signaling pathways involved in mechanostranduction is  the Rho GTPase pathway, responsible for many important cellular processes, e.g.~motility, cell adhesion, polarisation, differentiation, remodelling of the exoskeleton, and the ECM 
~\cite{xie_cell_2023}. The RhoA signalling pathway is activated through the activation of FAK in response to tension on integrins,  which depends on  ECM stiffness~\cite{Tomar_2009, Young_2023}. 

Our model for mechanotransduction related to the RhoA-mediated intercellular signalling pathway is based on models developed in~\cite{scott_spatial_2021} and ~\cite{eroume_exploring_2021, sun_computational_2016}.  To incorporate the interactions between mechanics and signalling processes, we extend the model proposed in ~\cite{scott_spatial_2021} by considering elastic deformations of the cell.  Activated FAK is downstream in the RhoA GTPase signalling pathway and hence the activation of RhoA is a function of activated FAK.
The activation of RhoA results in ECM remodelling and deposition of new fibres, increasing ECM stiffness and hence activation of FAK~\cite{Di_2023}.  FAK is expressed in the cytoplasm of the cell and is activated on the cell membrane.  
To simplify the model and focus only on the most significant aspects from the perspective of qualitative behaviour, we reduce the model for the RhoA signalling pathway of~\cite{scott_spatial_2021} that includes the dynamics of  FAK, RhoA, ROCK, Myo, LIMK, mDia, Cofilin, F-actin and YAP/TAZ by considering only the dynamics of  FAK and activated RhoA. Such a reduction is possible since other molecules considered in the full model of~\cite{scott_spatial_2021} do not influence the dynamics of FAK and RhoA. Our rationale behind considering a simplified model is to more clearly elucidate the emergent features that arise when mechanics is coupled with signalling. It is not challenging to incorporate other biochemical species or different reaction kinetics within the framework we propose.  

We let $Y\subset \mathbb{R}^3$ denote the cytoplasm and $\Gamma=\partial Y$ the cell membrane.  We denote by    $\phi_d$ and $\phi_a$ the concentrations of inactive and active FAK,  and  by $\rho_a$ the concentration of active RhoA. We recall that FAK, both active and inactive, is assumed to be cytoplasm resident and activated RhoA membrane resident. Our reduced model for the biochemistry consists of the following system of reaction-diffusion equations
\begin{eqnarray}\label{reduced_RhoA_model}
\begin{aligned}
    \partial_t \phi_d - D_1\Delta \phi_d &= \,\, k_1\phi_a &&\mbox{in }Y, \; t>0, \\
     \partial_t \phi_a - D_2\Delta \phi_a &= -k_1\phi_a &&\mbox{in }Y, \; t>0, \\
    D_1\nabla \phi_d \cdot \nu &= -n_rk_2\phi_d-n_rk_3\frac{E}{C+E}\phi_d &&\mbox{on }\Gamma,\;  t>0, \\
    D_2\nabla \phi_a \cdot \nu &= \, \, n_rk_2\phi_d+ n_rk_3\frac{E}{C+E}\phi_d &&\mbox{on }\Gamma,\;  t>0,\\
    \partial_t  \rho_a - D_3\Delta_\Gamma \rho_a &= -k_4\rho_a+  n_rk_5\left((\gamma\phi_a)^n+1\right)\Big(\frac{M_\rho}{|Y|}-\frac{\rho_a}{n_r}\Big) && \mbox{on }\Gamma, \; t>0,\\
    \phi_d (0, x) = \phi_d^0&(x), \quad \phi_a (0, x) = \phi_a^0(x)  \quad \mbox{ in } Y, \qquad \rho_a (0, x) = \rho_a^0(x)  &&\mbox{on } \Gamma, 
\end{aligned}
\end{eqnarray}
where $\Delta_\Gamma$ is the Laplace Beltrami operator modelling diffusion on the surface  $\Gamma$, see e.g.~\cite{dziuk2013finite}, $n_r=|Y|/|\Gamma|$ is the ratio between the volume of the cytoplasm and the area of the cell membrane, $k_1, k_4$ are deactivation and $k_2, k_3, k_5$ are activation constants, $E$ is the substrate stiffness, $D_1,D_2,D_3$ are the diffusion constants, $C$, $n$ and $\gamma$ are positive constants, $\phi_d^0(x)$, $\phi_a^0(x)$ and $\rho_a^0(x)$ are bounded nonnegative functions, and $\tfrac{M_\rho}{|Y|}-\tfrac{\rho_a}{n_r}$ is an approximation for deactivated RhoA ($\rho_d$) on the surface with $M_\rho=\int_Y\rho_d^0\dd{x}+\int_\Gamma\rho_a^0\dd{s}$ the total mass of RhoA, a quantity conserved in the full model of \cite{scott_spatial_2021} and assumed to be conserved here. Simulations illustrating the agreement between results obtained using the reduced model~\eqref{reduced_RhoA_model} with those of \cite{scott_spatial_2021} for the full model are presented in \nameref{S1_Appendix}, Section A.1. 

\subsection*{Mathematical model for mechanotransduction}\label{sec:modelelasticity}
As a starting point for the mechanics, we consider small deformations and hence, assume a linear elastic constitutive law for the mechanics of the cell. Although viscoelastic or poroelastic behaviour of cells is proposed in many works \cite{kasza_cell_2007, moeendarbary_cytoplasm_2013}, linear elasticity is often chosen for modelling simplicity and it can yield results consistent with experimental observations \cite{banerjee_controlling_2013, oakes_geometry_2014, chojowski_reversible_2020}.  Our focus is modelling a bidirectional coupling between cell stiffness and signalling processes, a minimal model assuming a linear elastic law for the cell mechanics is therefore sufficient for this work, as it avoids unnecessary complications that arise from the consideration of viscous stresses. To demonstrate that our results remain relevant under more complex assumptions, a linear viscoelastic model is presented and simulated in~Section A.7 in \nameref{S1_Appendix}. We note that for the corresponding simulation results, the qualitative behaviour in both cases, linear elastic and linear viscoelastic, is the same.

An important simplification that arises under the small deformations assumption  inherent in this work  is that the dynamics of the signalling molecules may be effectively considered on the reference configuration with no additional terms arising due to the deformation. Models where the assumption of small deformations is relaxed will be addressed in future studies. 

The cell nucleus plays an important role in governing the mechanical properties of the cell~\cite{chengPredicting2023, Graham_2016}, whilst we predominantly neglect this in the present work, in \nameref{S1_Appendix}, Section A.5, we have included simulations of the model with a `passive' nucleus that is considered to be more rigid than the cytoplasm.

It has been shown that the stiffness of the cell increases as F-actin increases~\cite{scott_spatial_2021}. F-actin is downstream from activated FAK and, as apparent from the results in \cite{scott_spatial_2021}, we can use activated FAK as a proxy for F-actin. Therefore, we assume that the Young's modulus $E_c$ of the cell is a function of the activated FAK concentration. Based on experimental observations \cite{gardel_elastic_2004} and numerical simulations \cite{scott_spatial_2021}, we propose 
\begin{eqnarray}\label{eq:Ec}
    E_c = E_c(\phi_a) = k_7(1 + (k_8\phi_a)^{p}),
\end{eqnarray}
where $k_7$, $k_8$ and $p$ are non-negative constants.
Then for elastic deformations of the cell, we have 
\begin{eqnarray} \label{eq:coupled_elast}
      -\nabla \cdot \sigma(u) = 0  \quad \text{in } Y, 
\end{eqnarray}
with  
$$
 \sigma(u) = \lambda(\phi_a)(\nabla \cdot u)I +  2\mu(\phi_a) \big( \nabla u + (\nabla u)^T\big)
$$
and the Lame constants $\lambda$ and $\mu$ are given by 
$$
    \lambda(\phi_a) = \frac{E_c(\phi_a) \nu_c}{(1+\nu_c)(1-2\nu_c)},\quad 
    \mu(\phi_a) = \frac{E_c(\phi_a)}{2(1+\nu_c)},
$$
where $\nu_c$ is the Poisson ratio of the cell. 

Activated RhoA regulates remodellling of stress fibres inside the cell and stabilisation of actin filaments \cite{Burridge_2016,  Doyle_2015, Zhao_2007, chrzanowska-wodnicka_rho-stimulated_1996}. This mechanism is modelled by the stress on the boundary being dependent on activated RhoA concentration 
\begin{eqnarray} \label{eq:coupled_force_boundary}
    \sigma(u) \nu =k_6\mathbb{P}(\rho_a \nu) \quad \text{on }\Gamma, 
\end{eqnarray}
where  $k_6$ is a positive constant, $\mathbb{P}$ is a projection on  the space orthogonal to the space of rigid deformations, i.e.~rotations and translations. 
Alongside models where the cell is allowed to deform freely, to model a typical experimental set-up where cells are placed on a rigid substrate, we consider 
\begin{eqnarray}\label{eq:coupled_normal_zero}
    u \cdot \nu =0,  \qquad   \Pi_\tau (\sigma(u)\nu) = 0 \qquad \text{ on }\Gamma_0, 
\end{eqnarray}
 together with condition Eq~(\ref{eq:coupled_force_boundary}) on $\Gamma \setminus \Gamma_0$, where $\Pi_\tau (w)= w - (w\cdot \nu) \nu $ denotes  the tangential projection of  vector $w$. In this work, we choose $\Gamma_0=\Gamma \cap \lbrace x\in \mathbb{R}^3|x_3=0\rbrace$ to simulate a rigid substrate. 

It has been shown that an increased contractility is associated with increased activated FAK, see e.g.~\cite{Burridge_2016}. Thus we assume that FAK is activated by the stress of the cell and as a proxy for the cytosolic stress we use the positive part of trace of the Cauchy stress tensor ${\rm tr}(\sigma)_+$, where   ${\rm tr}(\sigma)$ is the first stress invariant and the positive part  reflects
the fact that extension rather than compression causes the activation of FAK. This modifies the system in Eq~(\ref{reduced_RhoA_model}) to 
\begin{eqnarray} \label{eq:coupled_reactions}
\begin{aligned}
    \partial_t \phi_d - D_1\Delta \phi_d &= \phantom{-} k_1\phi_a - C_1\mbox{tr}(\sigma)_+\phi_d &&\mbox{in }Y, \;t>0,  \\
    \partial_t \phi_a - D_2\Delta \phi_a &= -k_1\phi_a + C_1\mbox{tr}(\sigma)_+\phi_d &&\mbox{in }Y,\; t>0,  \\
       D_1\nabla \phi_d \cdot \nu &= -n_rk_2\phi_d-n_rk_3\frac{E}{C+E}\phi_d &&\mbox{on }\Gamma,\; t>0,  \\
    D_2\nabla \phi_a \cdot \nu &= \phantom{-} n_rk_2\phi_d+n_rk_3\frac{E}{C+E}\phi_d &&\mbox{on }\Gamma, \; t>0,  \\
    \partial_t  \rho_a - D_3\Delta_\Gamma \rho_a &= -k_4\rho_a+ n_rk_5 \big((\gamma\phi_a)^n+1\big)\Big(\frac{M_\rho}{|Y|}-\frac{\rho_a}{n_r}\Big) && \mbox{on }\Gamma,\; t>0,\\
    \phi_d (0, x) = \phi_d^0&(x), \quad \phi_a (0, x) = \phi_a^0(x)  \quad \mbox{ in } Y, \qquad \rho_a (0, x) = \rho_a^0(x)  &&\mbox{on } \Gamma, 
    \end{aligned} 
\end{eqnarray}
where $v_+ = \max\{ v, 0\}$. We can prove existence, uniqueness and boundedness of solutions to the system in Eqs~(\ref{eq:coupled_elast})-(\ref{eq:coupled_reactions})  which we intend to report on elsewhere. 

\subsection*{Initial and boundary conditions and parametrization}
Using the model in Eqs~(\ref{eq:coupled_elast})-(\ref{eq:coupled_reactions}) we investigate different scenarios demonstrating the interactions between mechanics and signalling processes. First,  we consider the impact of the cell Young's modulus  $E_c$  and compare the dynamics when considering  a constant $E_c$ versus the case where $E_c$ depends on activated FAK as defined in Eq~(\ref{eq:Ec}). We also model the effect of the stress on the signalling molecules FAK and simulate equations Eq~(\ref{eq:coupled_reactions}) for $C_1=0$~$(\mbox{kPa s})^{-1}$ and  $C_1=0.1$~$(\mbox{kPa s})^{-1}$, respectively. Additionally we consider two experimental scenarios: (i) the cell is placed on a rigid substrate, modelled by the boundary conditions Eq~(\ref{eq:coupled_force_boundary})  on $\Gamma \setminus \Gamma_0$ and Eq~(\ref{eq:coupled_normal_zero}) on $\Gamma_0$ or (ii) the cell is embedded in an agar substrate and we apply the force boundary condition Eq~(\ref{eq:coupled_force_boundary}) on the entire cell membrane. We also distinguish between two different stimuli, similar to~\cite{scott_spatial_2021}, (i) the so called `$2$xD stimulus', where the 
substrate stiffness is only applied to the bottom of the cell, i.e.\ $E$ is nonzero only on $\Gamma_0$, and  (ii) the `$3$D stimulus' where the cell is embedded in an agar (substrate) and the impact of the substrate stiffness on the signalling processes is considered on the whole cell membrane.  To analyse the impact of the cell shape on the dynamics of signalling molecules and mechanical deformations we  consider both axisymmetric cells and  polarised cells with a lamellipodium like shape. The diameter of the cell is larger for the lamellipodium cells such that the volume is similar to the axisymmetric cells. 

The parameters are chosen as in Table~\ref{tab:parameters}. For numerical simulations, we use a Finite Element Method to discretize in space and a semi-implicit Euler method to discretize in time, with the mesh size $h=2.94~\mu$m and time step $\Delta t=0.5$~s. Details on the numerical scheme and benchmark computations demonstrating the accuracy of the approach for a problem with a known solution are given in \nameref{S1_Appendix}, Section A.3.

\begin{table}[!ht]
  \flushleft
  \caption{
  {\bf Parameter values for  simulations of the model in Eqs~(\ref{eq:coupled_elast})-(\ref{eq:coupled_reactions}).}}

\begin{subtable}{0.38\textwidth}
    \subcaption{Parameters inherited from the model in Eq~\eqref{reduced_RhoA_model} that are identical to \cite{scott_spatial_2021}. Because the goal is to compare and extend upon the model of \cite{scott_spatial_2021}, we choose the values to be the same as \cite{scott_spatial_2021}, which are based  on literature and fitting to data. The exception are $D_1$ and $D_2$, which is discussed in \nameref{S1_Appendix}, Section A.1. }
    \begin{tabular}{|l|l|l|}
  \hline
    \multicolumn{2}{|l|}{Parameters} & Value \\ \hline
    \multicolumn{2}{|l|}{$\phi_d^0$} & $0.7~\mu\mbox{mol}/\mbox{dm}^3$  \\ \hline
    \multicolumn{2}{|l|}{$\phi_a^0$} & $0.3~\mu\mbox{mol}/\mbox{dm}^3$  \\ \hline
    \multicolumn{2}{|l|}{$\rho_a^0$} & $6\cdot 10^{-7}~\mu\mbox{mol}/\mbox{dm}^2$  \\ \hline
    \multicolumn{2}{|l|}{$\rho_d^0$} & $1~\mu\mbox{mol}/\mbox{dm}^3$ \\ \hline
    \multicolumn{2}{|l|}{$D_1$} & $4~\mu\mbox{m}^2/$s  \\ \hline
    \multicolumn{2}{|l|}{$D_2$} & $4~\mu\mbox{m}^2/$s  \\ \hline
    \multicolumn{2}{|l|}{$D_3$} & $0.3~\mu\mbox{m}^2/$s  \\ \hline
    \multicolumn{2}{|l|}{$k_1$} & $0.035~\mbox{s}^{-1}$  \\ \hline
    \multicolumn{2}{|l|}{$k_2$} & $0.015~\mbox{s}^{-1}$  \\ \hline
    \multicolumn{2}{|l|}{$k_3$} & $0.379~\mbox{s}^{-1}$  \\ \hline
    \multicolumn{2}{|l|}{$k_4$} & $0.625~\mbox{s}^{-1}$  \\ \hline
    \multicolumn{2}{|l|}{$k_5$} & $0.0168~\mbox{s}^{-1}$  \\ \hline
    \multicolumn{2}{|l|}{$E$} & $0.1,5.7,7\cdot 10^6$~kPa  \\ \hline
    \multicolumn{2}{|l|}{$C$} & $3.25$~kPa  \\ \hline
    \multicolumn{2}{|l|}{$n$} & $5$  \\ \hline
    \multicolumn{2}{|l|}{$\gamma$} & $8.8068~\mbox{dm}^3/\mu\mbox{mol}$  \\ \hline
    \multirow{2}{*}{\makecell{axisymmetric \\ shape}} & $|Y|$ & $1193~\mu$m$^3$ \\ \cline{2-3} 
    & $|\Gamma|$ & $1020~\mu$m$^2$ \\ \hline
    \multirow{2}{*}{\makecell{lamellipodium \\ shape}} & $|Y|$ & $1099~\mu$m$^3$ \\ \cline{2-3}
    & $|\Gamma|$ & $1115~\mu$m$^2$ \\ \hline
  \end{tabular}
\end{subtable}\hspace{0.4cm}
\begin{subtable}{0.59\textwidth}
    \subcaption{Parameter values for parameters introduced in this paper. A brief robustness analysis is performed on all parameters introduced in the coupled model in Eqs~\eqref{eq:coupled_elast}-\eqref{eq:coupled_reactions}, see \nameref{S1_Appendix}, Section A.6.}
  \begin{tabular}{|l|l|l|}
  \hline
    Parameters & Value & Reference/ Justification \\ \hline
    $C_1$ & $0.1$~$(\mbox{kPa s})^{-1}$ & \makecell{range $0-2$~$(\mbox{kPa s})^{-1}$ \\ is explored in results} \\ \hline
    $k_6 $ & $0.1~\mbox{s}^{-1}$ & \makecell{fitted to yield magnitude of \\  deformation range $0-10~\mu$m \cite{scott_spatial_2021}} \\ \hline
    $k_7$  & $0.2$~kPa & \makecell{fitted to results \\  for $\phi_a$ in \cite{scott_spatial_2021}} \\ \hline
    $k_8$  & $2.4245~\mbox{dm}^3/\mu\mbox{mol}$ & \makecell{fitted to results \\  for $\phi_a$ in \cite{scott_spatial_2021}} \\ \hline
    $p$ & $2.6$ & \makecell{fitted to results for $\phi_a$ \\   and F-actin in \cite{scott_spatial_2021} and \cite{gardel_elastic_2004}} \\ \hline
    $\nu_c$ & $0.3$ & \makecell{estimated $0.17-0.66$  \cite{mokbelPoisson2020}} \\ \hline
  \end{tabular}
\end{subtable}
  \label{tab:parameters}
\end{table}

\newpage
\section*{Results}\label{sec:sims}
\subsection*{Numerical simulations with 2xD stimulus}\label{subsec:2xD}
\begin{figure}[!ht]
\begin{minipage}{\textwidth}
\begin{minipage}{0.05\textwidth}\centering
\vspace{-0.5cm}
    \textbf{(A)}\\ \vspace{0.3cm}
        $\rho_a$ \\ \vspace{0.65cm}
        $\phi_d$ \\ \vspace{0.65cm}
        $\phi_a$ \\ \vspace{0.65cm}
        $\vert u \vert$ \\ \vspace{0.65cm}
    \textbf{(C)}\\ \vspace{0.3cm}
        $\rho_a$ \\ \vspace{0.65cm}
        $\phi_d$ \\ \vspace{0.65cm}
        $\phi_a$ \\ \vspace{0.65cm}
        $\vert u \vert$ 
\end{minipage}\hfill
\begin{minipage}{0.45\textwidth}
\begin{minipage}{0.33\textwidth} \centering
    {\footnotesize $0.1$kPa } \\ 
        \includegraphics[width=\textwidth]{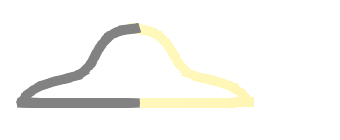}
        \includegraphics[width=\textwidth]{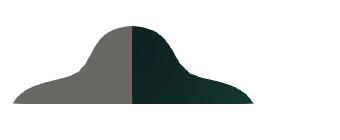}
        \includegraphics[width=\textwidth]{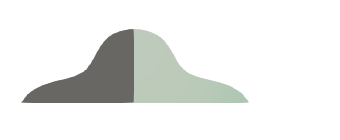}
        \includegraphics[width=\textwidth]{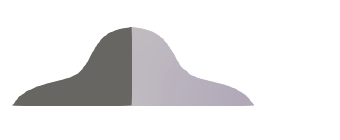}\\
        \vspace{0.5cm}
        \includegraphics[width=\textwidth]{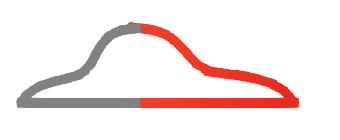}
        \includegraphics[width=\textwidth]{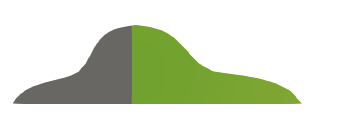}
        \includegraphics[width=\textwidth]{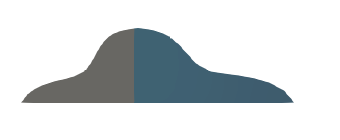}
        \includegraphics[width=\textwidth]{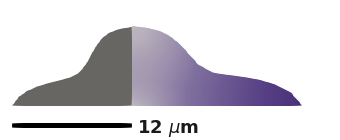}
\end{minipage}\hfill
\begin{minipage}{0.33\textwidth} \centering
    {\footnotesize $5.7$kPa } \\
        \includegraphics[width=\textwidth]{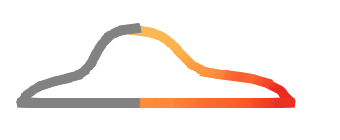}
        \includegraphics[width=\textwidth]{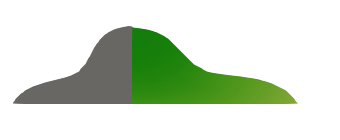}
        \includegraphics[width=\textwidth]{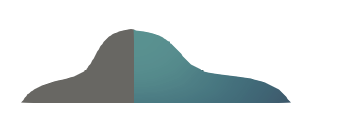}
        \includegraphics[width=\textwidth]{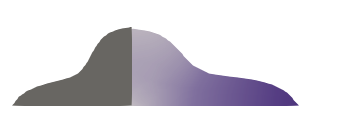}\\
        \vspace{0.5cm}
        \includegraphics[width=\textwidth]{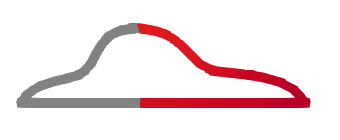}
        \includegraphics[width=\textwidth]{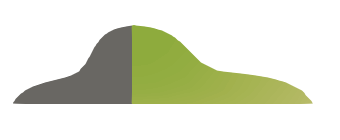}
        \includegraphics[width=\textwidth]{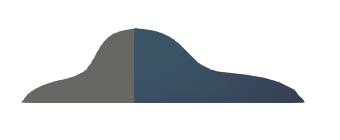}
        \includegraphics[width=\textwidth]{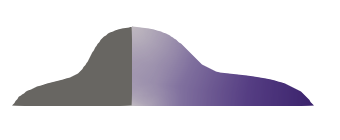}
\end{minipage}\hfill
\begin{minipage}{0.33\textwidth} \centering
    {\footnotesize $7$GPa} \\
        \includegraphics[width=\textwidth]{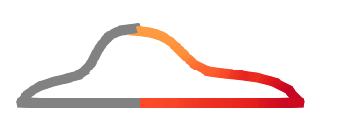}
        \includegraphics[width=\textwidth]{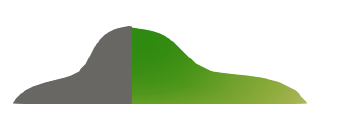}
        \includegraphics[width=\textwidth]{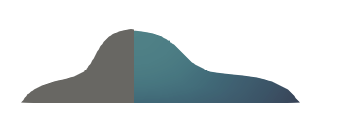}
        \includegraphics[width=\textwidth]{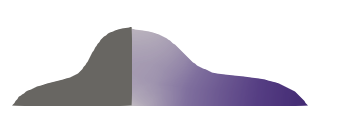}\\
        \vspace{0.5cm}
        \includegraphics[width=\textwidth]{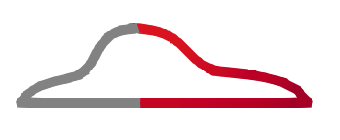}
        \includegraphics[width=\textwidth]{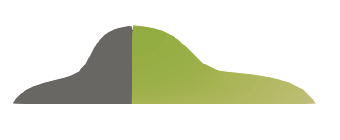}
        \includegraphics[width=\textwidth]{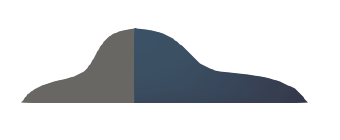}
        \includegraphics[width=\textwidth]{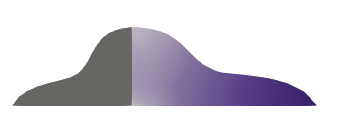}
\end{minipage}\hfill 
\end{minipage}\hfill 
\begin{minipage}{0.05\textwidth}\centering
\vspace{-0.5cm}
    \textbf{(B)}\\ \vspace{4.5cm}
    \textbf{(D)}\\ \vspace{3.8cm}
\end{minipage}\hfill
\begin{minipage}{0.45\textwidth}
\begin{minipage}{0.33\textwidth} \centering
    {\footnotesize $0.1$kPa } \\ 
        \includegraphics[width=\textwidth]{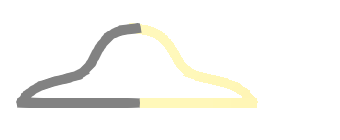}
        \includegraphics[width=\textwidth]{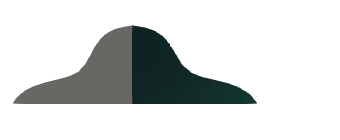}
        \includegraphics[width=\textwidth]{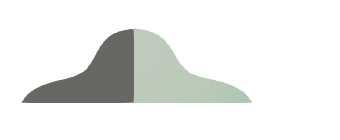}
        \includegraphics[width=\textwidth]{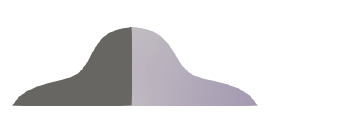}\\
        \vspace{0.5cm}
        \includegraphics[width=\textwidth]{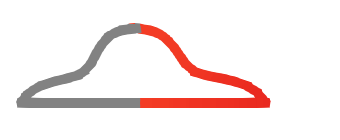}
        \includegraphics[width=\textwidth]{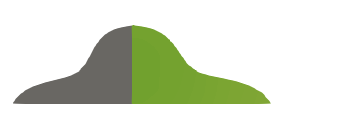}
        \includegraphics[width=\textwidth]{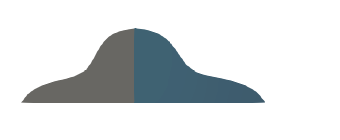}
        \includegraphics[width=\textwidth]{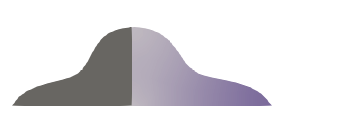}
\end{minipage}\hfill
\begin{minipage}{0.33\textwidth} \centering
    {\footnotesize $5.7$kPa } \\
        \includegraphics[width=\textwidth]{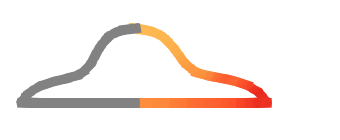}
        \includegraphics[width=\textwidth]{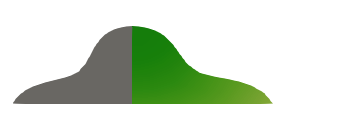}
        \includegraphics[width=\textwidth]{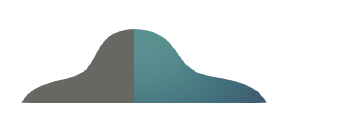}
        \includegraphics[width=\textwidth]{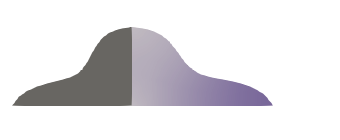}\\
        \vspace{0.5cm}
        \includegraphics[width=\textwidth]{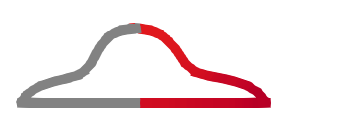}
        \includegraphics[width=\textwidth]{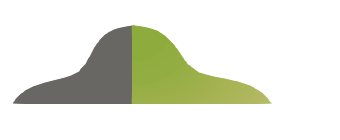}
        \includegraphics[width=\textwidth]{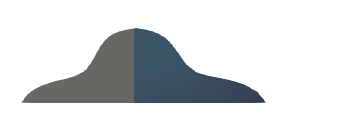}
        \includegraphics[width=\textwidth]{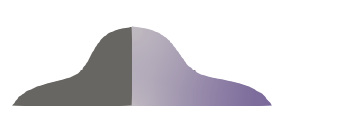}
\end{minipage}\hfill
\begin{minipage}{0.33\textwidth} \centering
    {\footnotesize $7$GPa} \\
        \includegraphics[width=\textwidth]{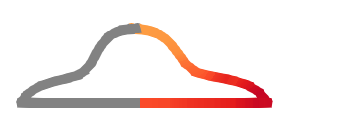}
        \includegraphics[width=\textwidth]{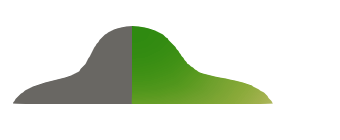}
        \includegraphics[width=\textwidth]{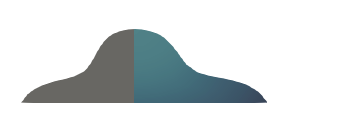}
        \includegraphics[width=\textwidth]{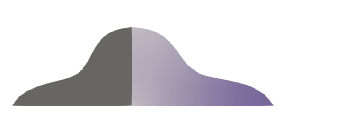}\\
        \vspace{0.5cm}
        \includegraphics[width=\textwidth]{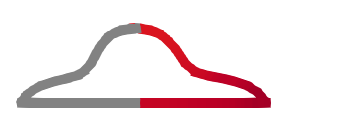}
        \includegraphics[width=\textwidth]{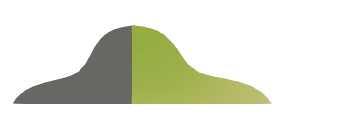}
        \includegraphics[width=\textwidth]{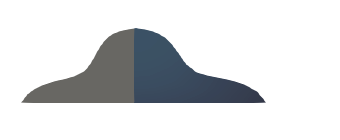}
        \includegraphics[width=\textwidth]{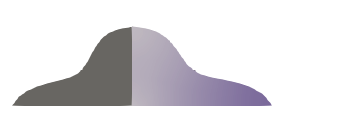}
\end{minipage}\hfill
\end{minipage}
\end{minipage}\hfill
\begin{minipage}{\textwidth}\raggedright
    \hspace{1cm}
    \includegraphics[width=0.22\textwidth]{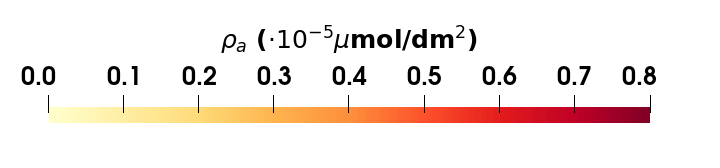}
    \includegraphics[width=0.22\textwidth]{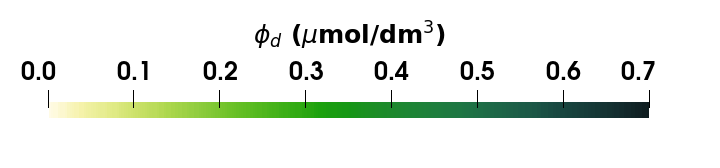}
    \includegraphics[width=0.22\textwidth]{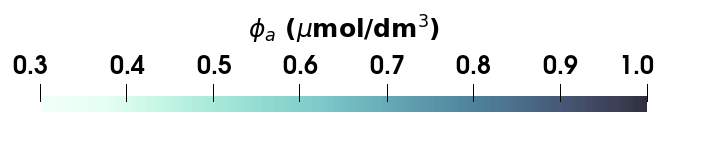}
    \includegraphics[width=0.22\textwidth]{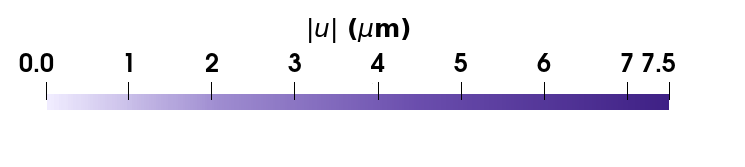}
\end{minipage}
  \caption{{\bf Numerical simulation results showing $\rho_a$, $\phi_d$, $\phi_a$ and $\vert u\vert$  for the model in Eqs~\eqref{eq:coupled_elast}-\eqref{eq:coupled_reactions} for the axisymmetric shape and in the case of $2$xD stimulus at a steady state at $T=100$~s.}
    Four different scenarios are considered: \textbf{(A)} $C_1=0$~$(\mbox{kPa s})^{-1}$ $(\sigma\not\rightarrow \phi_a)$ and $E_c=0.6$~kPa $(\phi_a\not\rightarrow E_c)$; \textbf{(B)} $C_1=0$~$(\mbox{kPa s})^{-1}$ $(\sigma\not\rightarrow \phi_a)$ and $E_c=f(\phi_a)$ $(\phi_a\rightarrow E_c)$; \textbf{(C)} $C_1=0.1$~$(\mbox{kPa s})^{-1}$ $(\sigma\rightarrow \phi_a)$ and $E_c=0.6$~kPa $(\phi_a\not\rightarrow E_c)$; \textbf{(D)} $C_1=0.1$~$(\mbox{kPa s})^{-1}$ $(\sigma\rightarrow \phi_a)$ and $E_c=f(\phi_a)$ $(\phi_a\rightarrow E_c)$. Within each subfigure, the rows represent $\rho_a$, $\phi_d$, $\phi_a$ and $\vert u\vert$ on a cross-section of the plane $x_1=0$ of the axisymmetric cell, and the columns represent $E=0.1, 5.7, 7\cdot 10^6$~kPa. Parameter values as in Table~\ref{tab:parameters}.}
  \label{fig:sim_subs_partfixed_2D}
\end{figure}

\begin{figure}[!ht]
\begin{minipage}{\textwidth}
\begin{minipage}{0.05\textwidth}\centering
\vspace{-0.5cm}
    \textbf{(A)}\\ \vspace{0.3cm}
        $\rho_a$ \\ \vspace{1.1cm}
        $\phi_d$ \\ \vspace{1.1cm}
        $\phi_a$ \\ \vspace{1.1cm}
        $\vert u \vert$ \\ \vspace{1.1cm}
    \textbf{(C)}\\ \vspace{0.3cm}
        $\rho_a$ \\ \vspace{1.1cm}
        $\phi_d$ \\ \vspace{1.1cm}
        $\phi_a$ \\ \vspace{1.1cm}
        $\vert u \vert$ 
\end{minipage}\hfill
\begin{minipage}{0.45\textwidth}
\begin{minipage}{0.33\textwidth} \centering
    {\footnotesize $0.1$kPa } \\ 
        \includegraphics[width=\textwidth]{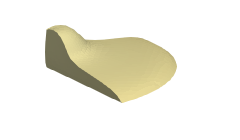}
        \includegraphics[width=\textwidth]{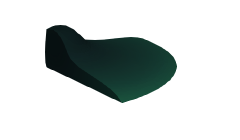}
        \includegraphics[width=\textwidth]{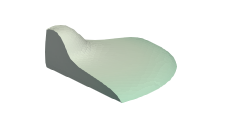}
        \includegraphics[width=\textwidth]{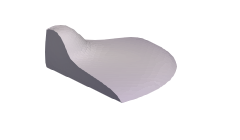}\\
        \vspace{0.5cm}
        \includegraphics[width=\textwidth]{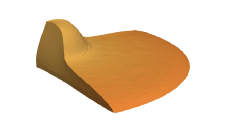}
        \includegraphics[width=\textwidth]{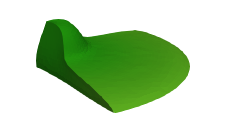}
        \includegraphics[width=\textwidth]{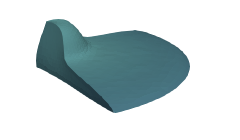}
        \includegraphics[width=\textwidth]{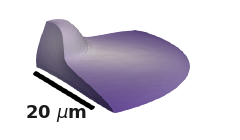}
\end{minipage}\hfill
\begin{minipage}{0.33\textwidth} \centering
    {\footnotesize $5.7$kPa } \\
        \includegraphics[width=\textwidth]{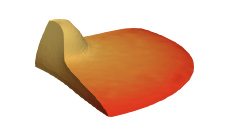}
        \includegraphics[width=\textwidth]{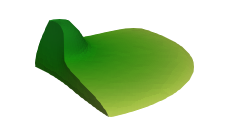}
        \includegraphics[width=\textwidth]{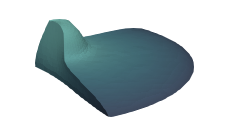}
        \includegraphics[width=\textwidth]{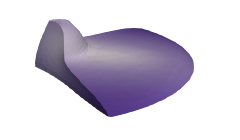}\\
        \vspace{0.5cm}
        \includegraphics[width=\textwidth]{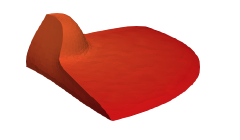}
        \includegraphics[width=\textwidth]{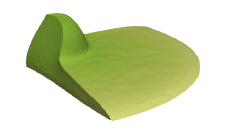}
        \includegraphics[width=\textwidth]{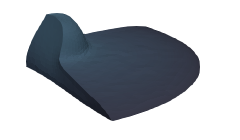}
        \includegraphics[width=\textwidth]{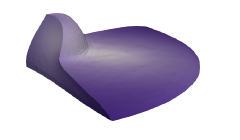}
\end{minipage}\hfill
\begin{minipage}{0.33\textwidth} \centering
    {\footnotesize $7$GPa} \\
        \includegraphics[width=\textwidth]{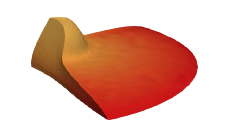}
        \includegraphics[width=\textwidth]{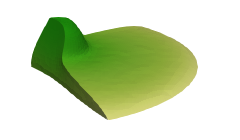}
        \includegraphics[width=\textwidth]{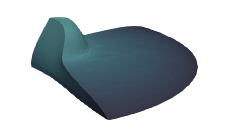}
        \includegraphics[width=\textwidth]{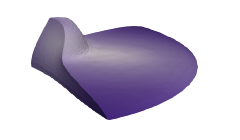}\\
        \vspace{0.5cm}
        \includegraphics[width=\textwidth]{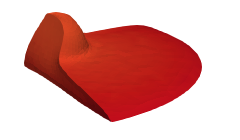}
        \includegraphics[width=\textwidth]{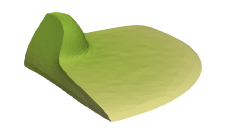}
        \includegraphics[width=\textwidth]{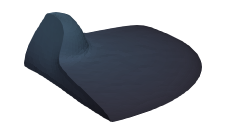}
        \includegraphics[width=\textwidth]{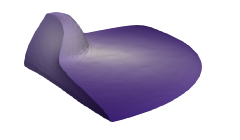}
\end{minipage}\hfill 
\end{minipage}\hfill 
\begin{minipage}{0.05\textwidth}\centering
\vspace{-0.5cm}
    \textbf{(B)}\\ \vspace{6.5cm}
    \textbf{(D)}\\ \vspace{5.3cm}
\end{minipage}\hfill
\begin{minipage}{0.45\textwidth}
\begin{minipage}{0.33\textwidth} \centering
    {\footnotesize $0.1$kPa } \\ 
        \includegraphics[width=\textwidth]{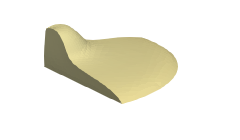}
        \includegraphics[width=\textwidth]{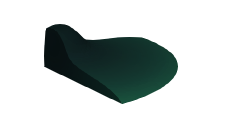}
        \includegraphics[width=\textwidth]{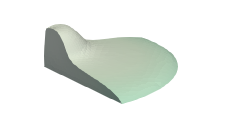}
        \includegraphics[width=\textwidth]{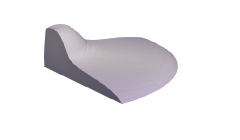}\\
        \vspace{0.5cm}
        \includegraphics[width=\textwidth]{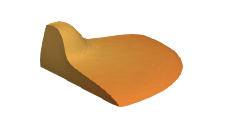}
        \includegraphics[width=\textwidth]{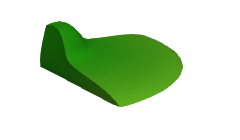}
        \includegraphics[width=\textwidth]{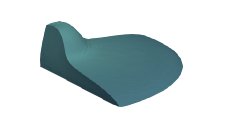}
        \includegraphics[width=\textwidth]{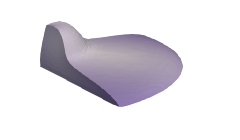}
\end{minipage}\hfill
\begin{minipage}{0.33\textwidth} \centering
    {\footnotesize $5.7$kPa } \\
        \includegraphics[width=\textwidth]{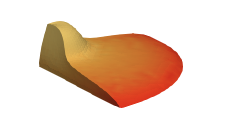}
        \includegraphics[width=\textwidth]{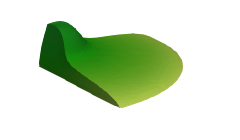}
        \includegraphics[width=\textwidth]{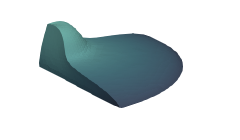}
        \includegraphics[width=\textwidth]{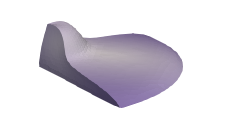}\\
        \vspace{0.5cm}
        \includegraphics[width=\textwidth]{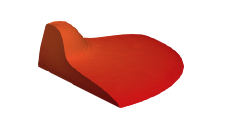}
        \includegraphics[width=\textwidth]{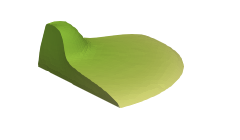}
        \includegraphics[width=\textwidth]{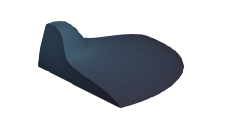}
        \includegraphics[width=\textwidth]{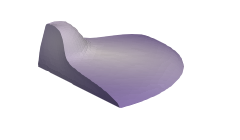}
\end{minipage}\hfill
\begin{minipage}{0.33\textwidth} \centering
    {\footnotesize $7$GPa} \\
        \includegraphics[width=\textwidth]{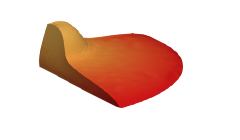}
        \includegraphics[width=\textwidth]{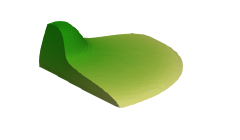}
        \includegraphics[width=\textwidth]{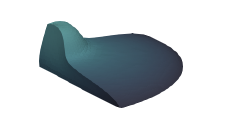}
        \includegraphics[width=\textwidth]{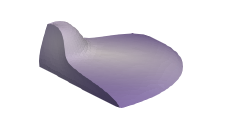}\\
        \vspace{0.5cm}
        \includegraphics[width=\textwidth]{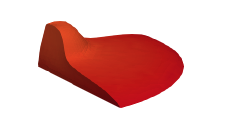}
        \includegraphics[width=\textwidth]{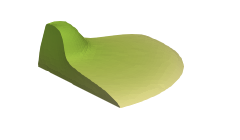}
        \includegraphics[width=\textwidth]{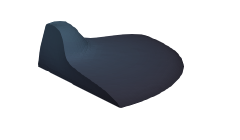}
        \includegraphics[width=\textwidth]{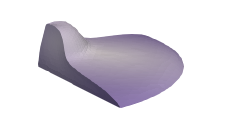}
\end{minipage}\hfill
\end{minipage}
\end{minipage}\hfill
\begin{minipage}{\textwidth}\raggedright
    \hspace{1cm}
    \includegraphics[width=0.22\textwidth]{figures/coupled_p_bar_grouped.png}
    \includegraphics[width=0.22\textwidth]{figures/coupled_cd_bar_grouped.png}
    \includegraphics[width=0.22\textwidth]{figures/coupled_ca_bar_grouped.png}
    \includegraphics[width=0.22\textwidth]{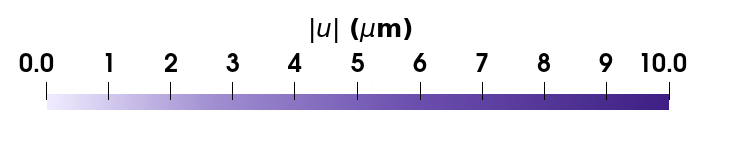}
\end{minipage}
  \caption{{\bf Numerical simulation results showing $\rho_a$, $\phi_d$, $\phi_a$ and $\vert u\vert$  for the model in Eqs~\eqref{eq:coupled_elast}-\eqref{eq:coupled_reactions} for the  lamellipodium shape and in the case of $2$xD stimulus at a steady state at $T=100$~s.}
    Four different scenarios are considered: \textbf{(A)} $C_1=0$~$(\mbox{kPa s})^{-1}$ $(\sigma\not\rightarrow \phi_a)$ and $E_c=0.6$~kPa $(\phi_a\not\rightarrow E_c)$; \textbf{(B)} $C_1=0$~$(\mbox{kPa s})^{-1}$ $(\sigma\not\rightarrow \phi_a)$ and $E_c=f(\phi_a)$ $(\phi_a\rightarrow E_c)$; \textbf{(C)} $C_1=0.1$~$(\mbox{kPa s})^{-1}$ $(\sigma\rightarrow \phi_a)$ and $E_c=0.6$~kPa $(\phi_a\not\rightarrow E_c)$; \textbf{(D)} $C_1=0.1$~$(\mbox{kPa s})^{-1}$ $(\sigma\rightarrow \phi_a)$ and $E_c=f(\phi_a)$ $(\phi_a\rightarrow E_c)$. Within each subfigure, the rows represent $\rho_a$, $\phi_d$, $\phi_a$ and $\vert u\vert$ on the surface of the cell, and the columns represent $E=0.1, 5.7, 7\cdot 10^6$~kPa. Parameter values as in Table~\ref{tab:parameters}.}
  \label{fig:sim_lamel_partfixed_2D}
\end{figure}

\begin{figure}[!ht]
  \begin{minipage}{\textwidth} 
    \hspace{3cm} axisymmetric shape \hspace{4.5cm} lamellipodium shape \\ 
    \vspace{-0.2cm}
    \hspace{1.8cm} {\footnotesize $E_c=0.6$~kPa } \hspace{2cm}  {\footnotesize $E_c=f(\phi_a)$ } \hspace{2cm} {\footnotesize $E_c=0.6$~kPa } \hspace{2cm} {\footnotesize $E_c=f(\phi_a)$ } \\ 
    \vspace{-0.2cm}
\end{minipage}\hfill
\begin{minipage}{0.95\textwidth} 
    \includegraphics[width=\textwidth]{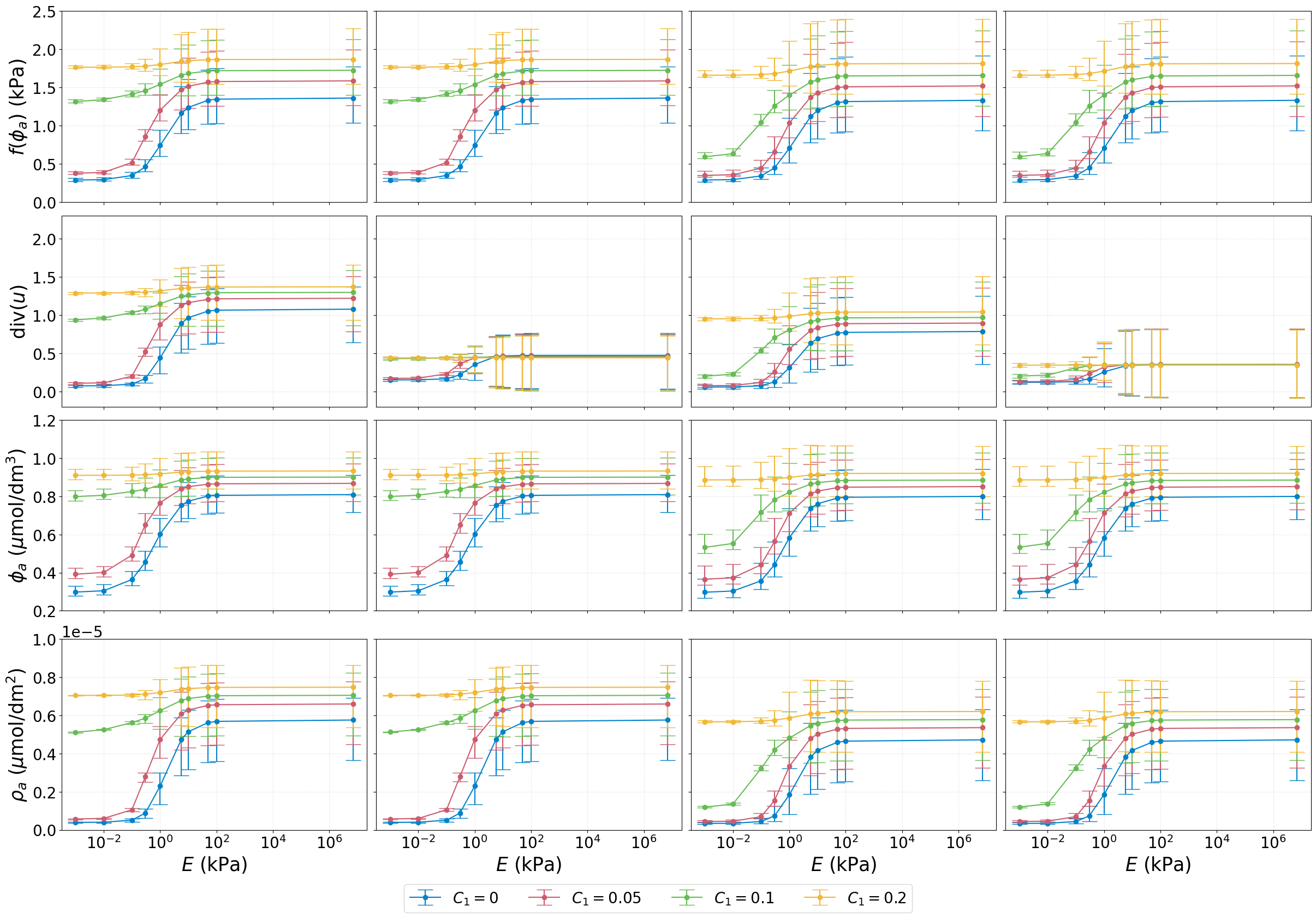}
\end{minipage}
  \caption{{\bf Simulation results showing the mean, $\frac{1}{|\Omega|}\int_\Omega \cdot \dd{x}$, min and max values of $f(\phi_a)$, ${\rm div}(u)$, $\phi_a$ and $\rho_a$ as  functions of substrate stiffness $E$.}
    We consider different couplings  with four different values for $C_1$ and two different shapes at $T=100$~s by which time the results are at a steady state. All other parameter values as in Table~\ref{tab:parameters}.}
  \label{fig:comp_partfixed_2D}
\end{figure}

First we look at the results that would most reflect a cell on a substrate in vitro. Here, the substrate stiffness appears as a stimulus only on the bottom boundary of the cell, i.e.~$E$ is nonzero only on $\Gamma_0$, and deformation is restricted in the vertical direction at the bottom boundary of the cell. The results for the axisymmetric shape of the cell are found in Fig~\ref{fig:sim_subs_partfixed_2D}, whereas results for the lamellipodium shape are presented in Fig~\ref{fig:sim_lamel_partfixed_2D}. Note that results for $\phi_a$ and $\rho_a$ when $E_c=0.6$~kPa and $C_1=0$~$(\mbox{kPa s})^{-1}$ are identical to the one without mechanics in \nameref{S1_Appendix}, Section A.1. In this case, we see that the magnitude of the deformation $|u|$ is largest at the edge of the cell. The cell expands axisymmetrically at the base. As expected, the expansion is larger for higher concentrations of $\phi_a$. For a lower substrate stiffness, $E=0.1$~kPa, the cell barely expands. When $C_1=0.1$~$(\mbox{kPa s})^{-1}$, the concentrations of $\phi_a$ and $\rho_a$ and the magnitude of the deformation $|u|$ increase, with a bigger increase for lower substrate stiffness and a smaller increase for larger substrate stiffness. 
When comparing $E_c=0.6$~kPa and $E_c=f(\phi_a)$,  the deformations show similar patterns, expanding at the base of the cell,  however, the magnitude of the deformation is much lower in the case $E_c=f(\phi_a)$. This is probably because $E_c=f(\phi_a)\approx 0.6$~kPa for a small substrate stiffness $E$, but is doubled in magnitude for larger substrate stiffness, see Fig~\ref{fig:comp_partfixed_2D}. The larger cell Young's modulus $E_c$ means it is harder for the cell to deform, resulting in a lower magnitude of deformation. This difference illustrates that, unlike the constant Young's modulus case,  a concentration-dependent Young's modulus allows for potential homeostasis and adaptation of cell mechanics to different values of the substrate stiffness~\cite{Grolleman_2023}.

For the two-way couplings between the mechanics and chemistry, i.e.~$E_c=f(\phi_a)$ and $C_1=0.1$~$(\mbox{kPa s})^{-1}$, we see similar results for the deformation as when $E_c=f(\phi_a)$ and $C_1=0$~$(\mbox{kPa s})^{-1}$. The main difference is that the deformation for $E=0.1$~kPa is now at a similar magnitude as for the larger substrate stiffnesses, demonstrating the importance of the signalling processes in the adaptation of cell mechanics to changing environmental conditions. 

Comparing the simulation results for the two different shapes in Fig~\ref{fig:sim_subs_partfixed_2D} and~\ref{fig:sim_lamel_partfixed_2D}, the concentration of activated RhoA, $\rho_a$, is slightly lower for the lamellipodium shape. For the lamellipodium shape, we observe the largest deformations at the corners furthest from the nucleus. 

Fig~\ref{fig:comp_partfixed_2D} summarises the results at time $T=100$~s by plotting the mean, $\tfrac{1}{|\Omega|}\int_\Omega \cdot \dd{x}$, of $E_c=f(\phi_a)$, the volume change ${\rm div}(u)$, $\phi_a$, and $\rho_a$ as functions of the substrate stiffness $E$, with the bars being the range of these variables, for different values of the constant $C_1$ in the activation of FAK by the cell stress.  As expected, an increase in $C_1$ results in an increase in the concentration of activated FAK, $\phi_a$. 
The dependence of $\phi_a$ on the substrate stiffness $E$, especially for $C_1=0$~$(\mbox{kPa s})^{-1}$ and $C_1=0.5$~$(\mbox{kPa s})^{-1}$,  resembles a Hill function  representing a threshold response.  This agrees with simulations in~\cite{scott_spatial_2021} which themselves fit experimental observations presented in~\cite{beamish_engineered_2017}. 
For most of the cases, the results for the lamellipodium shape are very similar to the results for the axisymmetric shape. However, for $C_1=0.1$~$(\mbox{kPa s})^{-1}$, the magnitude of the threshold-like response in all variables is bigger in the lamellipodium case. In terms of the Young's modulus, when $E_c=0.6$~kPa we  observe much larger volume changes than when $E_c=f(\phi_a)$ in all the numerical experiments.

\clearpage
\subsection*{Numerical simulations for the 3D stimulus case on a rigid substrate}\label{3Dstimulus_1}
\begin{figure}[!ht]
\begin{minipage}{\textwidth}
\begin{minipage}{0.05\textwidth}\centering
\vspace{-0.5cm}
    \textbf{(A)}\\ \vspace{0.3cm}
        $\rho_a$ \\ \vspace{0.65cm}
        $\phi_d$ \\ \vspace{0.65cm}
        $\phi_a$ \\ \vspace{0.65cm}
        $\vert u \vert$ \\ \vspace{0.65cm}
    \textbf{(C)}\\ \vspace{0.3cm}
        $\rho_a$ \\ \vspace{0.65cm}
        $\phi_d$ \\ \vspace{0.65cm}
        $\phi_a$ \\ \vspace{0.65cm}
        $\vert u \vert$ 
\end{minipage}\hfill
\begin{minipage}{0.45\textwidth}
\begin{minipage}{0.33\textwidth} \centering
    {\footnotesize $0.1$kPa } \\ 
        \includegraphics[width=\textwidth]{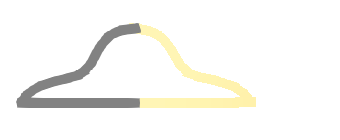}
        \includegraphics[width=\textwidth]{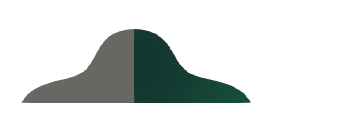}
        \includegraphics[width=\textwidth]{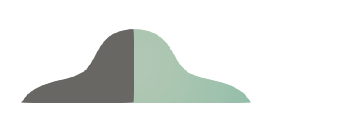}
        \includegraphics[width=\textwidth]{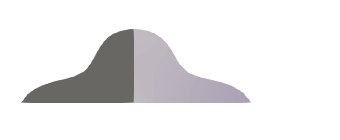}\\
        \vspace{0.5cm}
        \includegraphics[width=\textwidth]{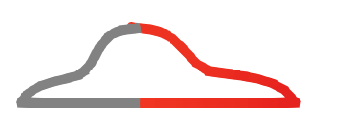}
        \includegraphics[width=\textwidth]{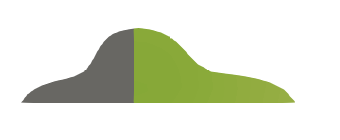}
        \includegraphics[width=\textwidth]{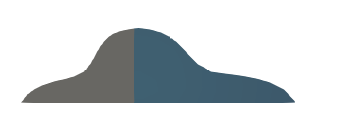}
        \includegraphics[width=\textwidth]{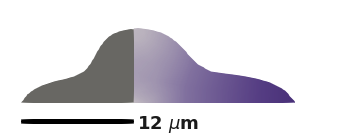}
\end{minipage}\hfill
\begin{minipage}{0.33\textwidth} \centering
    {\footnotesize $5.7$kPa } \\
        \includegraphics[width=\textwidth]{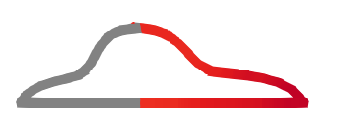}
        \includegraphics[width=\textwidth]{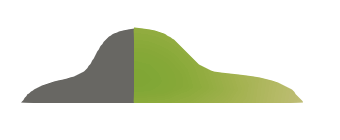}
        \includegraphics[width=\textwidth]{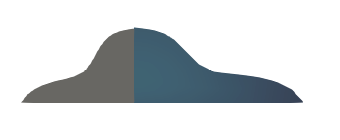}
        \includegraphics[width=\textwidth]{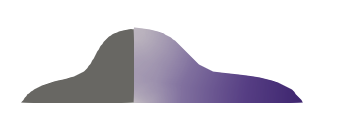}\\
        \vspace{0.5cm}
        \includegraphics[width=\textwidth]{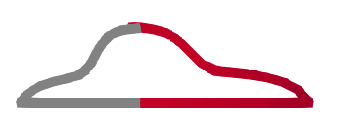}
        \includegraphics[width=\textwidth]{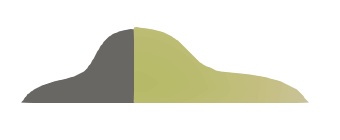}
        \includegraphics[width=\textwidth]{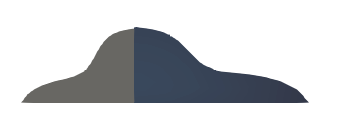}
        \includegraphics[width=\textwidth]{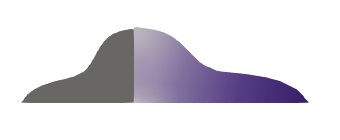}
\end{minipage}\hfill
\begin{minipage}{0.33\textwidth} \centering
    {\footnotesize $7$GPa} \\
        \includegraphics[width=\textwidth]{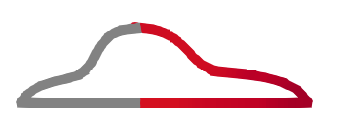}
        \includegraphics[width=\textwidth]{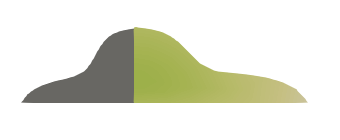}
        \includegraphics[width=\textwidth]{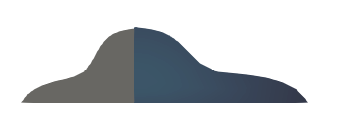}
        \includegraphics[width=\textwidth]{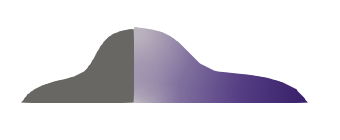}\\
        \vspace{0.5cm}
        \includegraphics[width=\textwidth]{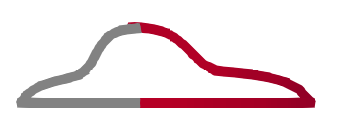}
        \includegraphics[width=\textwidth]{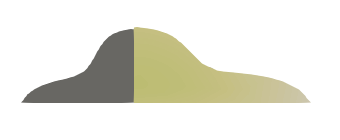}
        \includegraphics[width=\textwidth]{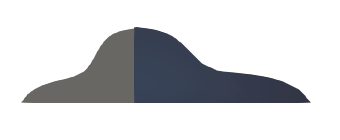}
        \includegraphics[width=\textwidth]{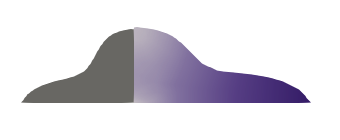}
\end{minipage}\hfill 
\end{minipage}\hfill 
\begin{minipage}{0.05\textwidth}\centering
\vspace{-0.5cm}
    \textbf{(B)}\\ \vspace{4.5cm}
    \textbf{(D)}\\ \vspace{3.8cm}
\end{minipage}\hfill
\begin{minipage}{0.45\textwidth}
\begin{minipage}{0.33\textwidth} \centering
    {\footnotesize $0.1$kPa } \\ 
        \includegraphics[width=\textwidth]{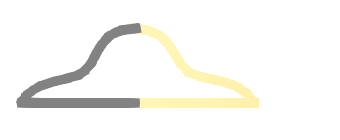}
        \includegraphics[width=\textwidth]{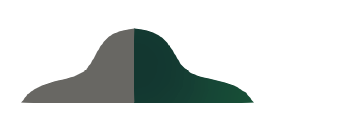}
        \includegraphics[width=\textwidth]{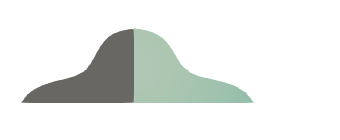}
        \includegraphics[width=\textwidth]{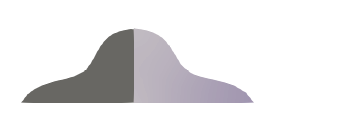}\\
        \vspace{0.5cm}
        \includegraphics[width=\textwidth]{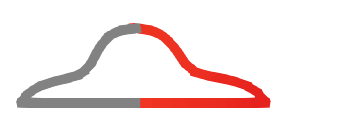}
        \includegraphics[width=\textwidth]{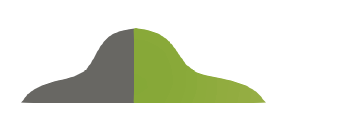}
        \includegraphics[width=\textwidth]{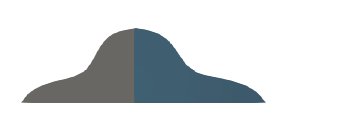}
        \includegraphics[width=\textwidth]{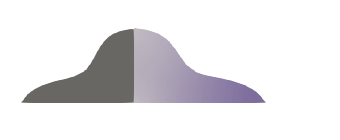}
\end{minipage}\hfill
\begin{minipage}{0.33\textwidth} \centering
    {\footnotesize $5.7$kPa } \\
        \includegraphics[width=\textwidth]{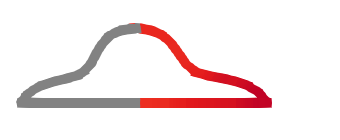}
        \includegraphics[width=\textwidth]{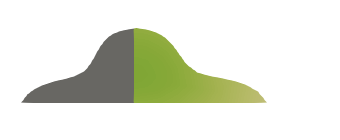}
        \includegraphics[width=\textwidth]{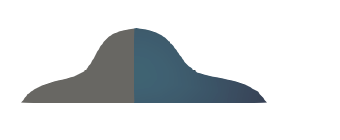}
        \includegraphics[width=\textwidth]{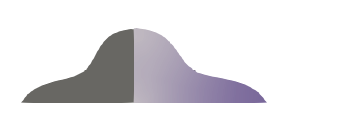}\\
        \vspace{0.5cm}
        \includegraphics[width=\textwidth]{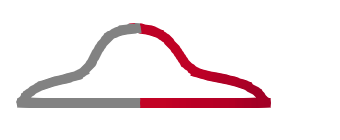}
        \includegraphics[width=\textwidth]{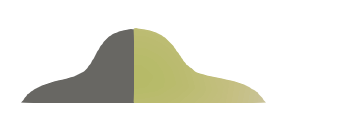}
        \includegraphics[width=\textwidth]{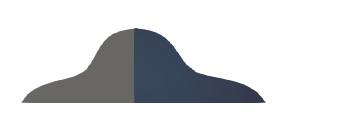}
        \includegraphics[width=\textwidth]{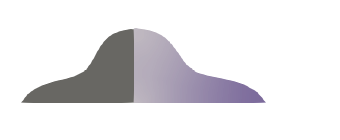}
\end{minipage}\hfill
\begin{minipage}{0.33\textwidth} \centering
    {\footnotesize $7$GPa} \\
        \includegraphics[width=\textwidth]{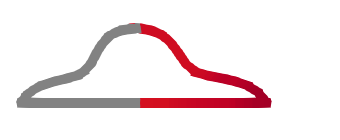}
        \includegraphics[width=\textwidth]{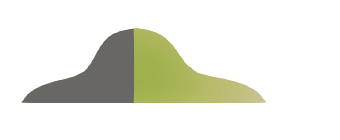}
        \includegraphics[width=\textwidth]{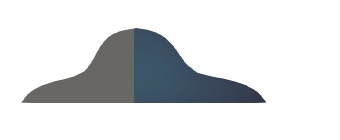}
        \includegraphics[width=\textwidth]{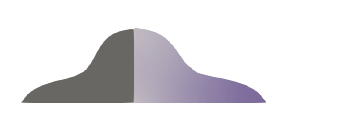}\\
        \vspace{0.5cm}
        \includegraphics[width=\textwidth]{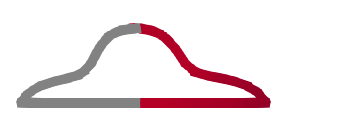}
        \includegraphics[width=\textwidth]{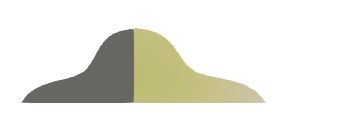}
        \includegraphics[width=\textwidth]{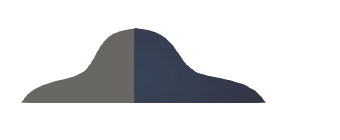}
        \includegraphics[width=\textwidth]{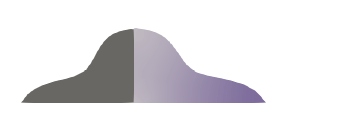}
\end{minipage}\hfill
\end{minipage}
\end{minipage}\hfill
\begin{minipage}{\textwidth}\raggedright
    \hspace{1cm}
    \includegraphics[width=0.22\textwidth]{figures/coupled_p_bar_grouped.png}
    \includegraphics[width=0.22\textwidth]{figures/coupled_cd_bar_grouped.png}
    \includegraphics[width=0.22\textwidth]{figures/coupled_ca_bar_grouped.png}
    \includegraphics[width=0.22\textwidth]{figures/coupled_u_bar_grouped.png}
\end{minipage}
  \caption{{\bf Numerical simulation results showing $\rho_a$, $\phi_d$, $\phi_a$ and $\vert u\vert$  for the model in Eqs~\eqref{eq:coupled_elast}-\eqref{eq:coupled_reactions}  for the axisymmetric shape  and in the case of the $3$D stimulus at a steady state at $T=100$~s.}
    Four different scenarios are considered: \textbf{(A)} $C_1=0$~$(\mbox{kPa s})^{-1}$ $(\sigma\not\rightarrow \phi_a)$ and $E_c=0.6$~kPa $(\phi_a\not\rightarrow E_c)$; \textbf{(B)} $C_1=0$~$(\mbox{kPa s})^{-1}$ $(\sigma\not\rightarrow \phi_a)$ and $E_c=f(\phi_a)$ $(\phi_a\rightarrow E_c)$; \textbf{(C)} $C_1=0.1$~$(\mbox{kPa s})^{-1}$ $(\sigma\rightarrow \phi_a)$ and $E_c=0.6$~kPa $(\phi_a\not\rightarrow E_c)$; \textbf{(D)} $C_1=0.1$~$(\mbox{kPa s})^{-1}$ $(\sigma\rightarrow \phi_a)$ and $E_c=f(\phi_a)$ $(\phi_a\rightarrow E_c)$. Within each subfigure, the rows represent $\rho_a$, $\phi_d$, $\phi_a$ and $\vert u\vert$ on a cross-section of the plane $x_1=0$ of the axisymmetric cell, and the columns represent $E=0.1, 5.7, 7\cdot 10^6$~kPa. Parameter values as in Table~\ref{tab:parameters}.}
  \label{fig:sim_subs_partfixed}
\end{figure}

\begin{figure}[!h]
  \begin{minipage}{\textwidth} 
    \hspace{3cm} axisymmetric shape \hspace{4.5cm} lamellipodium shape \\ 
    \vspace{-0.2cm}
    \hspace{1.8cm} {\footnotesize $E_c=0.6$~kPa } \hspace{2cm}  {\footnotesize $E_c=f(\phi_a)$ } \hspace{2cm} {\footnotesize $E_c=0.6$~kPa } \hspace{2cm} {\footnotesize $E_c=f(\phi_a)$ } \\ 
    \vspace{-0.2cm}
\end{minipage}\hfill
\begin{minipage}{0.95\textwidth} 
    \includegraphics[width=\textwidth]{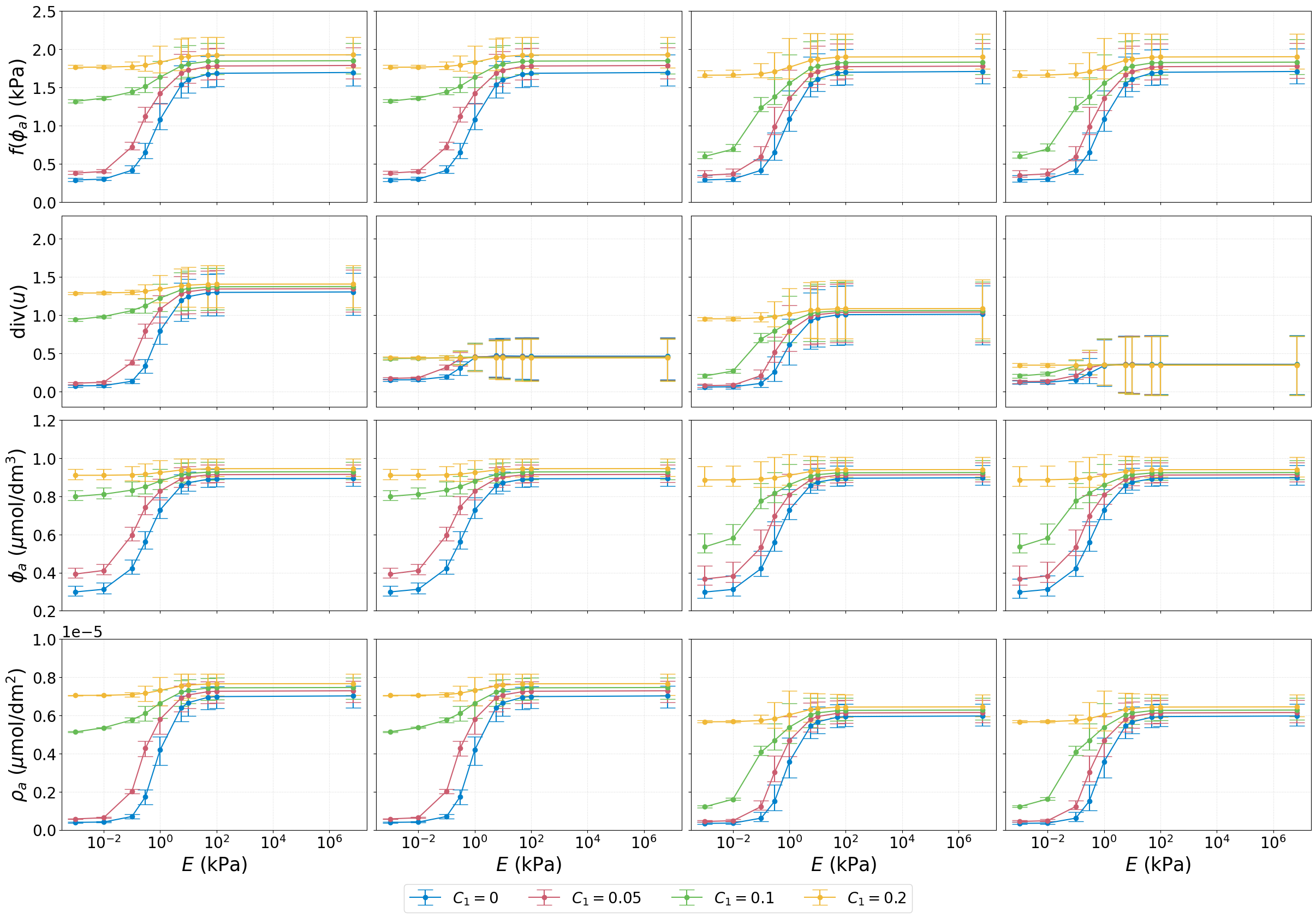}
\end{minipage}
  \caption{{\bf Simulation results showing the mean, $\frac{1}{|\Omega|}\int_\Omega \cdot \dd{x}$, min and max values of $f(\phi_a)$, ${\rm div}(u)$, $\phi_a$ and $\rho_a$ as  functions of substrate stiffness $E$.}
    We consider different couplings  with four different values for $C_1$ and two different shapes at $T=100$~s by which time the results are at a steady state. All other parameter values as in Table~\ref{tab:parameters}.}
  \label{fig:comp_partfixed}
\end{figure}

In numerical simulations for a $3$D stimulus on a rigid substrate, the substrate stiffness affects the whole cell membrane and we consider the boundary conditions Eq~(\ref{eq:coupled_force_boundary}) on $\Gamma \setminus \Gamma_0$ and Eq~(\ref{eq:coupled_normal_zero}) on $\Gamma_0$. 
The results for  numerical experiments can be found in Fig~\ref{fig:sim_subs_partfixed}-\ref{fig:comp_partfixed}. 
Overall,  the concentrations $\phi_a$ and $\rho_a$ are larger than in the case of the $2$xD stimulus, which is in line with the results in~\cite{scott_spatial_2021}. The higher concentrations of $\rho_a$ results in  larger deformations, where the maximum magnitude of the deformation 
in the case of the $2$xD stimulus was $7~\mu$m, see Fig~\ref{fig:sim_subs_partfixed_2D}, while the maximum magnitude of the deformation in the case of the $3$D stimulus is $7.5~\mu$m, see Fig~\ref{fig:sim_subs_partfixed}. Similar behaviour is observed for the lamellipodium shape, see Fig~\ref{fig:sim_lamel_partfixed_2D} and in \nameref{S1_Fig}. 
Another difference between two cases are larger variations in concentration and a larger difference between maximal and minimal values in the case of the $2$xD stimulus than in the case of $3$D stimulus, see Fig~\ref{fig:comp_partfixed_2D} and~\ref{fig:comp_partfixed}. Similar behaviour is observed also in the model for the signalling processes without mechanics, see \nameref{S1_Appendix}, Section A.1. 

\newpage
\subsection*{Numerical simulations for the model in~Eqs\eqref{eq:coupled_elast},~\eqref{eq:coupled_force_boundary}, and \eqref{eq:coupled_reactions}.}
 To investigate a setting more close to a cell in vivo, we consider the coupled model in Eqs~\eqref{eq:coupled_elast}-\eqref{eq:coupled_force_boundary}, \eqref{eq:coupled_reactions} with force boundary conditions on the whole cell membrane, without restricting the deformation on the bottom of the cell. 
 
\subsubsection*{Numerical simulations  in the case of $3$D stimulus.}
\begin{figure}[!ht]
\begin{minipage}{\textwidth}
\begin{minipage}{0.05\textwidth}\centering
\vspace{-0.5cm}
    \textbf{(A)}\\ \vspace{0.3cm}
        $\rho_a$ \\ \vspace{0.65cm}
        $\phi_d$ \\ \vspace{0.65cm}
        $\phi_a$ \\ \vspace{0.65cm}
        $\vert u \vert$ \\ \vspace{0.65cm}
    \textbf{(C)}\\ \vspace{0.3cm}
        $\rho_a$ \\ \vspace{0.65cm}
        $\phi_d$ \\ \vspace{0.65cm}
        $\phi_a$ \\ \vspace{0.65cm}
        $\vert u \vert$ 
\end{minipage}\hfill
\begin{minipage}{0.45\textwidth}
\begin{minipage}{0.33\textwidth} \centering
    {\footnotesize $0.1$kPa } \\ 
        \includegraphics[width=\textwidth]{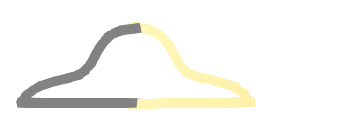}
        \includegraphics[width=\textwidth]{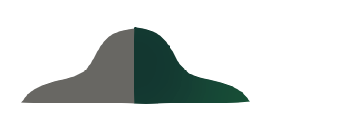}
        \includegraphics[width=\textwidth]{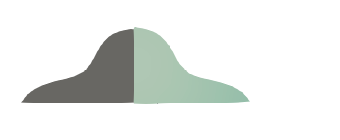}
        \includegraphics[width=\textwidth]{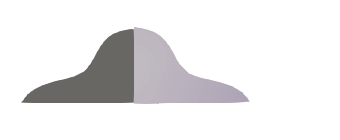}\\
        \vspace{0.5cm}
        \includegraphics[width=\textwidth]{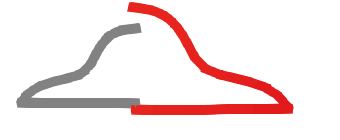}
        \includegraphics[width=\textwidth]{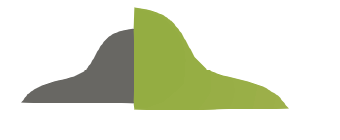}
        \includegraphics[width=\textwidth]{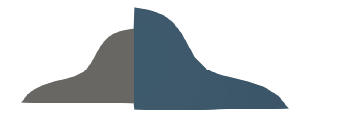}
        \includegraphics[width=\textwidth]{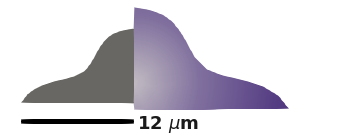}
\end{minipage}\hfill
\begin{minipage}{0.33\textwidth} \centering
    {\footnotesize $5.7$kPa } \\
        \includegraphics[width=\textwidth]{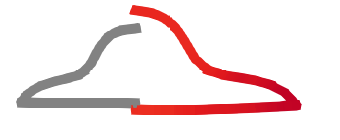}
        \includegraphics[width=\textwidth]{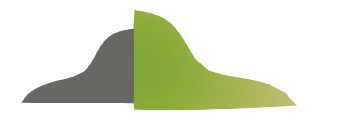}
        \includegraphics[width=\textwidth]{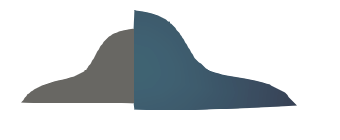}
        \includegraphics[width=\textwidth]{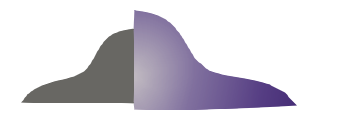}\\
        \vspace{0.5cm}
        \includegraphics[width=\textwidth]{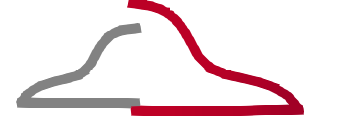}
        \includegraphics[width=\textwidth]{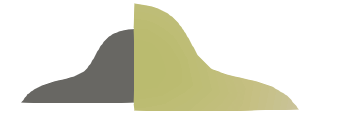}
        \includegraphics[width=\textwidth]{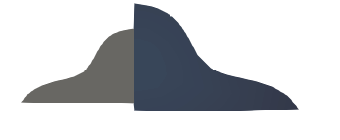}
        \includegraphics[width=\textwidth]{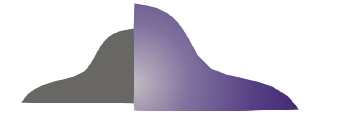}
\end{minipage}\hfill
\begin{minipage}{0.33\textwidth} \centering
    {\footnotesize $7$GPa} \\
        \includegraphics[width=\textwidth]{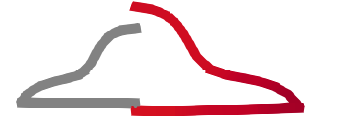}
        \includegraphics[width=\textwidth]{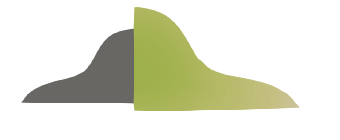}
        \includegraphics[width=\textwidth]{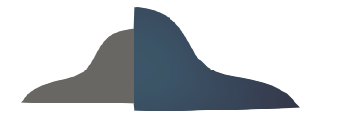}
        \includegraphics[width=\textwidth]{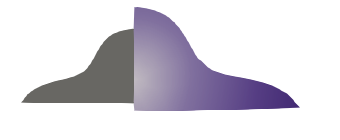}\\
        \vspace{0.5cm}
        \includegraphics[width=\textwidth]{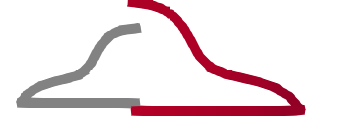}
        \includegraphics[width=\textwidth]{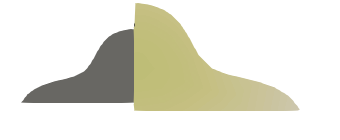}
        \includegraphics[width=\textwidth]{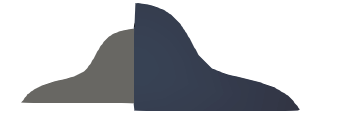}
        \includegraphics[width=\textwidth]{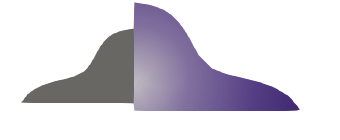}
\end{minipage}\hfill 
\end{minipage}\hfill 
\begin{minipage}{0.05\textwidth}\centering
\vspace{-0.5cm}
    \textbf{(B)}\\ \vspace{4.5cm}
    \textbf{(D)}\\ \vspace{3.8cm}
\end{minipage}\hfill
\begin{minipage}{0.45\textwidth}
\begin{minipage}{0.33\textwidth} \centering
    {\footnotesize $0.1$kPa } \\ 
        \includegraphics[width=\textwidth]{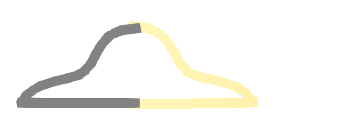}
        \includegraphics[width=\textwidth]{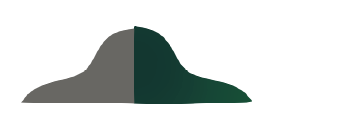}
        \includegraphics[width=\textwidth]{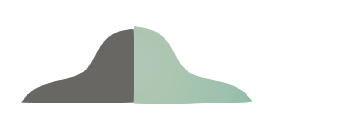}
        \includegraphics[width=\textwidth]{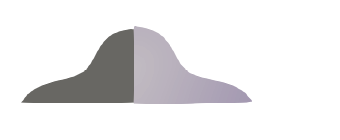}\\
        \vspace{0.5cm}
        \includegraphics[width=\textwidth]{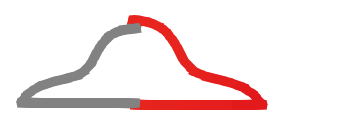}
        \includegraphics[width=\textwidth]{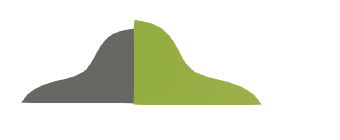}
        \includegraphics[width=\textwidth]{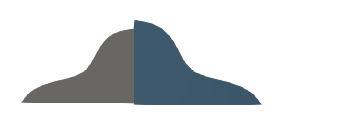}
        \includegraphics[width=\textwidth]{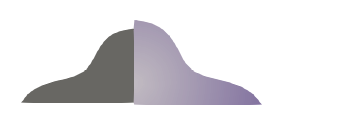}
\end{minipage}\hfill
\begin{minipage}{0.33\textwidth} \centering
    {\footnotesize $5.7$kPa } \\
        \includegraphics[width=\textwidth]{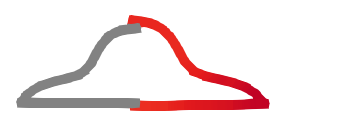}
        \includegraphics[width=\textwidth]{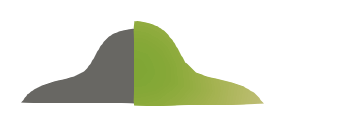}
        \includegraphics[width=\textwidth]{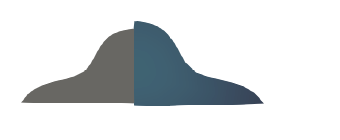}
        \includegraphics[width=\textwidth]{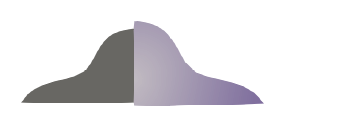}\\
        \vspace{0.5cm}
        \includegraphics[width=\textwidth]{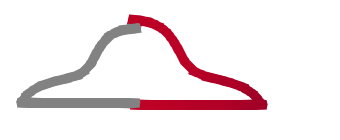}
        \includegraphics[width=\textwidth]{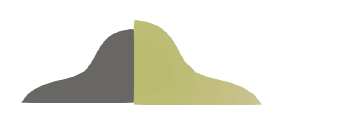}
        \includegraphics[width=\textwidth]{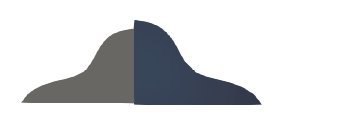}
        \includegraphics[width=\textwidth]{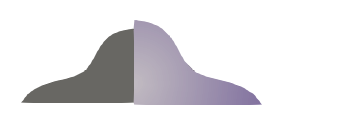}
\end{minipage}\hfill
\begin{minipage}{0.33\textwidth} \centering
    {\footnotesize $7$GPa} \\
        \includegraphics[width=\textwidth]{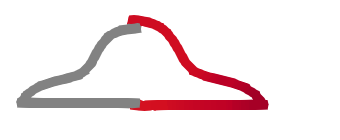}
        \includegraphics[width=\textwidth]{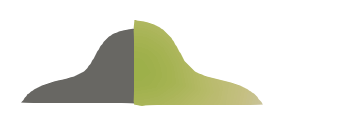}
        \includegraphics[width=\textwidth]{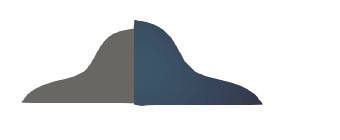}
        \includegraphics[width=\textwidth]{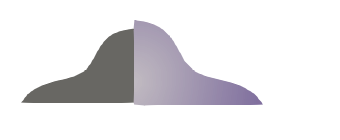}\\
        \vspace{0.5cm}
        \includegraphics[width=\textwidth]{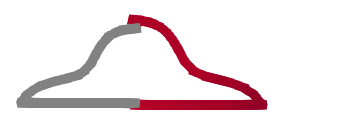}
        \includegraphics[width=\textwidth]{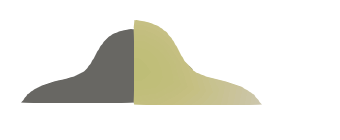}
        \includegraphics[width=\textwidth]{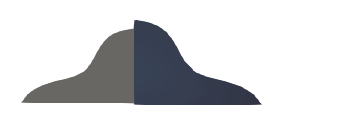}
        \includegraphics[width=\textwidth]{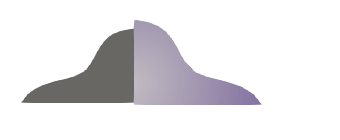}
\end{minipage}\hfill
\end{minipage}
\end{minipage}\hfill
\begin{minipage}{\textwidth}\raggedright
    \hspace{1cm}
    \includegraphics[width=0.22\textwidth]{figures/coupled_p_bar_grouped.png}
    \includegraphics[width=0.22\textwidth]{figures/coupled_cd_bar_grouped.png}
    \includegraphics[width=0.22\textwidth]{figures/coupled_ca_bar_grouped.png}
    \includegraphics[width=0.22\textwidth]{figures/coupled_u_bar_grouped.png}
\end{minipage}
  \caption{{\bf Numerical simulation results showing $\rho_a$, $\phi_d$, $\phi_a$ and $\vert u\vert$  for the model in Eqs~\eqref{eq:coupled_elast}, \eqref{eq:coupled_force_boundary}, and \eqref{eq:coupled_reactions}  for the axisymmetric shape  and in the case of the $3$D stimulus at a steady state at $T=100$~s.}
    Four different scenarios are considered: \textbf{(A)} $C_1=0$~$(\mbox{kPa s})^{-1}$ $(\sigma\not\rightarrow \phi_a)$ and $E_c=0.6$~kPa $(\phi_a\not\rightarrow E_c)$; \textbf{(B)} $C_1=0$~$(\mbox{kPa s})^{-1}$ $(\sigma\not\rightarrow \phi_a)$ and $E_c=f(\phi_a)$ $(\phi_a\rightarrow E_c)$; \textbf{(C)} $C_1=0.1$~$(\mbox{kPa s})^{-1}$ $(\sigma\rightarrow \phi_a)$ and $E_c=0.6$~kPa $(\phi_a\not\rightarrow E_c)$; \textbf{(D)} $C_1=0.1$~$(\mbox{kPa s})^{-1}$ $(\sigma\rightarrow \phi_a)$ and $E_c=f(\phi_a)$ $(\phi_a\rightarrow E_c)$. Within each subfigure, the rows represent $\rho_a$, $\phi_d$, $\phi_a$ and $\vert u\vert$ on a cross-section of the plane $x_1=0$ of the axisymmetric cell, and the columns represent $E=0.1, 5.7, 7\cdot 10^6$~kPa. Parameter values as in Table~\ref{tab:parameters}.}
  \label{fig:sim_subs}
\end{figure}

\begin{figure}[!h]
  \begin{minipage}{\textwidth}
    \hspace{3cm} axisymmetric shape \hspace{4.5cm} lamellipodium shape \\ 
    \vspace{-0.2cm}
    \hspace{1.8cm} {\footnotesize $E_c=0.6$~kPa } \hspace{2cm}  {\footnotesize $E_c=f(\phi_a)$ } \hspace{2cm} {\footnotesize $E_c=0.6$~kPa } \hspace{2cm} {\footnotesize $E_c=f(\phi_a)$ } \\ 
    \vspace{-0.2cm}
\end{minipage}\hfill
\begin{minipage}{0.95\textwidth} 
    \includegraphics[width=\textwidth]{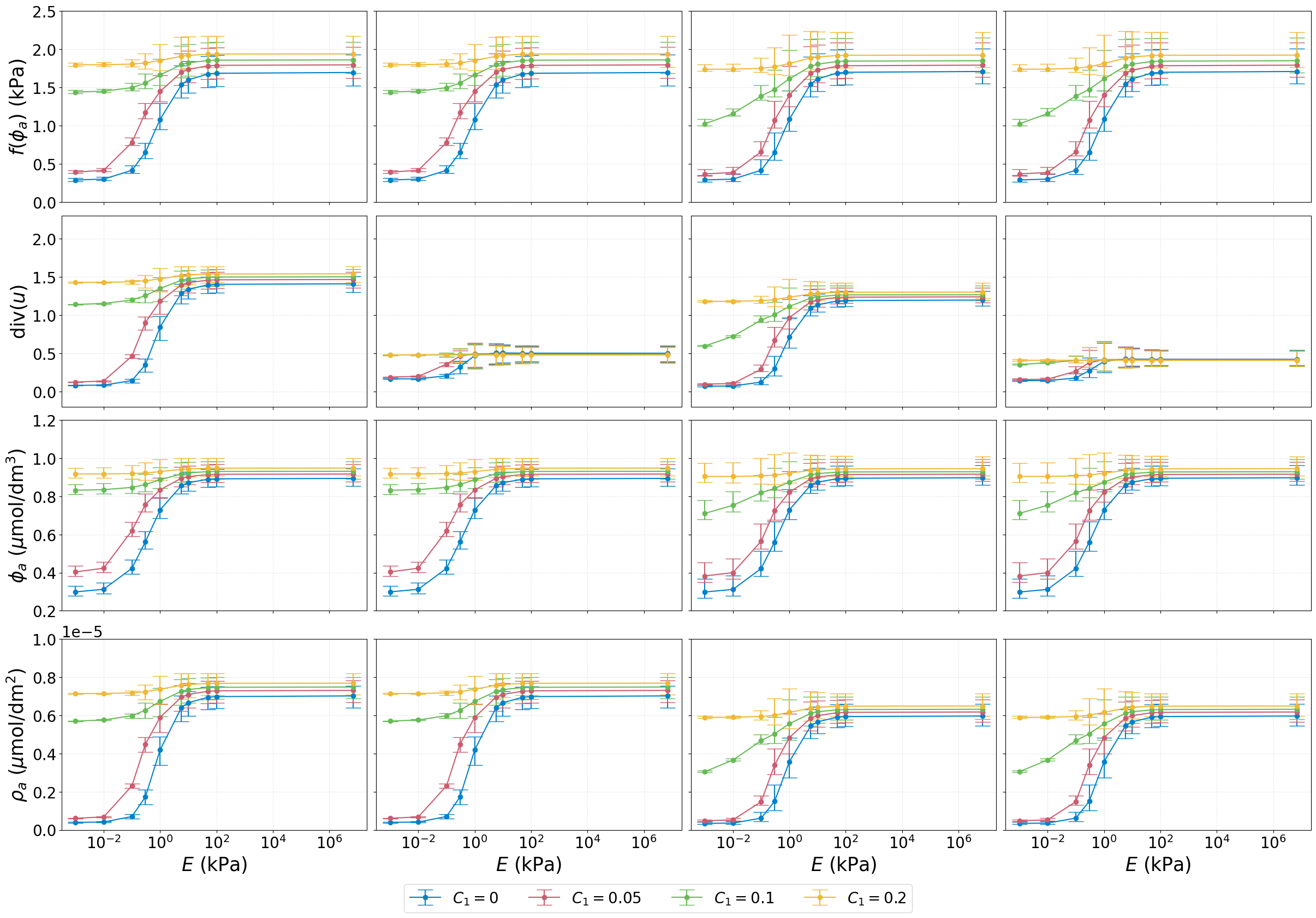}
\end{minipage}
  \caption{{\bf Simulation results showing the mean, $\frac{1}{|\Omega|}\int_\Omega \cdot \dd{x}$, min and max values of $f(\phi_a)$, ${\rm div}(u)$, $\phi_a$ and $\rho_a$ as  functions of substrate stiffness $E$, in the case of the model in Eqs~\eqref{eq:coupled_elast}, \eqref{eq:coupled_force_boundary} and \eqref{eq:coupled_reactions} and $3$D stimulus.}
    We consider different couplings, four different values for $C_1$, and two different shapes at $T=100$~s by which time the results are at a steady state. All other parameter values as in Table~\ref{tab:parameters}.}
  \label{fig:comp}
\end{figure}

Simulation results for a $3$D stimulus that models a cell surrounded by the  extracellular matrix are presented in Fig~\ref{fig:sim_subs}-\ref{fig:comp}. Corresponding results of the evolution of the mean of $f(\phi_a)$, ${\rm div}(u)$, $\phi_a$ and $\rho_a$ can be found in \nameref{S1_Appendix}, Section A.2. The results show the same differences between the different couplings as in Fig~\ref{fig:sim_subs_partfixed}-\ref{fig:comp_partfixed}. Comparing Fig~\ref{fig:sim_subs_partfixed} and~\ref{fig:sim_subs}, the results for the concentrations $\phi_a$ and $\rho_a$ are indistinguishable, however there is a clear difference in deformation of the bottom of the cell and in the case of the fixed vertical deformations the magnitude of the deformation at the base of the cell is slightly lower than in the case of force boundary conditions. The same differences are observed for the lamellipodium shape case, see \nameref{S1_Fig} and \nameref{S2_Fig}. The concentrations of the signalling molecules are also less sensitive to parameter changes than the volume change, as can be seen from the parameter analysis in \nameref{S1_Appendix}, Section A.6.
Comparing  Fig~\ref{fig:comp_partfixed} and~\ref{fig:comp}, the main  difference is in the behaviour of  ${\rm div}(u)$ as a function of $E$. Even though the average volume change is the same, we see differences in the maximum and minimum values of the local volume change across the domain. In particular, the maximum local volume change when considering  the model with a partially fixed boundary is larger and is located  on the base of the cell, while the maximum local volume change when considering the model with the force boundary conditions is smaller, but the cell deforms more evenly in all directions. 

\clearpage
\subsubsection*{Numerical simulations in the case of  $2$xD stimulus}\label{subsec:2xD_1}
\begin{figure}[!ht]
\begin{minipage}{\textwidth}
\begin{minipage}{0.05\textwidth}\centering
\vspace{-0.5cm}
    \textbf{(A)}\\ \vspace{0.3cm}
        $\rho_a$ \\ \vspace{0.65cm}
        $\phi_d$ \\ \vspace{0.65cm}
        $\phi_a$ \\ \vspace{0.65cm}
        $\vert u \vert$ \\ \vspace{0.65cm}
    \textbf{(C)}\\ \vspace{0.3cm}
        $\rho_a$ \\ \vspace{0.65cm}
        $\phi_d$ \\ \vspace{0.65cm}
        $\phi_a$ \\ \vspace{0.65cm}
        $\vert u \vert$ 
\end{minipage}\hfill
\begin{minipage}{0.45\textwidth}
\begin{minipage}{0.33\textwidth} \centering
    {\footnotesize $0.1$kPa } \\ 
        \includegraphics[width=\textwidth]{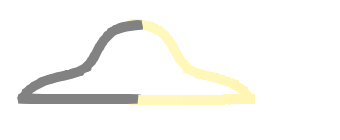}
        \includegraphics[width=\textwidth]{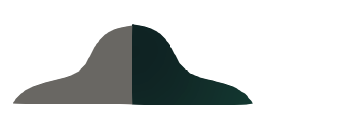}
        \includegraphics[width=\textwidth]{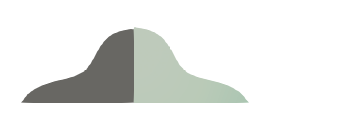}
        \includegraphics[width=\textwidth]{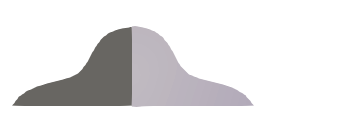}\\
        \vspace{0.5cm}
        \includegraphics[width=\textwidth]{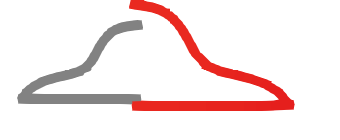}
        \includegraphics[width=\textwidth]{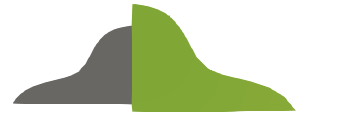}
        \includegraphics[width=\textwidth]{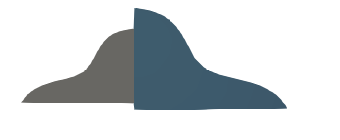}
        \includegraphics[width=\textwidth]{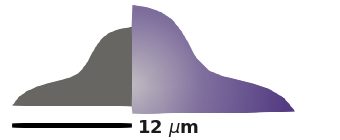}
\end{minipage}\hfill
\begin{minipage}{0.33\textwidth} \centering
    {\footnotesize $5.7$kPa } \\
        \includegraphics[width=\textwidth]{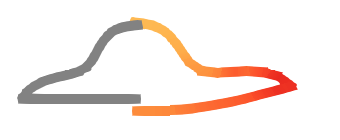}
        \includegraphics[width=\textwidth]{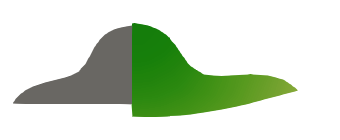}
        \includegraphics[width=\textwidth]{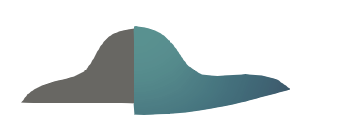}
        \includegraphics[width=\textwidth]{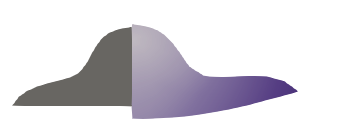}\\
        \vspace{0.5cm}
        \includegraphics[width=\textwidth]{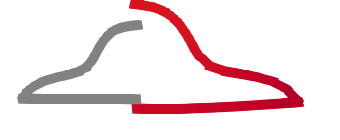}
        \includegraphics[width=\textwidth]{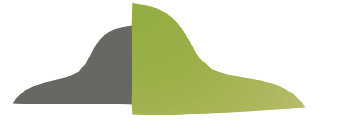}
        \includegraphics[width=\textwidth]{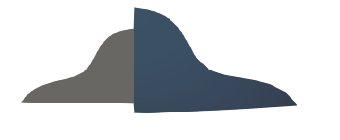}
        \includegraphics[width=\textwidth]{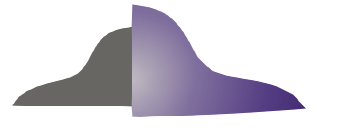}
\end{minipage}\hfill
\begin{minipage}{0.33\textwidth} \centering
    {\footnotesize $7$GPa} \\
        \includegraphics[width=\textwidth]{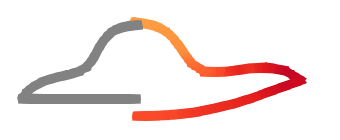}
        \includegraphics[width=\textwidth]{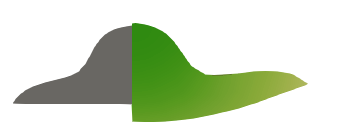}
        \includegraphics[width=\textwidth]{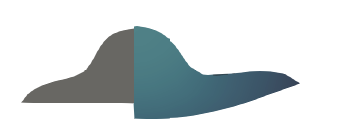}
        \includegraphics[width=\textwidth]{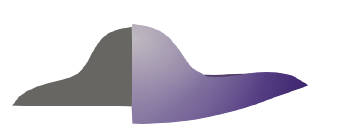}\\
        \vspace{0.5cm}
        \includegraphics[width=\textwidth]{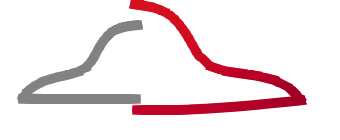}
        \includegraphics[width=\textwidth]{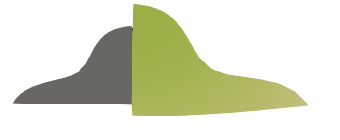}
        \includegraphics[width=\textwidth]{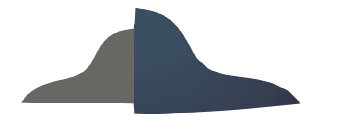}
        \includegraphics[width=\textwidth]{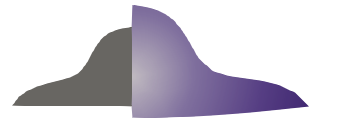}
\end{minipage}\hfill 
\end{minipage}\hfill 
\begin{minipage}{0.05\textwidth}\centering
\vspace{-0.5cm}
    \textbf{(B)}\\ \vspace{4.5cm}
    \textbf{(D)}\\ \vspace{3.8cm}
\end{minipage}\hfill
\begin{minipage}{0.45\textwidth}
\begin{minipage}{0.33\textwidth} \centering
    {\footnotesize $0.1$kPa } \\ 
        \includegraphics[width=\textwidth]{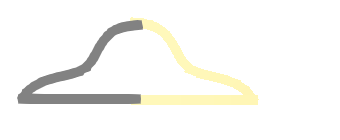}
        \includegraphics[width=\textwidth]{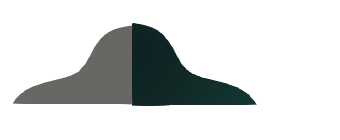}
        \includegraphics[width=\textwidth]{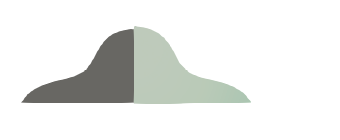}
        \includegraphics[width=\textwidth]{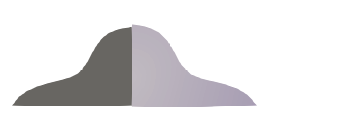}\\
        \vspace{0.5cm}
        \includegraphics[width=\textwidth]{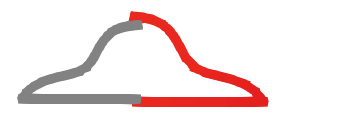}
        \includegraphics[width=\textwidth]{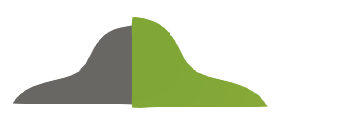}
        \includegraphics[width=\textwidth]{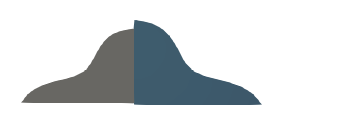}
        \includegraphics[width=\textwidth]{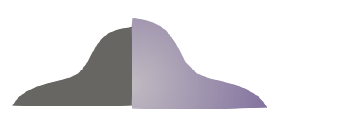}
\end{minipage}\hfill
\begin{minipage}{0.33\textwidth} \centering
    {\footnotesize $5.7$kPa } \\
        \includegraphics[width=\textwidth]{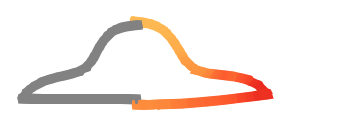}
        \includegraphics[width=\textwidth]{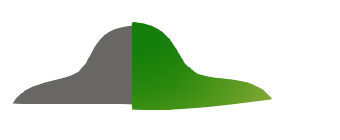}
        \includegraphics[width=\textwidth]{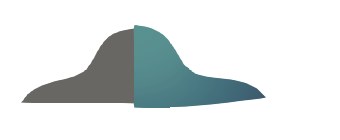}
        \includegraphics[width=\textwidth]{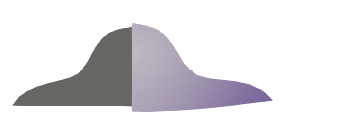}\\
        \vspace{0.5cm}
        \includegraphics[width=\textwidth]{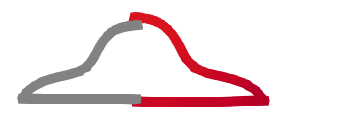}
        \includegraphics[width=\textwidth]{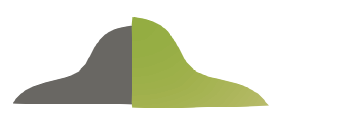}
        \includegraphics[width=\textwidth]{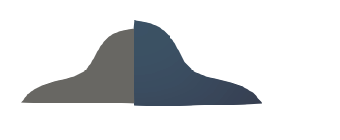}
        \includegraphics[width=\textwidth]{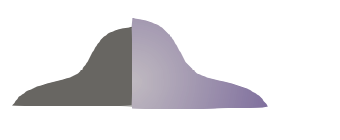}
\end{minipage}\hfill
\begin{minipage}{0.33\textwidth} \centering
    {\footnotesize $7$GPa} \\
        \includegraphics[width=\textwidth]{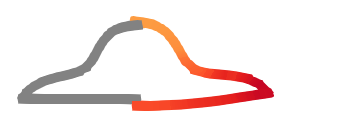}
        \includegraphics[width=\textwidth]{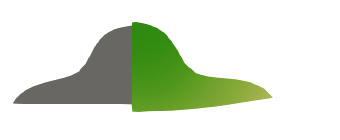}
        \includegraphics[width=\textwidth]{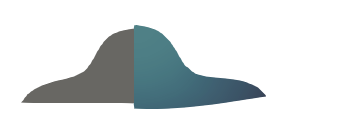}
        \includegraphics[width=\textwidth]{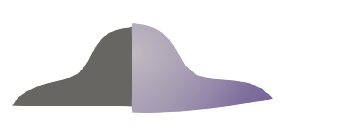}\\
        \vspace{0.5cm}
        \includegraphics[width=\textwidth]{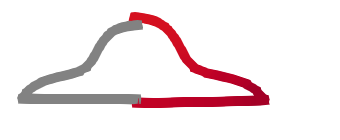}
        \includegraphics[width=\textwidth]{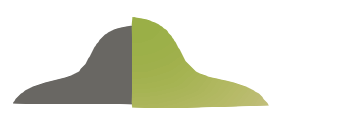}
        \includegraphics[width=\textwidth]{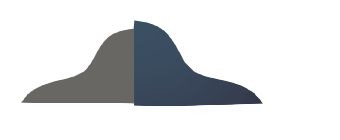}
        \includegraphics[width=\textwidth]{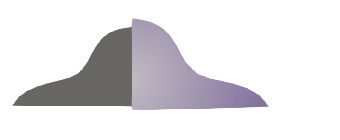}
\end{minipage}\hfill
\end{minipage}
\end{minipage}\hfill
\begin{minipage}{\textwidth}\raggedright
    \hspace{1cm}
    \includegraphics[width=0.22\textwidth]{figures/coupled_p_bar_grouped.png}
    \includegraphics[width=0.22\textwidth]{figures/coupled_cd_bar_grouped.png}
    \includegraphics[width=0.22\textwidth]{figures/coupled_ca_bar_grouped.png}
    \includegraphics[width=0.22\textwidth]{figures/coupled_u_bar_grouped.png}
\end{minipage}
  \caption{{\bf Numerical simulation results showing $\rho_a$, $\phi_d$, $\phi_a$ and $\vert u\vert$ simulation results  showing $\phi_a$ and $\rho_a$ for the model in Eqs~\eqref{eq:coupled_elast}, \eqref{eq:coupled_force_boundary}, and   \eqref{eq:coupled_reactions} for the axisymmetric shape and in the case of $2$xD stimulus at a steady state at $T=100$~s.}
    Four different scenarios are considered: \textbf{(A)} $C_1=0$~$(\mbox{kPa s})^{-1}$ $(\sigma\not\rightarrow \phi_a)$ and $E_c=0.6$~kPa $(\phi_a\not\rightarrow E_c)$; \textbf{(B)} $C_1=0$~$(\mbox{kPa s})^{-1}$ $(\sigma\not\rightarrow \phi_a)$ and $E_c=f(\phi_a)$ $(\phi_a\rightarrow E_c)$; \textbf{(C)} $C_1=0.1$~$(\mbox{kPa s})^{-1}$ $(\sigma\rightarrow \phi_a)$ and $E_c=0.6$~kPa $(\phi_a\not\rightarrow E_c)$; \textbf{(D)} $C_1=0.1$~$(\mbox{kPa s})^{-1}$ $(\sigma\rightarrow \phi_a)$ and $E_c=f(\phi_a)$ $(\phi_a\rightarrow E_c)$. Within each subfigure, the rows represent $\rho_a$, $\phi_d$, $\phi_a$ and $\vert u\vert$ on a cross-section of the plane $x_1=0$ of the axisymmetric cell, and the columns represent $E=0.1, 5.7, 7\cdot 10^6$~kPa. Parameter values as in Table~\ref{tab:parameters}.}
  \label{fig:sim_subs_2D}
\end{figure}

In Fig~\ref{fig:sim_subs_2D}, and in \nameref{S3_Fig} and \nameref{S4_Fig}  we report on simulation results in the case of $2$xD stimulus and force boundary conditions applied to the entire boundary. For the concentrations, the results are similar to the results in the case of $2$xD stimulus and no vertical deformation on the bottom of the cell, see Fig~\ref{fig:sim_subs_partfixed_2D} and~\ref{fig:sim_subs_2D}. 
However, the results for the deformation are different compared to the previous results. In Fig~\ref{fig:sim_subs_2D}, the cell does not just expand but changes shape as the edges of the cell deform upwards, which is not possible in the case of the partially fixed boundary as we assume no vertical deformation at the base. The deformation of the cell upwards can also be observed in the case of the $3$D stimulus, but it is smaller due to the impact of the ECM surrounding the cell, see Fig~\ref{fig:sim_subs}. 
We observe that for $C_1=0$~$(\mbox{kPa s})^{-1}$ the cell deforms upwards a little more than for $C_1=0.1$~$(\mbox{kPa s})^{-1}$. This is due to the larger variation in the concentration $\rho_a$ for $C_1=0$~$(\mbox{kPa s})^{-1}$ compared to $C_1=0.1$~$(\mbox{kPa s})^{-1}$. 
The same features are observed for the lamellipodium shape, see Fig~\ref{fig:sim_lamel_partfixed_2D}, and in \nameref{S2_Fig}, \nameref{S3_Fig} and \nameref{S4_Fig}. 

\section*{Discussion and  Conclusion}\label{sec:discussion}
We have derived a model for mechanotransduction via the RhoA signalling pathway with ECM stiffness and intracellular mechanical properties serving as the mechanical cues.  The modelling extends the work of~\cite{scott_spatial_2021} incorporating the explicit modelling of cell deformation based on an elastic constitutive assumption. We have extended on \cite{scott_spatial_2021, sun_computational_2016, eroume_exploring_2021} and introduced a two-way coupling between the mechanics of the cell and biochemical signalling processes. This  two-way coupling appears to be central to mechanical homeostasis  observed in biological experiments~\cite{Grolleman_2023}. We propose a robust numerical method, based on the bulk-surface finite element method (FEM), see e.g.~\cite{dziuk2013finite}, for the approximation of the model  and report on simulation results for different scenarios, validating the results by comparison with simulations presented in~\cite{scott_spatial_2021} and experimental observations in~\cite{beamish_engineered_2017}. Namely, we considered different levels of substrate stiffness for cells of different shapes that either sit on a rigid flat substrate or are embedded in a three-dimensional substrate.  

Our broad conclusions are that cell shape strongly influences the dynamics of the signalling molecules and the deformation of the cell, as seen in all figures comparing the axisymmetric and lamellipodium shape, where the emergent patterns differ, which is in line with experimental observations~\cite{chen_cell_2003, mcbeath_cell_2004}. 
Cell shape also affects experimentally observed features such as the threshold-like response to changes in substrate stiffness~\cite{beamish_engineered_2017} which is reproduced by the model. In Fig~\ref{fig:comp_partfixed_2D}, \ref{fig:comp_partfixed}, \ref{fig:comp} and in \nameref{S4_Fig}, we see that for certain parameters ($C_1=0.1$~$(\mbox{kPa s})^{-1}$ and low substrate stiffness), the cell shape affects the mean concentrations of the signalling molecules and the mean volume change of the cell, and thus changes the  threshold-like response. 

Our simulations exhibit novel emergent features, that are inaccessible without the framework we propose, such as the bidirectional coupling between mechanics and signalling processes through allowing the Young's modulus of the cell to depend on protein concentration that can allow for robustness in terms of the magnitude of deformation in response to differences in substrate stiffness. This is an example of a mechanical homeostasis mechanism that emerges only at this level of modelling complexity which is of relevance to biology~\cite{Grolleman_2023}. Other instances of mechanical homeostasis are the stress being maintained in the cardiovascular system under mechanical perturbations~\cite{kassab_biomechanical_2024} and the tensional homeostasis by the RhoA signalling pathway at the level of multiple cells~\cite{andersen_cell_2023, weaver_cellular_2016}, which is known to be governed by cellular stiffness sensing~\cite{chanduri_cellular_2024}. Another mechanism that experiences homeostatic response to substrate stiffness is that of the mechanical memory of the cell, describing the phenomenon of a cell responding less to substrates with lower stiffness if they have been cultured on stiff substrates~\cite{weaver_cellular_2016, cacopardo_characterizing_2022}. Due to the bidirectional coupling between the mechanics and the chemistry in our modelling framework, an extension of this work by changing the chosen couplings could be used to model these other mechanical homeostasis phenomena.  

Based on previous biological studies \cite{scott_spatial_2021, gardel_elastic_2004}, we considered cases in which the mechanical properties of the cell (cell stiffness) depend on the concentration of signalling molecules. This coupling yields less sensitivity of total deformation to substrate stiffness whilst leaving the dynamics of the signalling molecules themselves broadly unchanged, see Fig~\ref{fig:comp_partfixed_2D}, \ref{fig:comp_partfixed}, \ref{fig:comp} and in \nameref{S4_Fig}. The insensitivity of the dynamics of the signalling molecules to deformation levels arises since in the model proposed here they are influenced by the local stress rather than deformation.  
We note that the above constitutes another emergent homeostasis mechanism that the modelling framework allows us to explore.
We stress that our work serves as an example of how mechanotransduction may be modelled and more complicated models for the mechanics, biochemistry and couplings therefore are warranted based on the remarkable emergent features we observe even in our relatively simple setting. We expect such models to be particularly fruitful avenues for future work.

One such example of more complicated models for the mechanics could include a viscoelastic or poroelastic constitutive law. As presented in~Section A.7 in \nameref{S1_Appendix}, the assumption of a (linear) viscoelastic constitutive law leads to qualitatively similar results to those presented in this work for a purely (linear) elastic constitutive law.
Our current assumption of a linear elastic constitutive law for the mechanics of the cell is limiting as it assumes small deformations. This framework needs to be extended to study the effect of large deformations and shape changes, which would include the effect these deformations have on the signalling molecules. This would be especially interesting, as this study shows that cell shape is one of the determinants affecting the mean concentrations of the signalling molecules. 

The boundary conditions for the deformation we consider correspond to simple models of a cell in vitro (flat $2$D substrate) or in vivo (homogeneous $3$D matrix). We see that the cell on a $2$D substrate appears to spread radially with minimal deformation orthogonal to the substrate while the latter exhibits a more uniform although smaller in total magnitude $3$D deformation. Differences in deformation for different environments are in line with the literature as the effect of the substrate stiffness on cells is known to vary in 2D and 3D substrates~\cite{byfield_endothelial_2009}. An interesting extension that could be included in the above framework would be spatial variations in substrate stiffness or more complicated models for the substrate mechanics both of which are of much biological relevance~\cite{chaudhuri_substrate_2015, isomursuDirected2022,ross_physical_2012}.

This work shows how mechanistic modelling of mechanotransduction can reveal remarkable emergent properties. It lays the groundwork for future studies where further complexity can be added as required to model specific signalling pathways or  to reflect other mechanical models derived from different constitutive assumptions. We anticipate that choosing a viscoelastic or poroelastic constitutive law for the mechanics of the cell is an interesting direction for future studies, as this is in line with recent experimental observations~\cite{kasza_cell_2007, moeendarbary_cytoplasm_2013}. 
Given the fact that cell shape greatly influences the dynamics of the cell, as shown in this work, other reference geometries are also of interest as a subject for future work.  Extending the signalling model of \cite{scott_spatial_2021} further, we intend to couple the model of this work with a similar biomechanical model for the deformation of the nucleus coupled with the dynamics of signalling molecules within the nucleus, such as the YAP/TAZ pathway~\cite{jafarinia_insights_2024}.

\section*{Acknowledgments}
SV was supported by the EPSRC Centre for Doctoral Training in Mathematical Modelling, Analysis and Computation (MAC-MIGS) funded by the UK Engineering and Physical Sciences Research Council (grant EP/S023291/1), Heriot-Watt University and the 
University of Edinburgh. CV acknowledges support from the Dr Perry James (Jim) Browne Research Centre on Mathematics and its Applications (University of Sussex).\\
SV and MP would like to thank the Isaac Newton Institute for Mathematical Sciences, Cambridge, for support and hospitality during the research programme 'Uncertainty quantification and stochastic modelling of materials', EPSRC Grant Number EP/R014604/1, where some work on the manuscript was undertaken. 
\\
The authors would like to thank Padmini Rangamani for helpful discussions. 

%
%

\label{zzz}

\clearpage

\section*{S1 Fig}
\label{S1_Fig}
\begin{figure}[!ht]
\begin{minipage}{\textwidth}
\begin{minipage}{0.05\textwidth}\centering
\vspace{-0.5cm}
    \textbf{(A)}\\ \vspace{0.3cm}
        $\rho_a$ \\ \vspace{1.1cm}
        $\phi_d$ \\ \vspace{1.1cm}
        $\phi_a$ \\ \vspace{1.1cm}
        $\vert u \vert$ \\ \vspace{1.1cm}
    \textbf{(C)}\\ \vspace{0.3cm}
        $\rho_a$ \\ \vspace{1.1cm}
        $\phi_d$ \\ \vspace{1.1cm}
        $\phi_a$ \\ \vspace{1.1cm}
        $\vert u \vert$ 
\end{minipage}\hfill
\begin{minipage}{0.45\textwidth}
\begin{minipage}{0.33\textwidth} \centering
    {\footnotesize $0.1$kPa } \\ 
        \includegraphics[width=\textwidth]{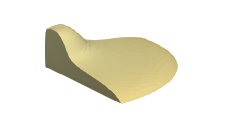}
        \includegraphics[width=\textwidth]{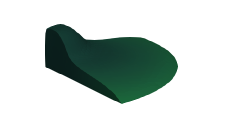}
        \includegraphics[width=\textwidth]{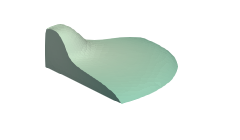}
        \includegraphics[width=\textwidth]{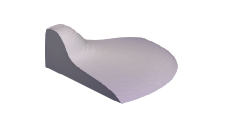}\\
        \vspace{0.5cm}
        \includegraphics[width=\textwidth]{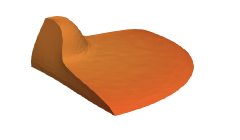}
        \includegraphics[width=\textwidth]{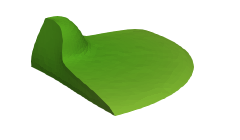}
        \includegraphics[width=\textwidth]{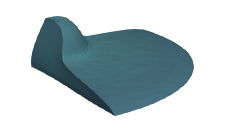}
        \includegraphics[width=\textwidth]{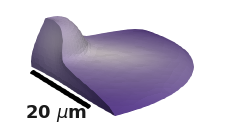}
\end{minipage}\hfill
\begin{minipage}{0.33\textwidth} \centering
    {\footnotesize $5.7$kPa } \\
        \includegraphics[width=\textwidth]{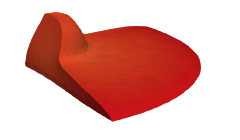}
        \includegraphics[width=\textwidth]{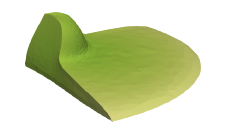}
        \includegraphics[width=\textwidth]{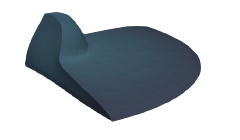}
        \includegraphics[width=\textwidth]{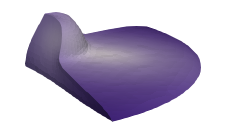}\\
        \vspace{0.5cm}
        \includegraphics[width=\textwidth]{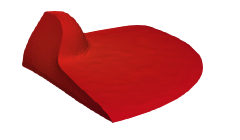}
        \includegraphics[width=\textwidth]{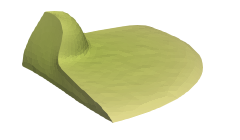}
        \includegraphics[width=\textwidth]{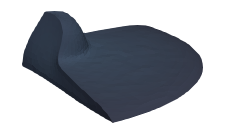}
        \includegraphics[width=\textwidth]{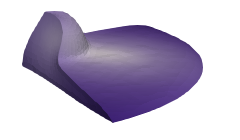}
\end{minipage}\hfill
\begin{minipage}{0.33\textwidth} \centering
    {\footnotesize $7$GPa} \\
        \includegraphics[width=\textwidth]{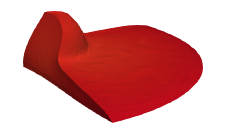}
        \includegraphics[width=\textwidth]{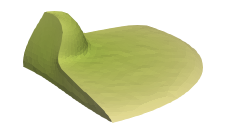}
        \includegraphics[width=\textwidth]{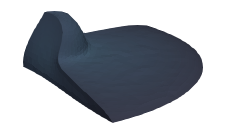}
        \includegraphics[width=\textwidth]{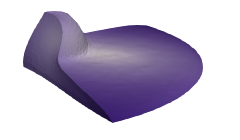}\\
        \vspace{0.5cm}
        \includegraphics[width=\textwidth]{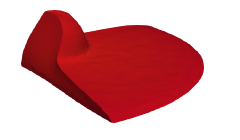}
        \includegraphics[width=\textwidth]{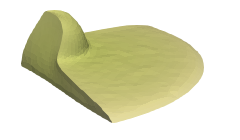}
        \includegraphics[width=\textwidth]{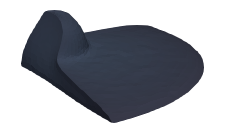}
        \includegraphics[width=\textwidth]{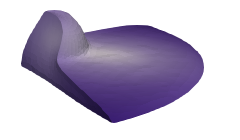}
\end{minipage}\hfill 
\end{minipage}\hfill 
\begin{minipage}{0.05\textwidth}\centering
\vspace{-0.5cm}
    \textbf{(B)}\\ \vspace{6.5cm}
    \textbf{(D)}\\ \vspace{5.3cm}
\end{minipage}\hfill
\begin{minipage}{0.45\textwidth}
\begin{minipage}{0.33\textwidth} \centering
    {\footnotesize $0.1$kPa } \\ 
        \includegraphics[width=\textwidth]{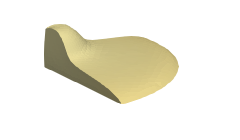}
        \includegraphics[width=\textwidth]{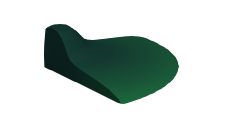}
        \includegraphics[width=\textwidth]{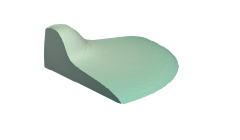}
        \includegraphics[width=\textwidth]{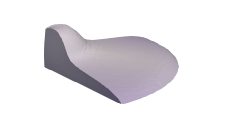}\\
        \vspace{0.5cm}
        \includegraphics[width=\textwidth]{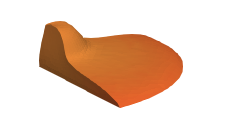}
        \includegraphics[width=\textwidth]{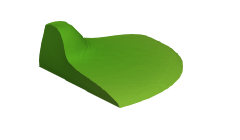}
        \includegraphics[width=\textwidth]{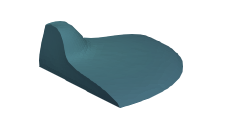}
        \includegraphics[width=\textwidth]{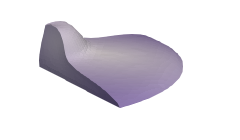}
\end{minipage}\hfill
\begin{minipage}{0.33\textwidth} \centering
    {\footnotesize $5.7$kPa } \\
        \includegraphics[width=\textwidth]{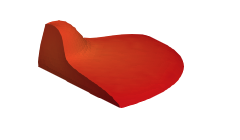}
        \includegraphics[width=\textwidth]{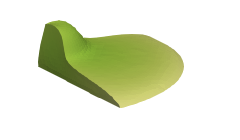}
        \includegraphics[width=\textwidth]{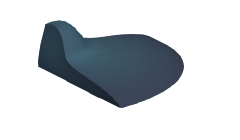}
        \includegraphics[width=\textwidth]{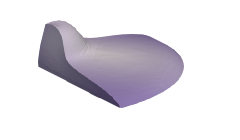}\\
        \vspace{0.5cm}
        \includegraphics[width=\textwidth]{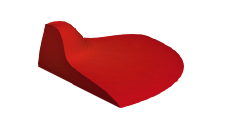}
        \includegraphics[width=\textwidth]{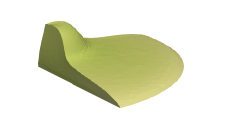}
        \includegraphics[width=\textwidth]{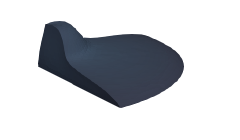}
        \includegraphics[width=\textwidth]{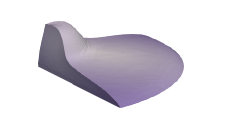}
\end{minipage}\hfill
\begin{minipage}{0.33\textwidth} \centering
    {\footnotesize $7$GPa} \\
        \includegraphics[width=\textwidth]{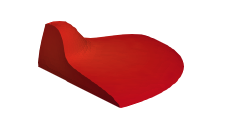}
        \includegraphics[width=\textwidth]{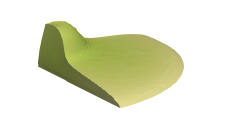}
        \includegraphics[width=\textwidth]{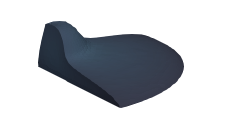}
        \includegraphics[width=\textwidth]{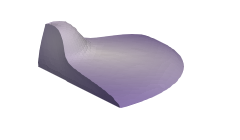}\\
        \vspace{0.5cm}
        \includegraphics[width=\textwidth]{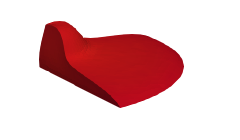}
        \includegraphics[width=\textwidth]{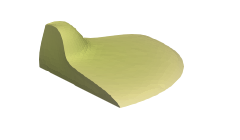}
        \includegraphics[width=\textwidth]{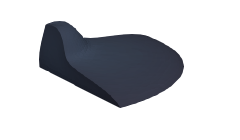}
        \includegraphics[width=\textwidth]{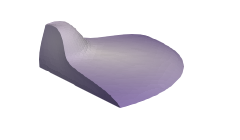}
\end{minipage}\hfill
\end{minipage}
\end{minipage}\hfill
\begin{minipage}{\textwidth}\raggedright
    \hspace{1cm}
    \includegraphics[width=0.22\textwidth]{figures/coupled_p_bar_grouped.png}
    \includegraphics[width=0.22\textwidth]{figures/coupled_cd_bar_grouped.png}
    \includegraphics[width=0.22\textwidth]{figures/coupled_ca_bar_grouped.png}
    \includegraphics[width=0.22\textwidth]{figures/coupled_u_lamel_grouped_bar.png}
\end{minipage}
  \caption{{\bf Numerical simulation results showing $\rho_a$, $\phi_d$, $\phi_a$ and $\vert u\vert$  for the model in Eqs~\eqref{eq:coupled_elast}-\eqref{eq:coupled_reactions}  for the lamellipodium shape  and in the case of the $3$D stimulus at a steady state at $T=100$~s.}
  Four different scenarios are considered: \textbf{(A)} $C_1=0$~$(\mbox{kPa s})^{-1}$ $(\sigma\not\rightarrow \phi_a)$ and $E_c=0.6$~kPa $(\phi_a\not\rightarrow E_c)$; \textbf{(B)} $C_1=0$~$(\mbox{kPa s})^{-1}$ $(\sigma\not\rightarrow \phi_a)$ and $E_c=f(\phi_a)$ $(\phi_a\rightarrow E_c)$; \textbf{(C)} $C_1=0.1$~$(\mbox{kPa s})^{-1}$ $(\sigma\rightarrow \phi_a)$ and $E_c=0.6$~kPa $(\phi_a\not\rightarrow E_c)$; \textbf{(D)} $C_1=0.1$~$(\mbox{kPa s})^{-1}$ $(\sigma\rightarrow \phi_a)$ and $E_c=f(\phi_a)$ $(\phi_a\rightarrow E_c)$. Within each subfigure, the rows represent $\rho_a$, $\phi_d$, $\phi_a$ and $\vert u\vert$ on the surface of the cell, and the columns represent $E=0.1, 5.7, 7\cdot 10^6$~kPa. Parameter values as in Table~\ref{tab:parameters}.}
\end{figure}

\clearpage 

\section*{S2 Fig}
\label{S2_Fig}
\begin{figure}[!ht]
\begin{minipage}{\textwidth}
\begin{minipage}{0.05\textwidth}\centering
\vspace{-0.5cm}
    \textbf{(A)}\\ \vspace{0.3cm}
        $\rho_a$ \\ \vspace{1.1cm}
        $\phi_d$ \\ \vspace{1.1cm}
        $\phi_a$ \\ \vspace{1.1cm}
        $\vert u \vert$ \\ \vspace{1.1cm}
    \textbf{(C)}\\ \vspace{0.3cm}
        $\rho_a$ \\ \vspace{1.1cm}
        $\phi_d$ \\ \vspace{1.1cm}
        $\phi_a$ \\ \vspace{1.1cm}
        $\vert u \vert$ 
\end{minipage}\hfill
\begin{minipage}{0.45\textwidth}
\begin{minipage}{0.33\textwidth} \centering
    {\footnotesize $0.1$kPa } \\ 
        \includegraphics[width=\textwidth]{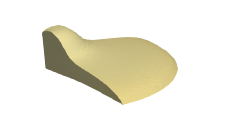}
        \includegraphics[width=\textwidth]{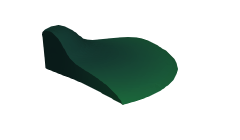}
        \includegraphics[width=\textwidth]{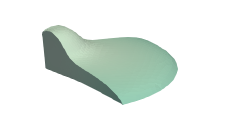}
        \includegraphics[width=\textwidth]{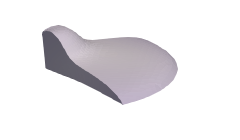}\\
        \vspace{0.5cm}
        \includegraphics[width=\textwidth]{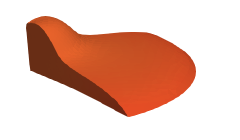}
        \includegraphics[width=\textwidth]{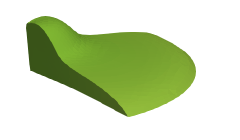}
        \includegraphics[width=\textwidth]{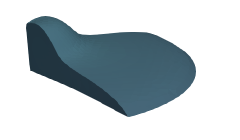}
        \includegraphics[width=\textwidth]{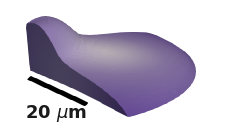}
\end{minipage}\hfill
\begin{minipage}{0.33\textwidth} \centering
    {\footnotesize $5.7$kPa } \\
        \includegraphics[width=\textwidth]{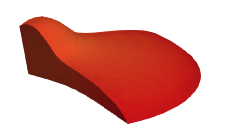}
        \includegraphics[width=\textwidth]{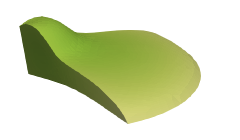}
        \includegraphics[width=\textwidth]{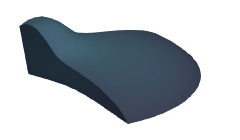}
        \includegraphics[width=\textwidth]{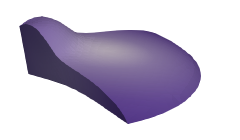}\\
        \vspace{0.5cm}
        \includegraphics[width=\textwidth]{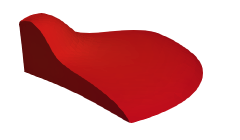}
        \includegraphics[width=\textwidth]{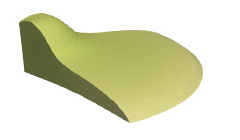}
        \includegraphics[width=\textwidth]{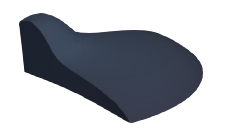}
        \includegraphics[width=\textwidth]{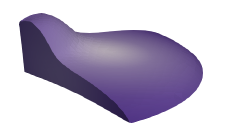}
\end{minipage}\hfill
\begin{minipage}{0.33\textwidth} \centering
    {\footnotesize $7$GPa} \\
        \includegraphics[width=\textwidth]{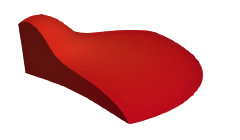}
        \includegraphics[width=\textwidth]{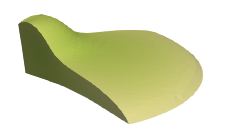}
        \includegraphics[width=\textwidth]{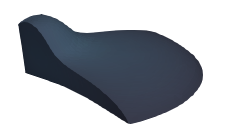}
        \includegraphics[width=\textwidth]{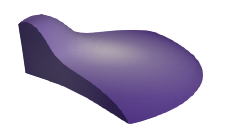}\\
        \vspace{0.5cm}
        \includegraphics[width=\textwidth]{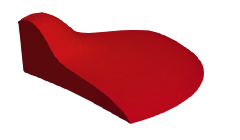}
        \includegraphics[width=\textwidth]{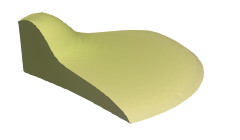}
        \includegraphics[width=\textwidth]{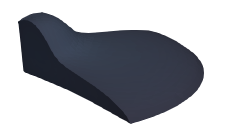}
        \includegraphics[width=\textwidth]{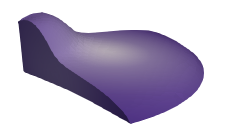}
\end{minipage}\hfill 
\end{minipage}\hfill 
\begin{minipage}{0.05\textwidth}\centering
\vspace{-0.5cm}
    \textbf{(B)}\\ \vspace{6.5cm}
    \textbf{(D)}\\ \vspace{5.3cm}
\end{minipage}\hfill
\begin{minipage}{0.45\textwidth}
\begin{minipage}{0.33\textwidth} \centering
    {\footnotesize $0.1$kPa } \\ 
        \includegraphics[width=\textwidth]{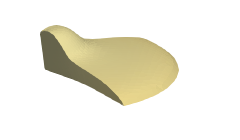}
        \includegraphics[width=\textwidth]{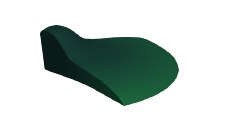}
        \includegraphics[width=\textwidth]{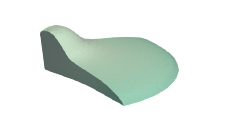}
        \includegraphics[width=\textwidth]{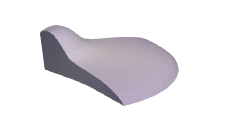}\\
        \vspace{0.5cm}
        \includegraphics[width=\textwidth]{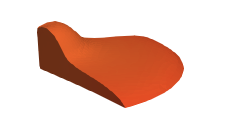}
        \includegraphics[width=\textwidth]{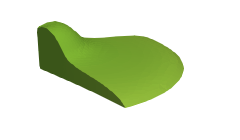}
        \includegraphics[width=\textwidth]{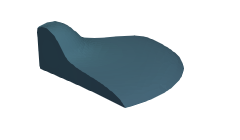}
        \includegraphics[width=\textwidth]{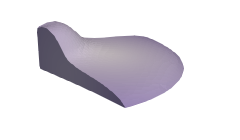}
\end{minipage}\hfill
\begin{minipage}{0.33\textwidth} \centering
    {\footnotesize $5.7$kPa } \\
        \includegraphics[width=\textwidth]{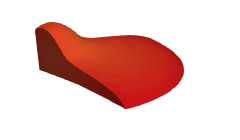}
        \includegraphics[width=\textwidth]{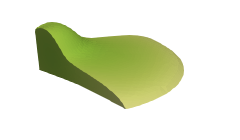}
        \includegraphics[width=\textwidth]{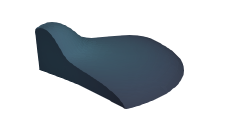}
        \includegraphics[width=\textwidth]{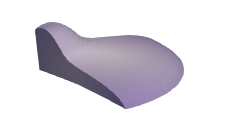}\\
        \vspace{0.5cm}
        \includegraphics[width=\textwidth]{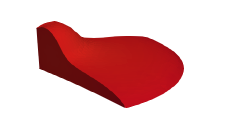}
        \includegraphics[width=\textwidth]{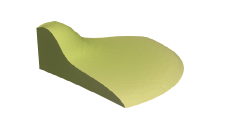}
        \includegraphics[width=\textwidth]{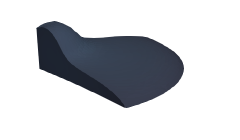}
        \includegraphics[width=\textwidth]{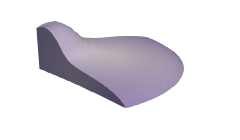}
\end{minipage}\hfill
\begin{minipage}{0.33\textwidth} \centering
    {\footnotesize $7$GPa} \\
        \includegraphics[width=\textwidth]{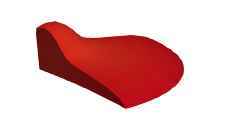}
        \includegraphics[width=\textwidth]{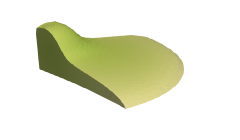}
        \includegraphics[width=\textwidth]{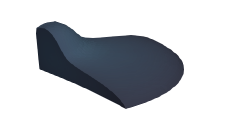}
        \includegraphics[width=\textwidth]{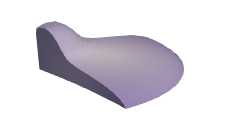}\\
        \vspace{0.5cm}
        \includegraphics[width=\textwidth]{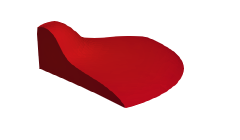}
        \includegraphics[width=\textwidth]{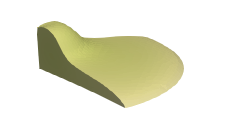}
        \includegraphics[width=\textwidth]{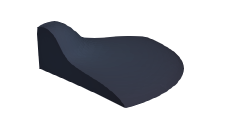}
        \includegraphics[width=\textwidth]{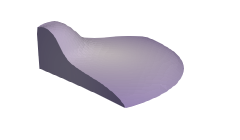}
\end{minipage}\hfill
\end{minipage}
\end{minipage}\hfill
\begin{minipage}{\textwidth}\raggedright
    \hspace{1cm}
    \includegraphics[width=0.22\textwidth]{figures/coupled_p_bar_grouped.png}
    \includegraphics[width=0.22\textwidth]{figures/coupled_cd_bar_grouped.png}
    \includegraphics[width=0.22\textwidth]{figures/coupled_ca_bar_grouped.png}
    \includegraphics[width=0.22\textwidth]{figures/coupled_u_lamel_grouped_bar.png}
\end{minipage}
  \caption{{\bf Numerical simulation results showing $\rho_a$, $\phi_d$, $\phi_a$ and $\vert u\vert$  for the model in Eqs~\eqref{eq:coupled_elast}, \eqref{eq:coupled_force_boundary}, and \eqref{eq:coupled_reactions}  for the lamellipodium shape  and in the case of the $3$D stimulus at a steady state at $T=100$~s.}
  Four different scenarios are considered: \textbf{(A)} $C_1=0$~$(\mbox{kPa s})^{-1}$ $(\sigma\not\rightarrow \phi_a)$ and $E_c=0.6$~kPa $(\phi_a\not\rightarrow E_c)$; \textbf{(B)} $C_1=0$~$(\mbox{kPa s})^{-1}$ $(\sigma\not\rightarrow \phi_a)$ and $E_c=f(\phi_a)$ $(\phi_a\rightarrow E_c)$; \textbf{(C)} $C_1=0.1$~$(\mbox{kPa s})^{-1}$ $(\sigma\rightarrow \phi_a)$ and $E_c=0.6$~kPa $(\phi_a\not\rightarrow E_c)$; \textbf{(D)} $C_1=0.1$~$(\mbox{kPa s})^{-1}$ $(\sigma\rightarrow \phi_a)$ and $E_c=f(\phi_a)$ $(\phi_a\rightarrow E_c)$. Within each subfigure, the rows represent $\rho_a$, $\phi_d$, $\phi_a$ and $\vert u\vert$ on the surface of the cell, and the columns represent $E=0.1, 5.7, 7\cdot 10^6$~kPa. Parameter values as in Table~\ref{tab:parameters}.}
\end{figure}

\clearpage 

\section*{S3 Fig}
\label{S3_Fig}
\begin{figure}[!ht]
\begin{minipage}{\textwidth}
\begin{minipage}{0.05\textwidth}\centering
\vspace{-0.5cm}
    \textbf{(A)}\\ \vspace{0.3cm}
        $\rho_a$ \\ \vspace{1.1cm}
        $\phi_d$ \\ \vspace{1.1cm}
        $\phi_a$ \\ \vspace{1.1cm}
        $\vert u \vert$ \\ \vspace{1.1cm}
    \textbf{(C)}\\ \vspace{0.3cm}
        $\rho_a$ \\ \vspace{1.1cm}
        $\phi_d$ \\ \vspace{1.1cm}
        $\phi_a$ \\ \vspace{1.1cm}
        $\vert u \vert$ 
\end{minipage}\hfill
\begin{minipage}{0.45\textwidth}
\begin{minipage}{0.33\textwidth} \centering
    {\footnotesize $0.1$kPa } \\ 
        \includegraphics[width=\textwidth]{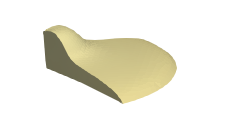}
        \includegraphics[width=\textwidth]{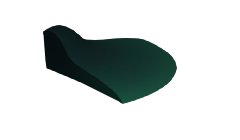}
        \includegraphics[width=\textwidth]{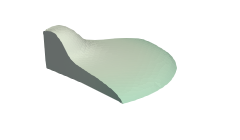}
        \includegraphics[width=\textwidth]{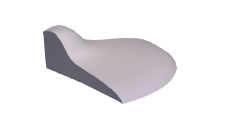}\\
        \vspace{0.5cm}
        \includegraphics[width=\textwidth]{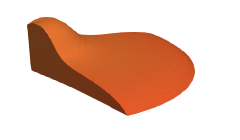}
        \includegraphics[width=\textwidth]{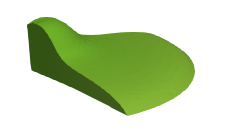}
        \includegraphics[width=\textwidth]{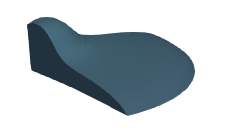}
        \includegraphics[width=\textwidth]{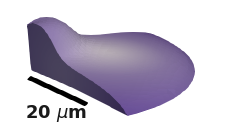}
\end{minipage}\hfill
\begin{minipage}{0.33\textwidth} \centering
    {\footnotesize $5.7$kPa } \\
        \includegraphics[width=\textwidth]{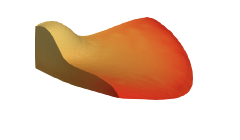}
        \includegraphics[width=\textwidth]{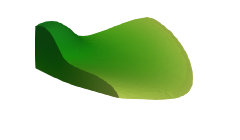}
        \includegraphics[width=\textwidth]{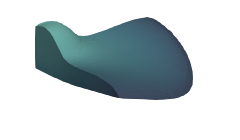}
        \includegraphics[width=\textwidth]{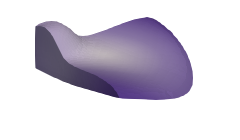}\\
        \vspace{0.5cm}
        \includegraphics[width=\textwidth]{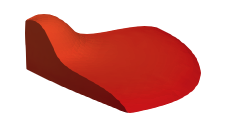}
        \includegraphics[width=\textwidth]{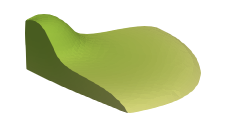}
        \includegraphics[width=\textwidth]{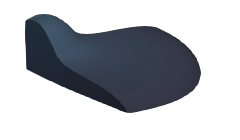}
        \includegraphics[width=\textwidth]{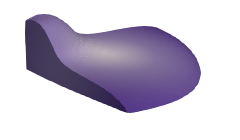}
\end{minipage}\hfill
\begin{minipage}{0.33\textwidth} \centering
    {\footnotesize $7$GPa} \\
        \includegraphics[width=\textwidth]{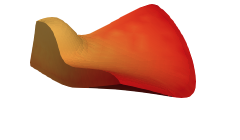}
        \includegraphics[width=\textwidth]{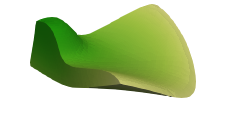}
        \includegraphics[width=\textwidth]{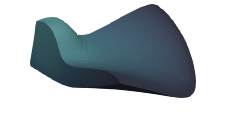}
        \includegraphics[width=\textwidth]{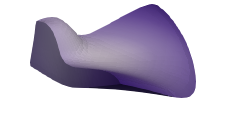}\\
        \vspace{0.5cm}
        \includegraphics[width=\textwidth]{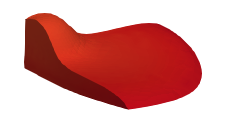}
        \includegraphics[width=\textwidth]{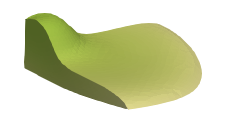}
        \includegraphics[width=\textwidth]{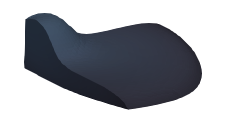}
        \includegraphics[width=\textwidth]{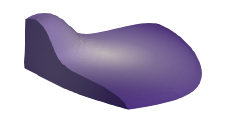}
\end{minipage}\hfill 
\end{minipage}\hfill 
\begin{minipage}{0.05\textwidth}\centering
\vspace{-0.5cm}
    \textbf{(B)}\\ \vspace{6.5cm}
    \textbf{(D)}\\ \vspace{5.3cm}
\end{minipage}\hfill
\begin{minipage}{0.45\textwidth}
\begin{minipage}{0.33\textwidth} \centering
    {\footnotesize $0.1$kPa } \\ 
        \includegraphics[width=\textwidth]{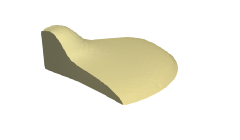}
        \includegraphics[width=\textwidth]{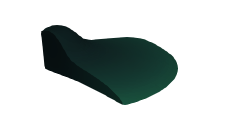}
        \includegraphics[width=\textwidth]{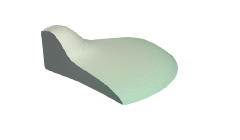}
        \includegraphics[width=\textwidth]{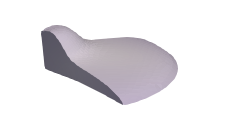}\\
        \vspace{0.5cm}
        \includegraphics[width=\textwidth]{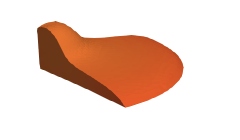}
        \includegraphics[width=\textwidth]{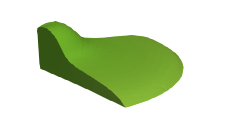}
        \includegraphics[width=\textwidth]{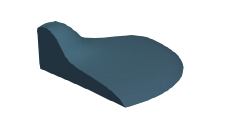}
        \includegraphics[width=\textwidth]{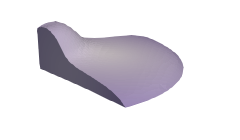}
\end{minipage}\hfill
\begin{minipage}{0.33\textwidth} \centering
    {\footnotesize $5.7$kPa } \\
        \includegraphics[width=\textwidth]{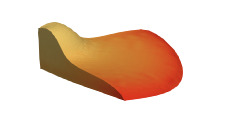}
        \includegraphics[width=\textwidth]{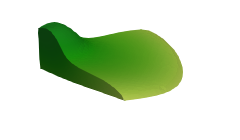}
        \includegraphics[width=\textwidth]{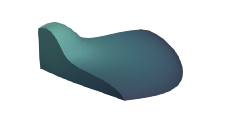}
        \includegraphics[width=\textwidth]{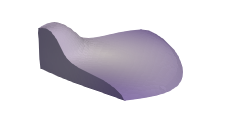}\\
        \vspace{0.5cm}
        \includegraphics[width=\textwidth]{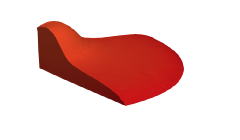}
        \includegraphics[width=\textwidth]{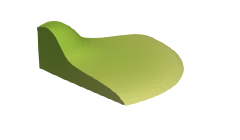}
        \includegraphics[width=\textwidth]{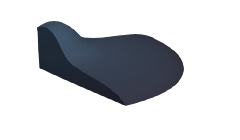}
        \includegraphics[width=\textwidth]{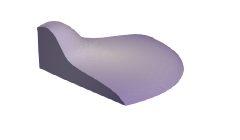}
\end{minipage}\hfill
\begin{minipage}{0.33\textwidth} \centering
    {\footnotesize $7$GPa} \\
        \includegraphics[width=\textwidth]{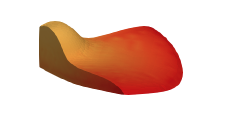}
        \includegraphics[width=\textwidth]{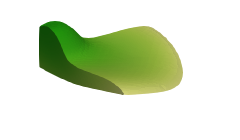}
        \includegraphics[width=\textwidth]{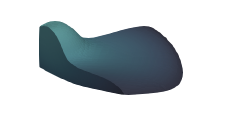}
        \includegraphics[width=\textwidth]{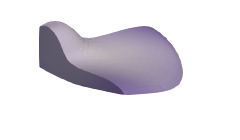}\\
        \vspace{0.5cm}
        \includegraphics[width=\textwidth]{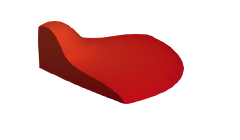}
        \includegraphics[width=\textwidth]{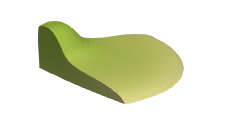}
        \includegraphics[width=\textwidth]{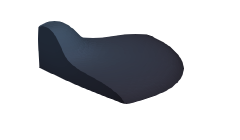}
        \includegraphics[width=\textwidth]{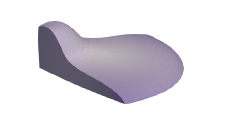}
\end{minipage}\hfill
\end{minipage}
\end{minipage}\hfill
\begin{minipage}{\textwidth}\raggedright
    \hspace{1cm}
    \includegraphics[width=0.22\textwidth]{figures/coupled_p_bar_grouped.png}
    \includegraphics[width=0.22\textwidth]{figures/coupled_cd_bar_grouped.png}
    \includegraphics[width=0.22\textwidth]{figures/coupled_ca_bar_grouped.png}
    \includegraphics[width=0.22\textwidth]{figures/coupled_u_lamel_grouped_bar.png}
\end{minipage}
  \caption{{\bf Numerical simulation results showing $\rho_a$, $\phi_d$, $\phi_a$ and $\vert u\vert$  for the model in Eqs~\eqref{eq:coupled_elast}, \eqref{eq:coupled_force_boundary}, and    \eqref{eq:coupled_reactions} for the  lamellipodium shape cells and $2$xD stimulus at a steady state at $T=100$~s.}
  Four different scenarios are considered: \textbf{(A)} $C_1=0$~$(\mbox{kPa s})^{-1}$ $(\sigma\not\rightarrow \phi_a)$ and $E_c=0.6$~kPa $(\phi_a\not\rightarrow E_c)$; \textbf{(B)} $C_1=0$~$(\mbox{kPa s})^{-1}$ $(\sigma\not\rightarrow \phi_a)$ and $E_c=f(\phi_a)$ $(\phi_a\rightarrow E_c)$; \textbf{(C)} $C_1=0.1$~$(\mbox{kPa s})^{-1}$ $(\sigma\rightarrow \phi_a)$ and $E_c=0.6$~kPa $(\phi_a\not\rightarrow E_c)$; \textbf{(D)} $C_1=0.1$~$(\mbox{kPa s})^{-1}$ $(\sigma\rightarrow \phi_a)$ and $E_c=f(\phi_a)$ $(\phi_a\rightarrow E_c)$. Within each subfigure, the rows represent $\rho_a$, $\phi_d$, $\phi_a$ and $\vert u\vert$ on the surface of the cell, and the columns represent $E=0.1, 5.7, 7\cdot 10^6$~kPa. Parameter values as in Table~\ref{tab:parameters}.}
\end{figure}

\clearpage 

\section*{S4 Fig}
\label{S4_Fig}
\begin{figure}[!h]
  \begin{minipage}{\textwidth} 
    \hspace{3cm} axisymmetric shape \hspace{4.5cm} lamellipodium shape \\ 
    \vspace{-0.2cm}
    \hspace{1.8cm} {\footnotesize $E_c=0.6$~kPa } \hspace{2cm}  {\footnotesize $E_c=f(\phi_a)$ } \hspace{2cm} {\footnotesize $E_c=0.6$~kPa } \hspace{2cm} {\footnotesize $E_c=f(\phi_a)$ } \\ 
    \vspace{-0.2cm}
\end{minipage}\hfill
\begin{minipage}{0.95\textwidth} 
    \includegraphics[width=\textwidth]{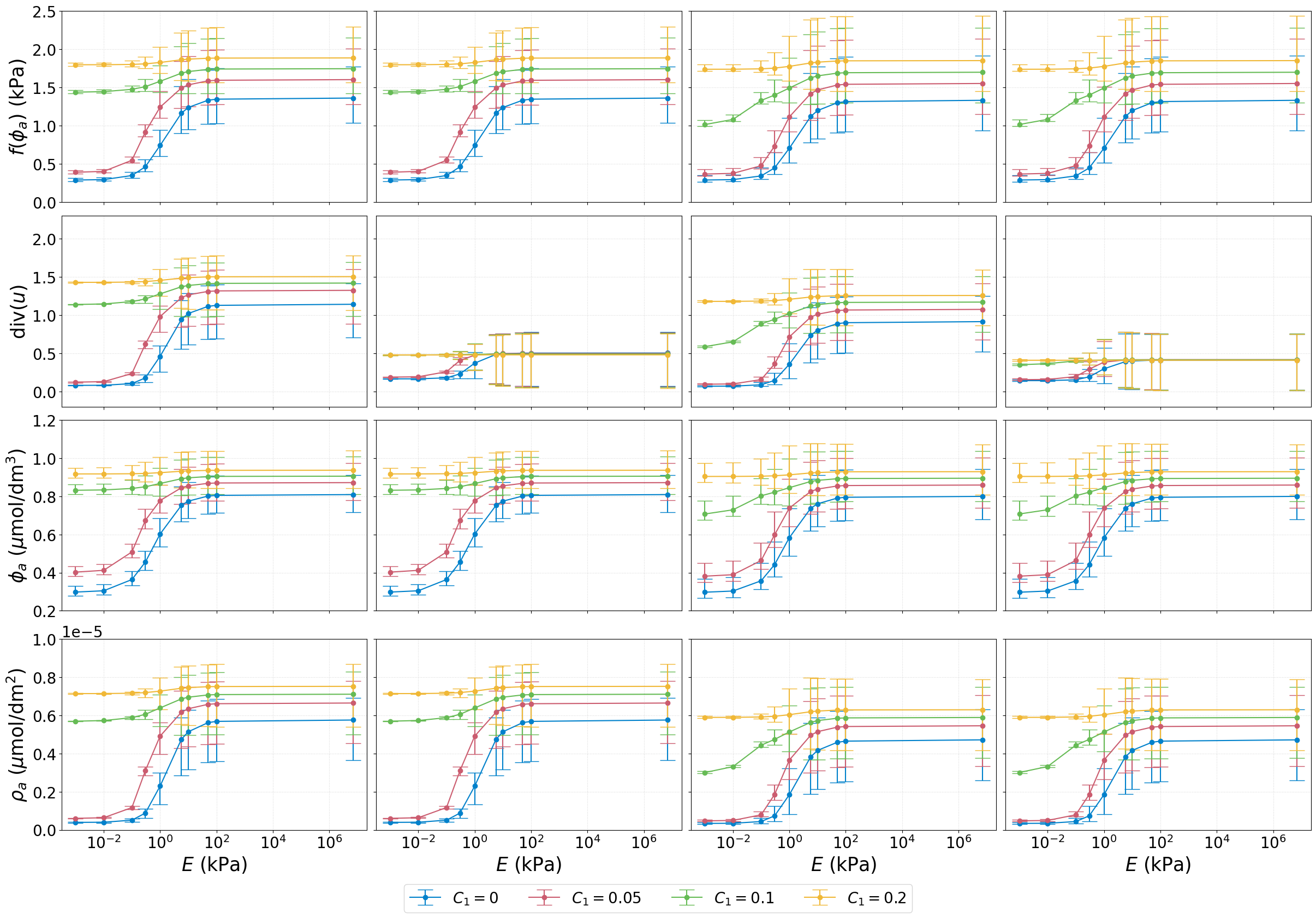}
\end{minipage}
  \caption{{\bf Simulation results showing the mean, $\frac{1}{|\Omega|}\int_\Omega \cdot \dd{x}$, min and max values of $f(\phi_a)$, ${\rm div}(u)$, $\phi_a$ and $\rho_a$ as  functions of substrate stiffness $E$, in the case of the model in Eqs~\eqref{eq:coupled_elast}, \eqref{eq:coupled_force_boundary} and \eqref{eq:coupled_reactions} and $2$xD stimulus.}
   We consider different couplings, four different values for $C_1$, and two different shapes at $T=100$~s by which time the results are at a steady state. All other parameter values as in Table~\ref{tab:parameters}.}
\end{figure}

\pagestyle{empty}
\newgeometry{letterpaper,top=2cm,bottom=2cm,left=2cm,right=2cm,marginparwidth=1.75cm}

\section*{S1 Appendix}
\label{S1_Appendix}
\renewcommand{\thefigure}{A\arabic{figure}}
\renewcommand{\thetable}{A\arabic{table}}
\renewcommand{\thesubsection}{A.\arabic{subsection}}
\renewcommand{\theequation}{A\arabic{equation}}
\setcounter{table}{0}
\setcounter{figure}{0}
\setcounter{equation}{0}

\subsection{Comparison of the reduced model and the full model of \cite{scott_spatial_2021}}\label{app:comparison}
We verify the reduced model \eqref{reduced_RhoA_model}  captures  the results obtained in \cite{scott_spatial_2021} for the full model for the RhoA signalling pathway. In numerical simulations of model~\eqref{reduced_RhoA_model} we use the same parameter values as in \cite{scott_spatial_2021}, except for the diffusion coefficients $D_1$ and $D_2$ for activated and deactivated FAK. It is  suggested in the literature that $D_1 = D_2= 4~\mu$m$^2/s$ \cite{le_devedec_residence_2012}, but \cite{scott_spatial_2021} uses $10~\mu$m$^2/s$ due to computational issues. Thus in our numerical simulations we consider  both diffusion coefficients. 

\begin{table}[!h]
    \centering
    \begin{tabular}{|l|l|l|l|}
        $\phi_d^0=0.7\mu\mbox{mol}/\mbox{dm}^3$ & $C=3.25$~kPa & $D_1=4\mbox{ or }10~\mu\mbox{m}^2/$s & $k_2=0.015~\mbox{s}^{-1}$  \\
        $\phi_a^0=0.3\mu\mbox{mol}/\mbox{dm}^3$ & $\gamma=8.8068~\mbox{dm}^3/\mu\mbox{mol}$ & $D_2=4\mbox{ or }10~\mu\mbox{m}^2/$s & $k_3=0.379~\mbox{s}^{-1}$  \\
        $\rho_a^0=33.6~\# /\mu\mbox{m}^2$ & $n=5$ & $D_3=0.3~\mu\mbox{m}^2/$s & $k_4=0.625~\mbox{s}^{-1}$  \\
        \quad $\approx 6\cdot 10^{-7}~\mu\mbox{mol}/\mbox{dm}^2$ &  $E=0.1,5.7,7\cdot 10^6$~kPa & $k_1=0.035~\mbox{s}^{-1}$ & $k_5=0.0168~\mbox{s}^{-1}$\\
        $\rho_d^0=1~\mu\mbox{mol}/\mbox{dm}^3$ 
        & $|Y|=1193~\mu$m$^3$ &
        $|\Gamma|=1020~\mu$m$^2$ &
    \end{tabular}
    \caption{Parameter values for  model~\eqref{reduced_RhoA_model}.}
    \label{tab:parameters_app}
\end{table}
The model \eqref{reduced_RhoA_model} is implemented in  FEniCS \cite{logg_automated_2012}, using a Finite Element Method for discretization in space and IMEX time-stepping method to discretize in time,  see Section A.3 for more details. Considering domain $Y\subset \mathbb{R}^3$, denoting the cytoplasm,  and $\Gamma=\partial Y$, defining the cell membrane, and times interval $(0,T)$, with $T = 100$~s, for the space discretisation we choose meshsize $h=2.94$ and time step $\Delta t=0.5$ for the backwards Euler discretisation in time. For our domain we have $n_r = |Y|/|\Gamma|= 1.17~\mu$m. We consider three different stimuli, similar to~\cite{scott_spatial_2021}, (i) the `$2$D stimulus', where the 
substrate stiffness is only applied to the bottom of the cell and any reaction terms of RhoA are nonzero only at the bottom of the cell, (ii) the `$2$xD stimulus', where the 
substrate stiffness is only applied to the bottom of the cell but the reaction terms of RhoA are nonzero on the whole cell membrane, and  (iii) the `$3$D stimulus' where the cell is embedded in an agar (substrate) and the impact of the substrate stiffness on the signalling processes is considered on the whole cell membrane.

\begin{figure}[!h]
\begin{minipage}{0.05\textwidth}\raggedleft
    \vspace{1cm}
    $\phi_a$ \\ \vspace{1cm}
    $\rho_a$
\end{minipage}\hfill
\begin{minipage}{0.29\textwidth} \centering
    2D stimulus \vspace{0.5cm} \\
\begin{minipage}{0.33\textwidth} \centering
    {\footnotesize $0.1$kPa } \\ 
        \includegraphics[width=\textwidth]{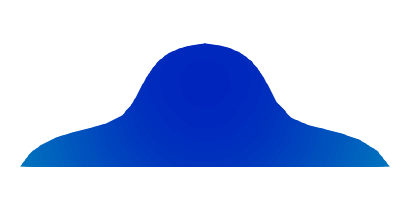}\\
        \vspace{0.5cm}
        \includegraphics[width=\textwidth]{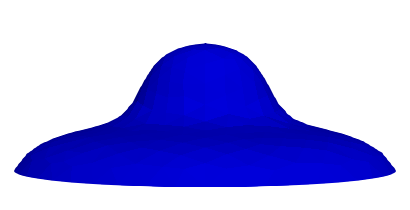}
\end{minipage}\hfill
\begin{minipage}{0.33\textwidth} \centering
    {\footnotesize $5.7$kPa } \\
        \includegraphics[width=\textwidth]{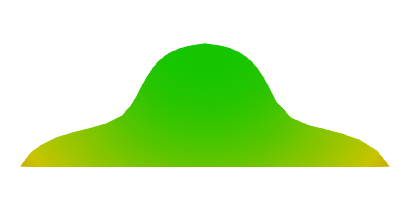}\\
        \vspace{0.5cm}
        \includegraphics[width=\textwidth]{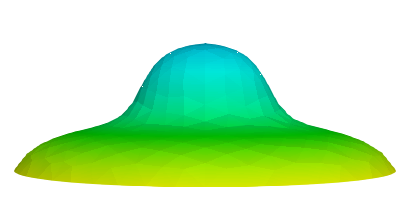}
\end{minipage}\hfill
\begin{minipage}{0.33\textwidth} \centering
    {\footnotesize $7$GPa} \\
        \includegraphics[width=\textwidth]{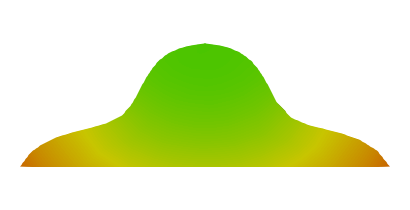}\\
        \vspace{0.5cm}
        \includegraphics[width=\textwidth]{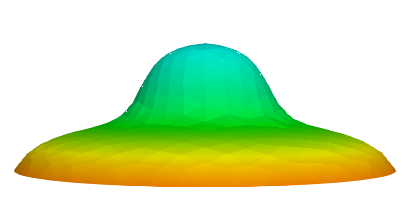}
\end{minipage}\hfill 
\end{minipage}\hfill \hspace{0.2cm}
\begin{minipage}{0.29\textwidth}  \centering
    2xD stimulus \vspace{0.5cm} \\
\begin{minipage}{0.33\textwidth} \centering
    {\footnotesize $0.1$kPa } \\ 
        \includegraphics[width=\textwidth]{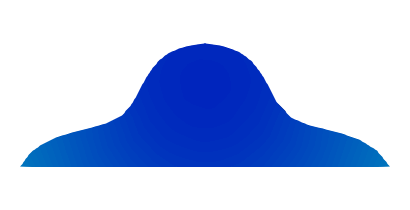}\\
        \vspace{0.5cm}
        \includegraphics[width=\textwidth]{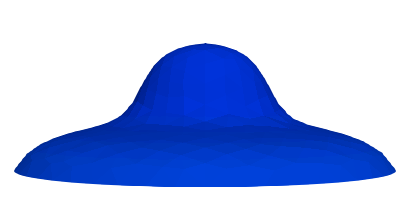}
\end{minipage}\hfill
\begin{minipage}{0.33\textwidth} \centering
    {\footnotesize $5.7$kPa } \\
        \includegraphics[width=\textwidth]{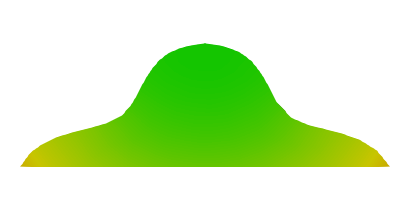}\\
        \vspace{0.5cm}
        \includegraphics[width=\textwidth]{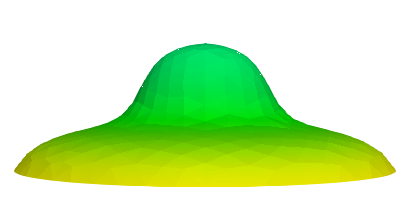}
\end{minipage}\hfill
\begin{minipage}{0.33\textwidth} \centering
    {\footnotesize $7$GPa} \\
        \includegraphics[width=\textwidth]{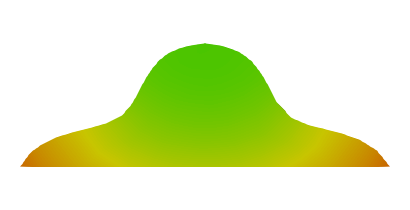}\\
        \vspace{0.5cm}
        \includegraphics[width=\textwidth]{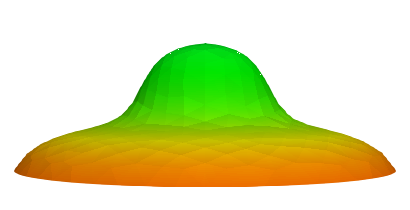}
\end{minipage}\hfill
\end{minipage}\hfill \hspace{0.2cm}
\begin{minipage}{0.29\textwidth}  \centering
    3D stimulus \vspace{0.5cm} \\
\begin{minipage}{0.33\textwidth} \centering
    {\footnotesize $0.1$kPa } \\ 
        \includegraphics[width=\textwidth]{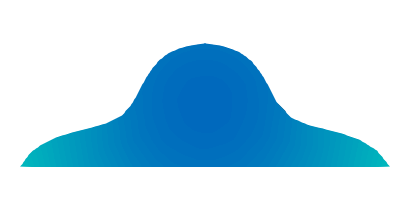}\\
        \vspace{0.5cm}
        \includegraphics[width=\textwidth]{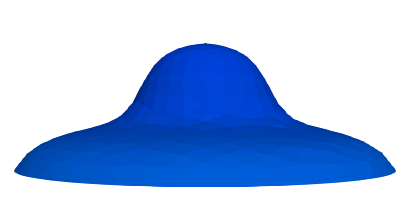}
\end{minipage}\hfill
\begin{minipage}{0.33\textwidth} \centering
    {\footnotesize $5.7$kPa } \\
        \includegraphics[width=\textwidth]{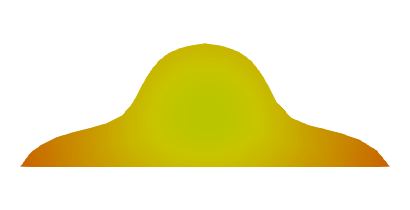}\\
        \vspace{0.5cm}
        \includegraphics[width=\textwidth]{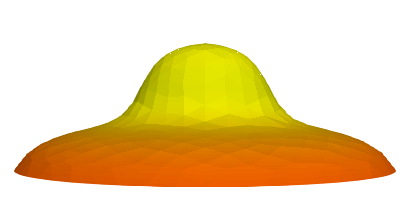}
\end{minipage}\hfill
\begin{minipage}{0.33\textwidth} \centering
    {\footnotesize $7$GPa} \\
        \includegraphics[width=\textwidth]{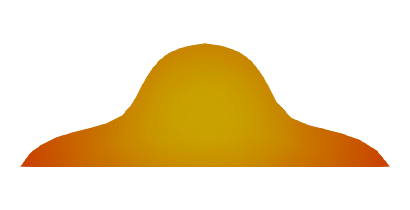}\\
        \vspace{0.5cm}
        \includegraphics[width=\textwidth]{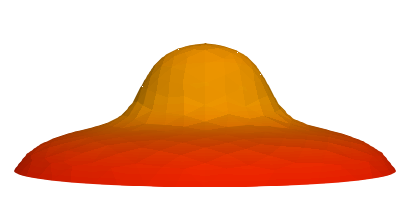}
\end{minipage}\hfill
\end{minipage}\hfill
\begin{minipage}{0.025\textwidth}\raggedright
    \vspace{1cm}
    \includegraphics[width=\textwidth]{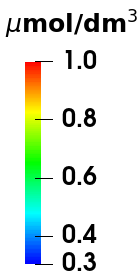}\\
    \vspace{0.1cm}
    \includegraphics[width=\textwidth]{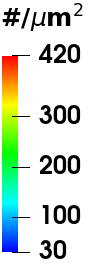}
\end{minipage}
    \caption{Numerical simulation results showing $\phi_a$ and $\rho_a$ for reduced model~\eqref{reduced_RhoA_model} at steady state for~$T=100$~s. Parameter values  as in Table~\ref{tab:parameters_app} and $D_1 = D_2 = 4~\mu\text{m}^2/\text{s}$.}
    \label{fig:rhomodel}
\end{figure}

\begin{figure}[!h]
\begin{minipage}{0.5\textwidth}
        \includegraphics[width=\textwidth]{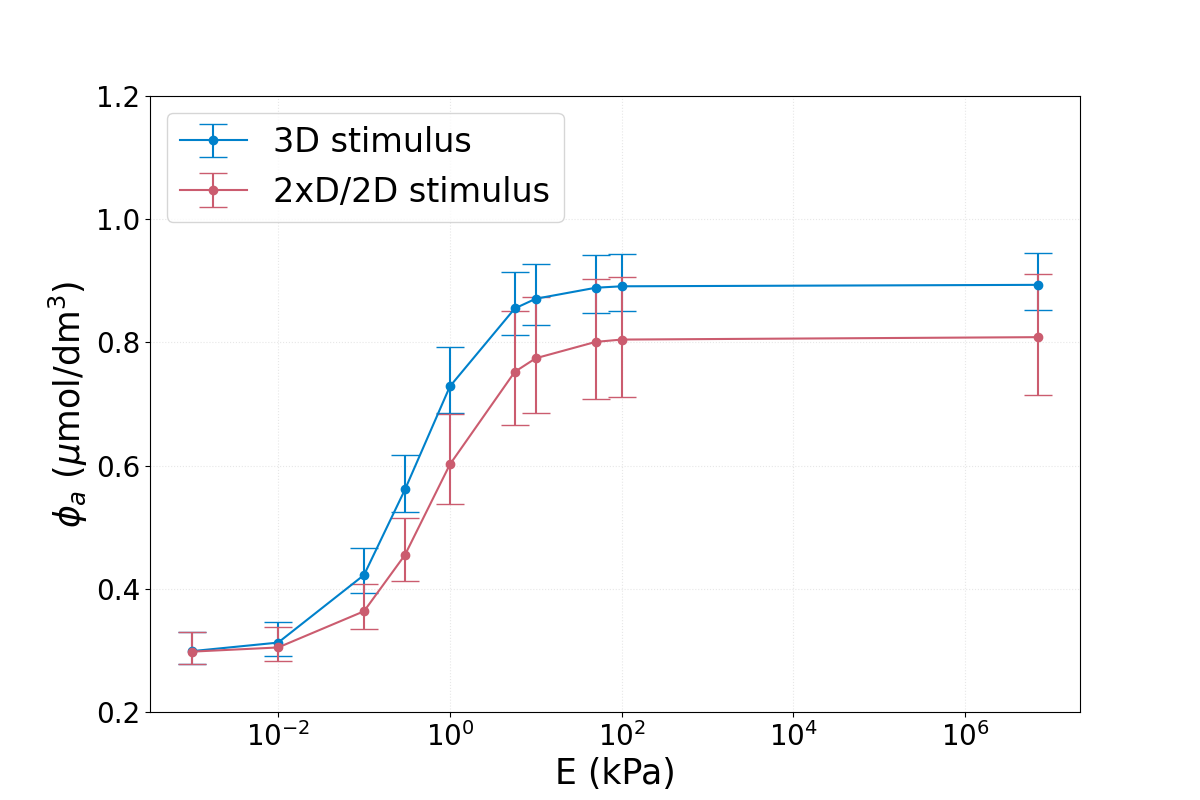}
\end{minipage}\hfill
\begin{minipage}{0.5\textwidth}
        \includegraphics[width=\textwidth]{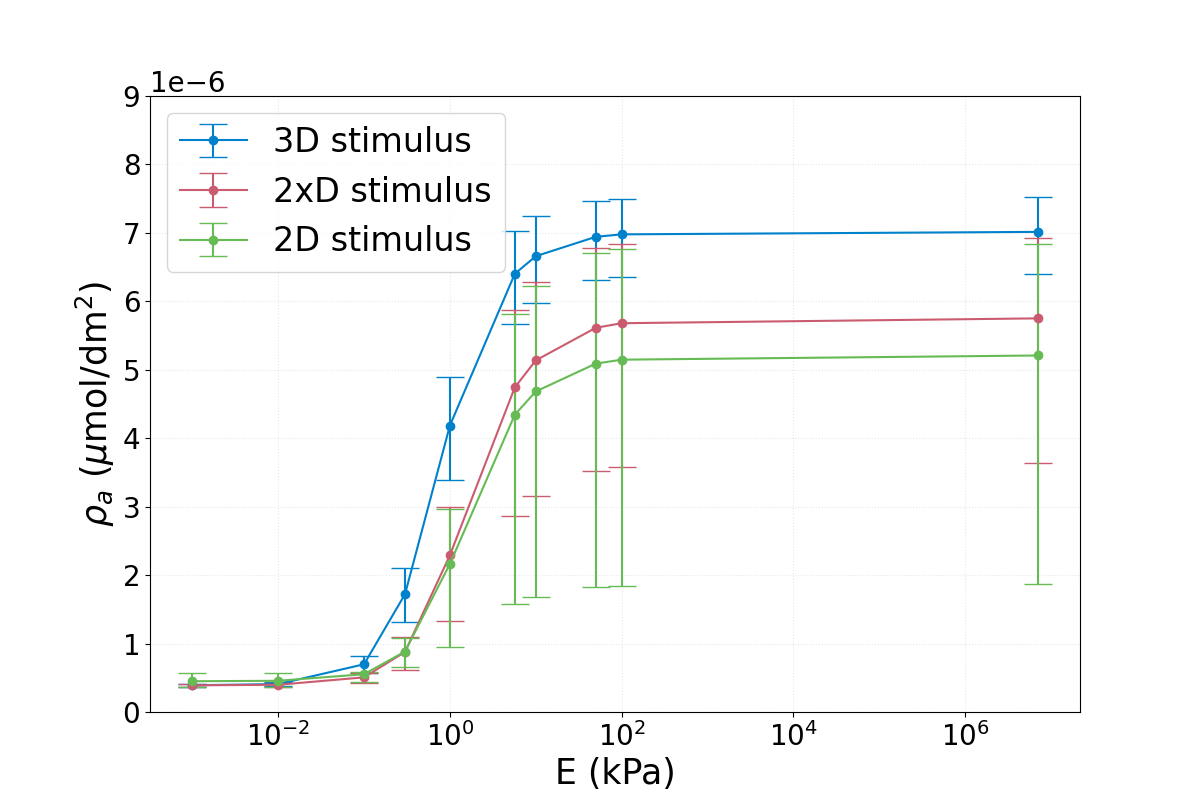}
\end{minipage}
    \caption{Results showing the effect of substrate stiffness $E$ on $\phi_a$ and $\rho_a$ for the reduced model~\eqref{reduced_RhoA_model} at $T=100$~s by which time the results are at a steady state. Parameter values as in Table~\ref{tab:parameters_app} with $D_1 = D_2 = 4~\mu\text{m}^2/\text{s}$. }
    \label{fig:rhomodel2}
\end{figure}

Comparing the simulation results for reduced model~\eqref{reduced_RhoA_model} in Figs~\ref{fig:rhomodel} and~\ref{fig:rhomodel2} to  the results presented in \cite[Fig~3]{scott_spatial_2021} for the full model, the dynamics of FAK and RhoA are almost identical qualitatively. Similar to the results in \cite[Fig~3]{scott_spatial_2021}, the highest concentration of both $\phi_a$ and $\rho_a$ is at the edges of the cell. Also, there is threshold value of $E\approx 1$kPa, below which   the concentrations of $\phi_a$ and $\rho_a$  stay close to the initial values  and then reaches high  steady states values, similar  for both $E=5.7$kPa and $E=7$GPa. Quantitatively, the values for $\phi_a$ are also close to the one reported in \cite[Fig~3]{scott_spatial_2021}. This suggests that the reduction of the  model as well as considering the whole cell domain without excluding a nucleus  does not have  significant effect on the dynamics of  FAK. However, for $\rho_a$ we obtain  slightly lower concentration, where the maximum concentration in our results is $420\#/\mu m^2$ and the maximum concentration in~\cite[Fig.3]{scott_spatial_2021} is $593\#/\mu m^2$. This difference could be related to the model reduction and approximation for the deactivated RhoA. The differences in the conversion from $\mu \mbox{mol}/dm^2$ to $\# /\mu m^2$ using our approach, see Section A.4 for details, and the approach of \cite{scott_spatial_2021} may also contribute.

\begin{figure}[!h]
\begin{minipage}{0.05\textwidth}\raggedleft
    \vspace{1cm}
    $\phi_a$ \\ \vspace{1cm}
    $\rho_a$
\end{minipage}\hfill
\begin{minipage}{0.29\textwidth} \centering
    2D stimulus \vspace{0.5cm} \\
\begin{minipage}{0.33\textwidth} \centering
    {\footnotesize $0.1$kPa } \\ 
        \includegraphics[width=\textwidth]{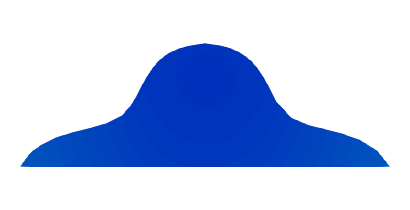}\\
        \vspace{0.5cm}
        \includegraphics[width=\textwidth]{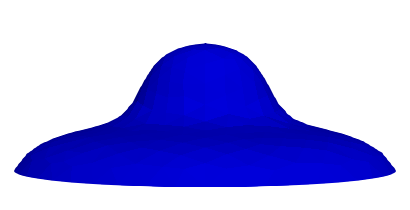}
\end{minipage}\hfill
\begin{minipage}{0.33\textwidth} \centering
    {\footnotesize $5.7$kPa } \\
        \includegraphics[width=\textwidth]{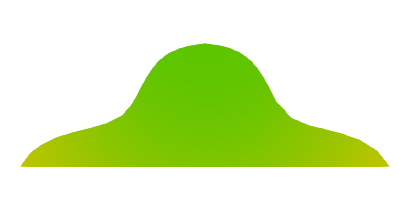}\\
        \vspace{0.5cm}
        \includegraphics[width=\textwidth]{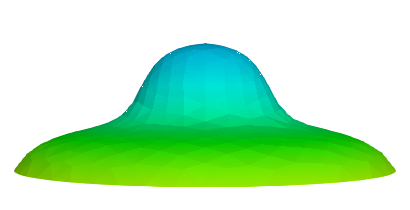}
\end{minipage}\hfill
\begin{minipage}{0.33\textwidth} \centering
    {\footnotesize $7$GPa} \\
        \includegraphics[width=\textwidth]{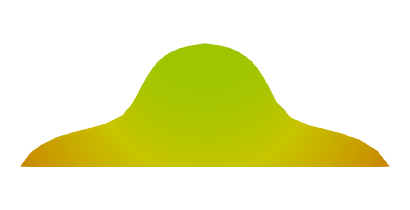}\\
        \vspace{0.5cm}
        \includegraphics[width=\textwidth]{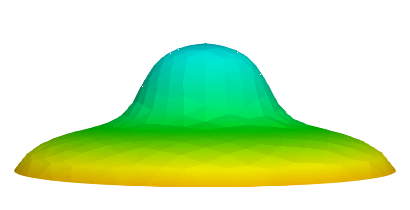}
\end{minipage}\hfill 
\end{minipage}\hfill \hspace{0.2cm}
\begin{minipage}{0.29\textwidth}  \centering
    2xD stimulus \vspace{0.5cm} \\
\begin{minipage}{0.33\textwidth} \centering
    {\footnotesize $0.1$kPa } \\ 
        \includegraphics[width=\textwidth]{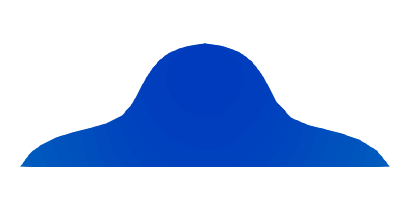}\\
        \vspace{0.5cm}
        \includegraphics[width=\textwidth]{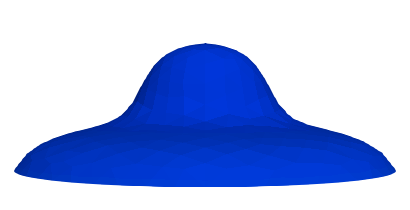}
\end{minipage}\hfill
\begin{minipage}{0.33\textwidth} \centering
    {\footnotesize $5.7$kPa } \\
        \includegraphics[width=\textwidth]{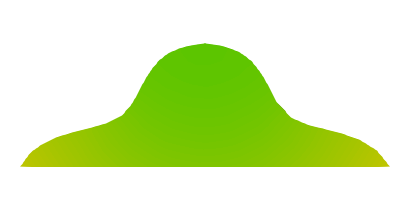}\\
        \vspace{0.5cm}
        \includegraphics[width=\textwidth]{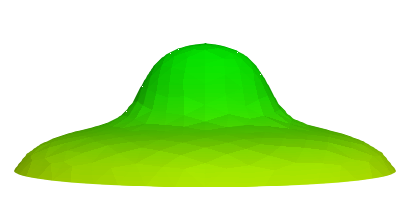}
\end{minipage}\hfill
\begin{minipage}{0.33\textwidth} \centering
    {\footnotesize $7$GPa} \\
        \includegraphics[width=\textwidth]{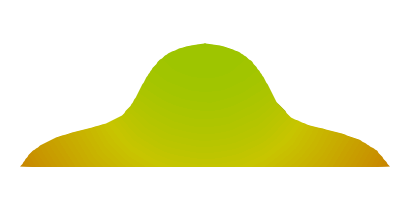}\\
        \vspace{0.5cm}
        \includegraphics[width=\textwidth]{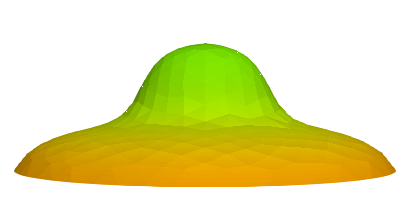}
\end{minipage}\hfill
\end{minipage}\hfill \hspace{0.2cm}
\begin{minipage}{0.29\textwidth}  \centering
    3D stimulus \vspace{0.5cm} \\
\begin{minipage}{0.33\textwidth} \centering
    {\footnotesize $0.1$kPa } \\ 
        \includegraphics[width=\textwidth]{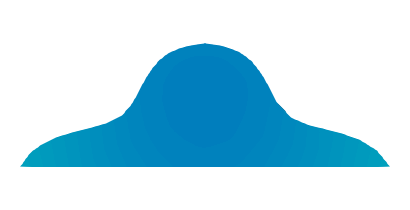}\\
        \vspace{0.5cm}
        \includegraphics[width=\textwidth]{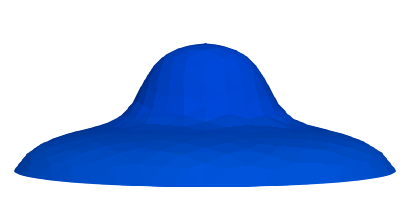}
\end{minipage}\hfill
\begin{minipage}{0.33\textwidth} \centering
    {\footnotesize $5.7$kPa } \\
        \includegraphics[width=\textwidth]{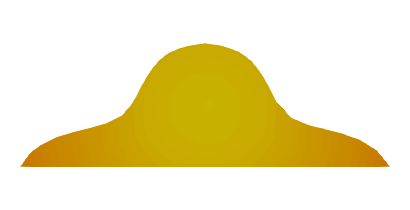}\\
        \vspace{0.5cm}
        \includegraphics[width=\textwidth]{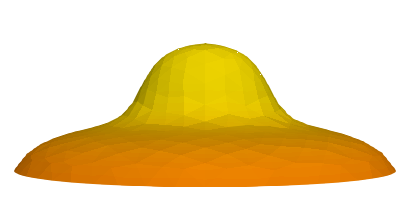}
\end{minipage}\hfill
\begin{minipage}{0.33\textwidth} \centering
    {\footnotesize $7$GPa} \\
        \includegraphics[width=\textwidth]{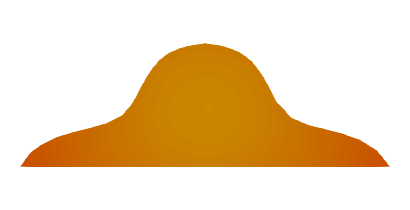}\\
        \vspace{0.5cm}
        \includegraphics[width=\textwidth]{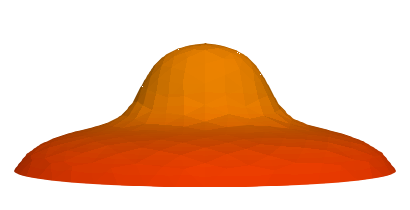}
\end{minipage}\hfill
\end{minipage}\hfill
\begin{minipage}{0.025\textwidth}\raggedright
    \vspace{1cm}
    \includegraphics[width=\textwidth]{figures/coupled_ca_bar.png}\\
    \vspace{0.1cm}
    \includegraphics[width=\textwidth]{figures/coupled_p_bar.png}
\end{minipage}
    \caption{Numerical simulation results showing $\phi_a$ and $\rho_a$ for reduced model~\eqref{reduced_RhoA_model} at $T=100$~s by which time the results are at a steady state. Parameter values as in Table~\ref{tab:parameters_app} with $D_1=D_2=10 \mu$m$/s$.}
    \label{fig:rhomodel_D1=10}
\end{figure}

\begin{figure}[!h]
\begin{minipage}{0.5\textwidth}
        \includegraphics[width=\textwidth]{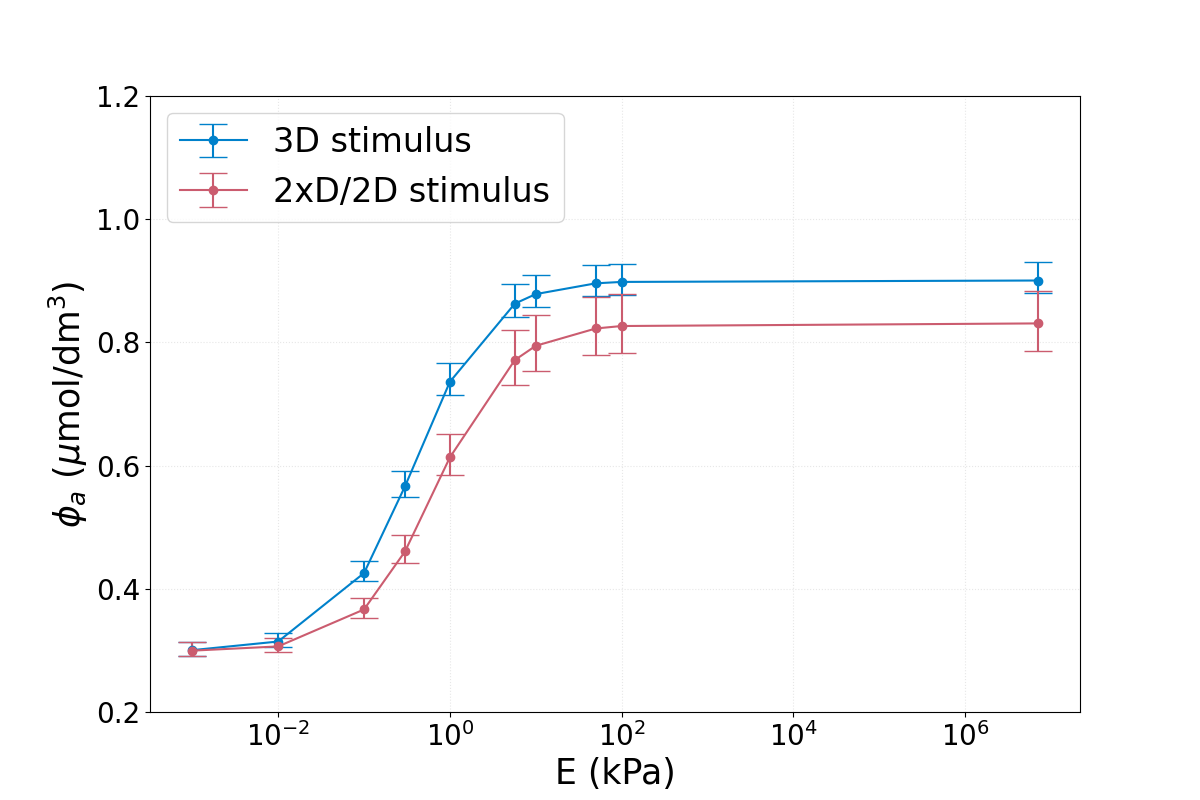}
\end{minipage}\hfill
\begin{minipage}{0.5\textwidth}
        \includegraphics[width=\textwidth]{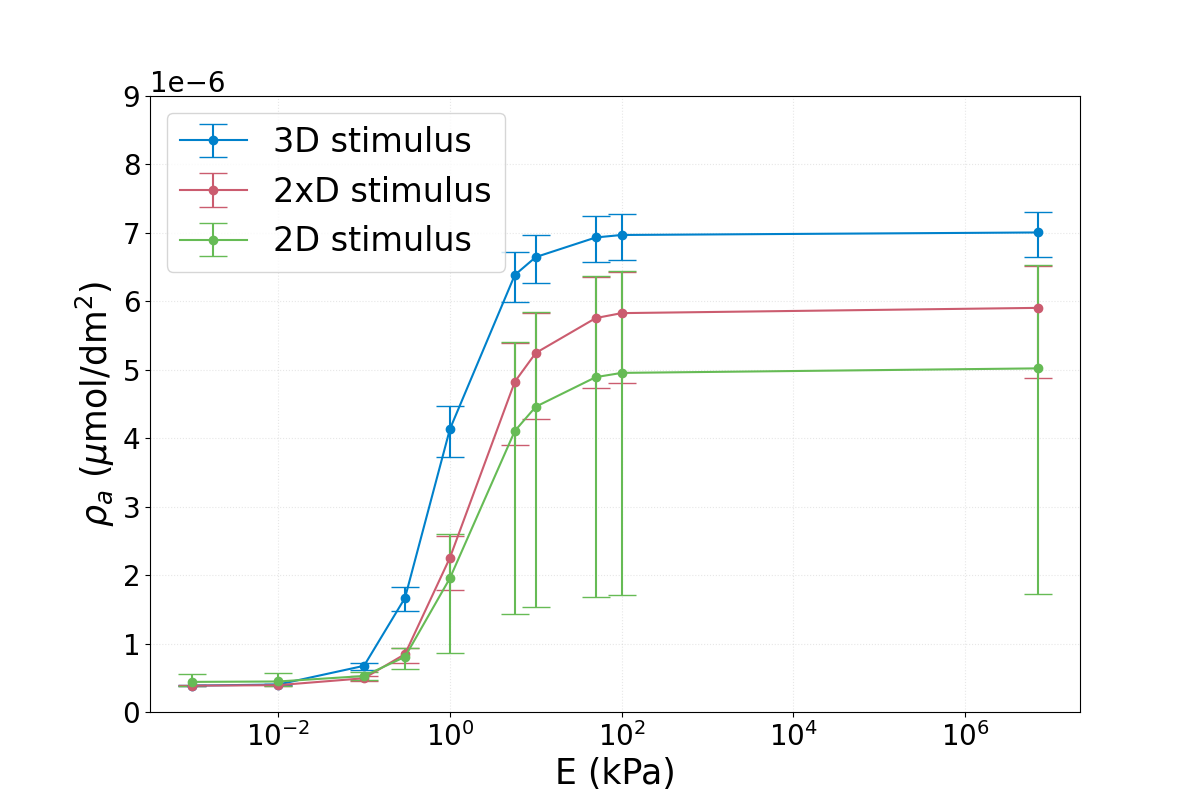}
\end{minipage}
    \caption{Results showing effect of substrate stiffness $E$ on $\phi_a$ and $\rho_a$ for the reduced model~\eqref{reduced_RhoA_model} at $T=100$~s by which time the results are at a steady state. Parameter values as in Table~\ref{tab:parameters_app} with $D_1=D_2=10 \mu$m$/s$.}
    \label{fig:rhomodel2_D1=10}
\end{figure}

\newpage
Simulation results for the reduced model with diffusion coefficients $D_1=D_2=10 \mu\mbox{m}^2/$s are presented in Figs~\ref{fig:rhomodel_D1=10} and~\ref{fig:rhomodel2_D1=10}. Comparing Figs~\ref{fig:rhomodel2} and~\ref{fig:rhomodel2_D1=10}, we see that the averaged over space dynamics are very similar, but numerical simulation results for a lower diffusion coefficient show lower minimum and higher maximum concentrations, which can be explained by the fact that slower diffusion of $\phi_a$  causes stronger heterogeneity across the cell domain. In Figs~\ref{fig:rhomodel} and~\ref{fig:rhomodel_D1=10}, we observe the same dynamics for both diffusion coefficients, where the maximum concentration is at the edges of the cell and the minimum concentration is in the middle. 
Since the dynamics of $\rho_a$ depends on $\phi_a$, similar results are obtained for $\rho_a$. 

\subsection{Temporal statistics}\label{app:temp_stats}
\begin{figure}[!h]
\begin{minipage}{\textwidth} 
    \hspace{4cm} $E=0.1$~kPa \hspace{6cm} $E=5.7$~kPa \\ 
    \vspace{-0.2cm}
    \hspace{1.5cm} {\footnotesize $E_c=0.6$~kPa } \hspace{2.3cm}  {\footnotesize $E_c=f(\phi_a)$ } \hspace{2.3cm} {\footnotesize $E_c=0.6$~kPa } \hspace{2.3cm} {\footnotesize $E_c=f(\phi_a)$ } \\ 
    \vspace{-0.1cm}
\end{minipage}\hfill
\begin{minipage}{\textwidth} 
    \includegraphics[width=\textwidth]{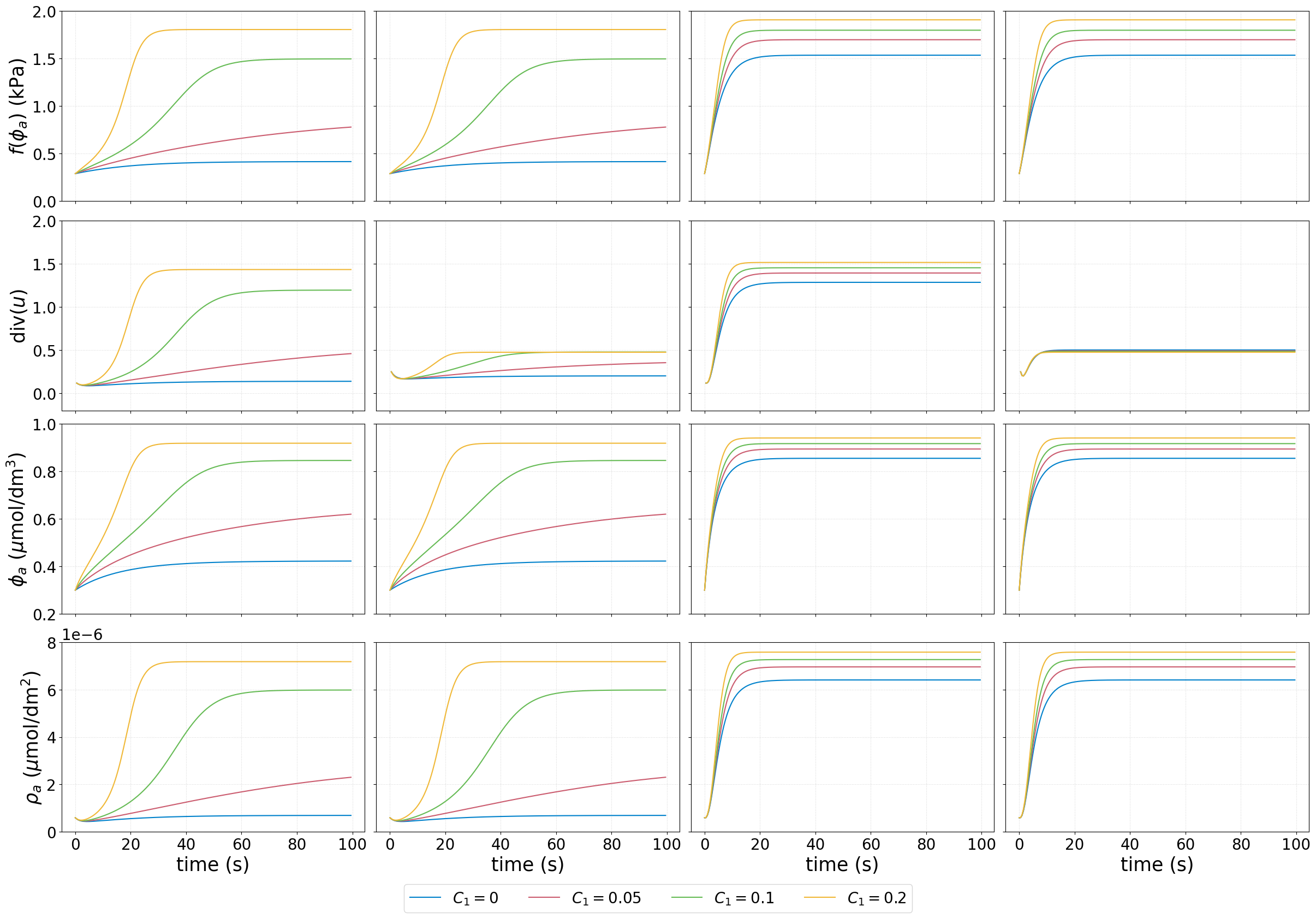}
\end{minipage}
    \caption{Simulations results showing the mean of  $f(\phi_a)$, ${\rm div}(u)$, $\phi_a$ and $\rho_a$ as  functions of time, in the case of model \eqref{eq:coupled_elast}, \eqref{eq:coupled_force_boundary} and \eqref{eq:coupled_reactions}, axisymmetric shape and 3D stimulus, for different couplings, four different values for $C_1$, and two different values for the substrate stiffness $E$. All other parameter values as in Table~\ref{tab:parameters}. The corresponding results can be found in Fig~\ref{fig:sim_subs}.}
    \label{fig:temp_subs}
\end{figure}

Fig~\ref{fig:temp_subs} shows the evolution of the mean of $f(\phi_a)$, ${\rm div}(u)$, $\phi_a$ and $\rho_a$ over time for different couplings and parameters for the axisymmetric cell shape.

\subsection{Numerical Scheme}\label{app:numeric}
For numerical simulations of model \eqref{eq:coupled_elast}-\eqref{eq:coupled_reactions}  we use FEM for discretization in space and backward Euler for discretization in time,  implemented in FEniCS~\cite{logg_automated_2012}. 
Consider the space
\begin{equation}\label{eq:W_def}
    \mathcal{W}(Y) = \lbrace u\in H^1(Y): \int_Y u_i \dd{x}=0, \; \int_Y \left(\partial_{x_j}u_i-\partial_{x_i} u_j \right)\dd{x} = 0 \; \;  \mbox{ for } \; i,j=1,2,3 \rbrace,
\end{equation}
such that $\mathcal{W}(Y)\cap \mathcal{R}(Y)=0$ with $\mathcal{R}(Y)$ the space of rigid motions. Then, the weak solution of the model \eqref{eq:coupled_elast}-\eqref{eq:coupled_reactions} is given  by  $(\phi_d,\phi_a) \in L^2(0,T;H^1(Y))$, $\rho_a\in L^2(0,T;H^1(\Gamma))$, with $(\partial_t\phi_d,\partial_t\phi_a)\in L^2(0,T;H^1(Y)')$ and $\partial_t \rho_a\in L^2(0,T;H^1(\Gamma)')$, and $u \in L^2(0,T;\mathcal{W}(Y))$ satisfying  
\begin{equation} \label{eq:weak_formul}
\begin{aligned}
   & \langle\partial_t \phi_d, \psi \rangle_{(H^1)^\prime,T} + \langle D_1 \nabla \phi_d, \nabla \psi \rangle_{Y_T} + \langle C_1\tr(\sigma(u))_+\phi_d, \psi \rangle_{Y_T} + \langle n_r\tilde{k}_3\phi_d, \psi\rangle_{\Gamma_T} = \langle k_1\phi_a,\psi \rangle_{Y_T}, \\
   & \langle \partial_t \phi_a, \varphi \rangle_{(H^1)^\prime,T} + \langle D_2 \nabla \phi_a, \nabla \varphi \rangle_{Y_T} + \langle k_1\phi_a,\varphi \rangle_{Y_T} =  \langle C_1\tr(\sigma(u))_+\phi_d, \varphi \rangle_{Y_T} + \langle n_r\tilde{k}_3\phi_d, \psi\rangle_{\Gamma_T} ,\\
  &  \langle\partial_t \rho_a, w \rangle_{(H^1)^\prime,T}+ \langle D_3 \nabla_\Gamma \rho_a, \nabla_\Gamma w \rangle_{\Gamma_T} + \langle \tilde{k}_4(\phi_a)\rho_a,w \rangle_{\Gamma_T}= \langle n_r\tilde{k}_5(\phi_a),w \rangle_{\Gamma_T},   \\
  &  \langle E(\phi_a)\epsilon(u),\epsilon(v) \rangle_{Y_T} = \langle k_6\mathbb{P}(\rho_a\nu),v \rangle_{\Gamma_T}, 
\end{aligned}
\end{equation}
for all $\psi,\varphi \in L^2\left(0,T;H^1(Y)\right)$, $w \in L^2\left(0,T;H^1(\Gamma)\right)$, and $v\in L^2(0,T;H^1(Y))$,  with initial conditions satisfied in the $L^2$-sense, and where $\tilde{k}_3=k_2+k_3\frac{E}{C+E}$, $\tilde{k}_4(\phi_a)=k_4+k_5((\gamma \phi_a)^n+1)$ and  $\tilde{k}_5(\phi_a)=k_5((\gamma \phi_a)^n+1)\frac{M_\rho}{|Y|}$.
Here $\langle\phi, \psi \rangle_{(H^1)^\prime,T}$ denotes the dual product between $\phi \in L^2(0,T; H^1(Y))$ and $\psi \in L^2(0,T; H^1(Y)^\prime)$ or between  $\phi \in L^2(0,T; H^1(\Gamma))$ and $\psi \in L^2(0,T; H^1(\Gamma)^\prime)$ and 
$$
\langle \phi, \psi\rangle_{Y_T} = \int_0^T \int_Y \phi \psi dx dt, \quad \langle \phi_1, \psi_1 \rangle_{\Gamma_T} = \int_0^T \int_\Gamma \phi_1 \psi_1 dx dt, \quad \text{ where } 1/p_1 + 1/p_2 =1, \; \; 1/q_1+ 1/q_2 =1, 
$$
for  $\phi \in L^{p_1}(0,T; L^{q_1}(Y))$, $\psi \in L^{p_2}(0,T; L^{q_2}(Y))$, $\phi_1 \in L^{p_1}(0,T; L^{q_1}(\Gamma))$, and  $\psi_1 \in L^{p_2}(0,T; L^{q_2}(\Gamma))
$. 
The discretization of the domain is given  by the polyhedral approximation of $Y$ such that $Y_h$ is the union of finitely many tetrahedrons in $\mathbb{R}^3$, and $S_h$ is  the set of these tetrahedrons $K$, such that
$$
Y_h = \bigcup_{K \in S_h} K.
$$
Then  the surface $\Gamma$ is approximated by $\Gamma_h$ such that $\Gamma_h=\partial Y_h$. The mesh size is defined by the maximum diameter of a simplex $h=\max\lbrace h_Y,h_\Gamma\rbrace$, where $ h_Y = \max_{K\in S_h}h_Y(K)$ and  $h_\Gamma = \max_{R=K\cap \Gamma_h\neq 0, K \in S_h}h_\Gamma(R)$ with $h_Y(K)$ being the diameter of a tetrahedron and $h_\Gamma(R)$  the diameter of a triangle on the surface.
\begin{figure}[!h]
    \begin{center}
    \begin{subfigure}{.3\textwidth}
    \includegraphics[width=\textwidth]{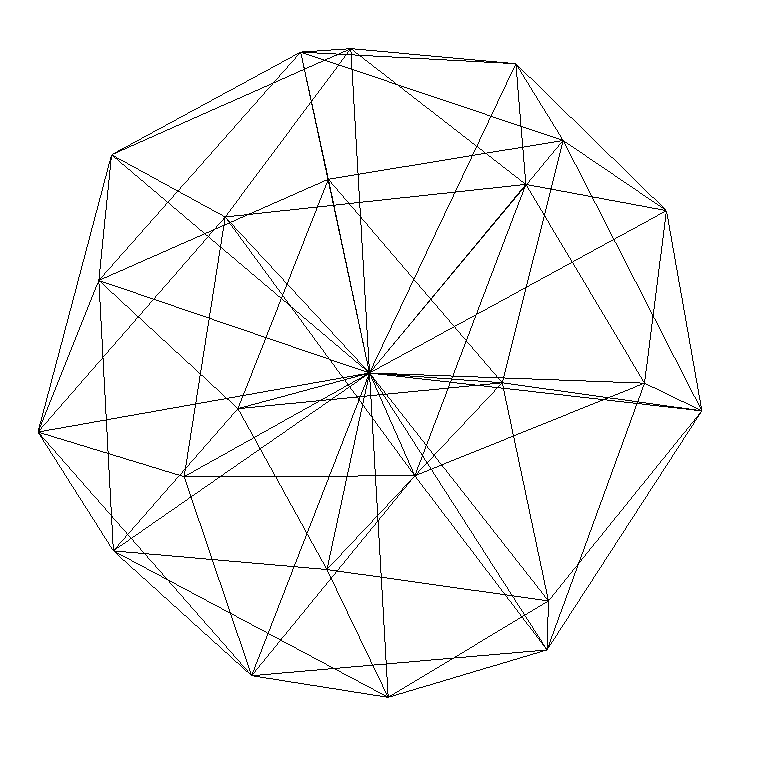}
    \end{subfigure}%
    \begin{subfigure}{.3\textwidth}
    \includegraphics[width=\textwidth]{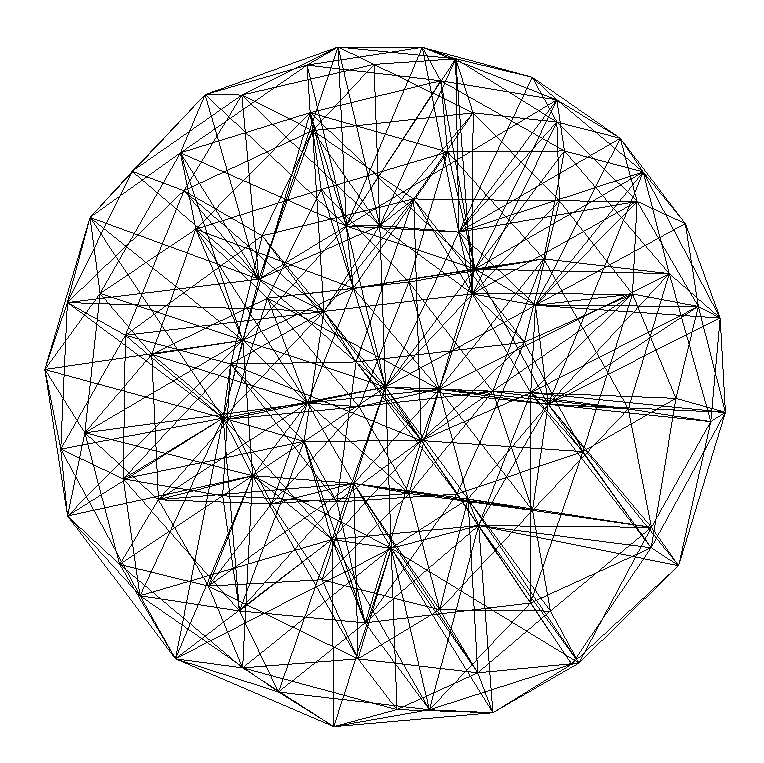}
    \end{subfigure}%
    \begin{subfigure}{.3\textwidth}
    \includegraphics[width=\textwidth]{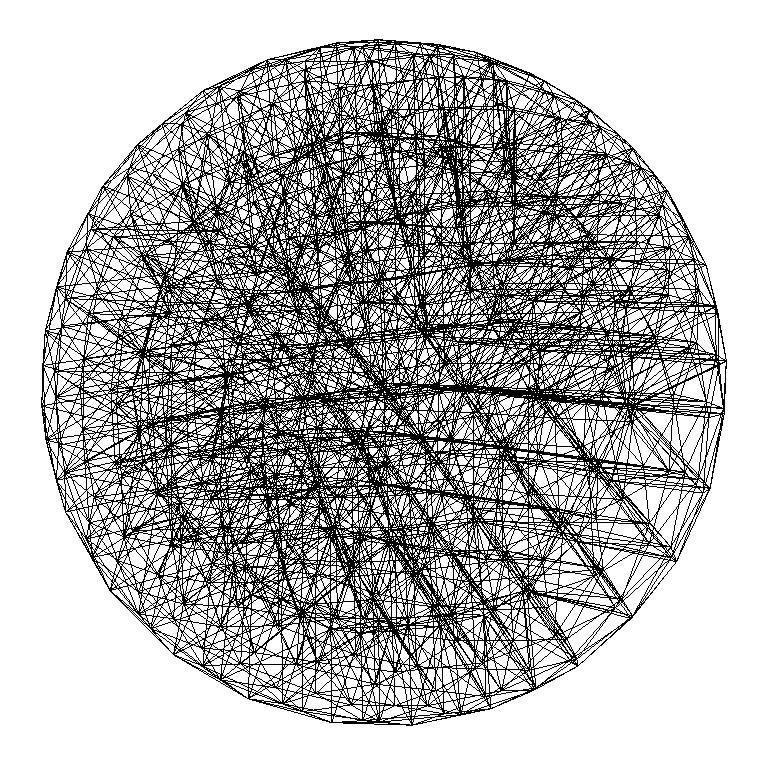}
    \end{subfigure}
    \end{center}
    \caption{Tetrahedral approximations of a unit sphere with, from left to right, decreasing values of $h$, created with Gmsh~\cite{geuzaine_gmsh_2009}.}
    \label{fig:mesh}
\end{figure}

\noindent The bulk and surface finite element spaces are given by 
\begin{align*}
    \mathbb{V}_{h,Y} &= \left\lbrace \Psi \in C(Y_h):\left.\Psi\right|_K \mbox{ is linear affine for each } K \in S_h \right\rbrace,\\
    \mathbb{V}_{h,\Gamma} &= \left\lbrace \Psi \in C(\Gamma_h):\left.\Psi\right|_R \mbox{ is linear affine for each } K \in S_h \mbox{ with } R=K\cap \Gamma_h\neq 0 \right\rbrace,\\
    \mathbb{W}_{h,Y} &= \left\lbrace V\in\mathcal{W}(Y_h) \mbox{ and } V_i \in C(Y_h):\left.V_i\right|_K \mbox{ is linear affine for each } K \in S_h \mbox{ for }i=1,2,3 \right\rbrace,
\end{align*}
where we recall the definition of $\mathcal{W}$ in Eq~(\ref{eq:W_def}). 
The bulk space is spanned by nodal basis functions defined by
$$
\chi_j\in\mathbb{V}_{h,Y}, \quad \chi_j(X_k)=\delta_{jk} \quad \text{ for } \; j,k=1,...,J,
$$
where $\delta_{jk}$ is the Kronecker delta and  $J$ is the number of nodes (vertices)  $X_j\in Y$ of the tetrahedrons $S_h$. Then  function $\Phi(t)\in \mathbb{V}_{h,Y}$ has the form
$$
\Phi(t, x)=\sum_{j=1}^J a_j(t)\chi_j(x) \quad  \text{ for } \; x\in Y_h, \; t \in (0,T),
$$
with real measurable functions $a_j$. Similarly  the surface finite element space is spanned by nodal basis functions 
$$
\mu_j\in\mathbb{V}_{h,\Gamma}, \quad \mu_j(Z_k)=\delta_{jk} \quad \text{ for } \; j,k=1,...,M,
$$
where $Z_j\in\Gamma$, with $j,k=1,...,M$,  are nodes of the triangulated surface such that $Z_k=X_k\bigcap \Gamma_h\neq 0$. Then  function $P(t)\in \mathbb{V}_{h,\Gamma}$ has the form
$$
P(t,x)=\sum_{j=1}^M b_j(t)\mu_j(x) \quad \text{ for } \; x\in \Gamma_h, \; t \in (0,T),
$$
with real measurable functions $b_j$. Thus, the semi-discretized problem corresponding to \eqref{eq:coupled_elast}-\eqref{eq:coupled_reactions} reads
\begin{equation} \label{eq:space_discr}
\begin{aligned}
    &\left<\partial_t \Phi_d, \Psi \right>_{Y_h} + \left< D_1 \nabla \Phi_d, \nabla \Psi \right>_{Y_h} + \left<C_1\tr(\sigma(U))_+\Phi_d, \Psi \right>_{Y_h} + \left<n_r\tilde{k}_3\Phi_d, \Psi\right>_{\Gamma_h} = \left<k_1\Phi_a,\Psi \right>_{Y_h}, \\
   & \left<\partial_t \Phi_a, \Psi \right>_{Y_h} + \left< D_2 \nabla \Phi_a, \nabla \Psi \right>_{Y_h} + \left<k_1\Phi_a,\Psi \right>_{Y_h} =  \left<C_1\tr(\sigma(U))_+\Phi_d, \Psi \right>_{Y_h} + \left<n_r\tilde{k}_3\Phi_d, \Psi\right>_{\Gamma_h} , \\
   & \left<\partial_t P_a, W \right>_{\Gamma_h} + \left< D_3 \nabla_{\Gamma_h} P_a, \nabla_{\Gamma_h} W \right>_{\Gamma_h} + \left<\tilde{k}_4(\Phi_a)P_a,W \right>_{\Gamma_h} = \left<n_r\tilde{k}_5(\Phi_a),W \right>_{\Gamma_h},\\
  &  \left<E(\Phi_a)\epsilon(U),\epsilon(V) \right>_{Y_h} = \left<k_6\mathbb{P}(P_a\hat{\nu}),V \right>_{\Gamma_h},
\end{aligned}
\end{equation}
for every test function $\Psi\in \mathbb{V}_{h,Y}$, $W\in \mathbb{V}_{h,\Gamma}$ and $V\in \mathbb{W}_{h,Y}$.

To obtain the fully discrete problem we discretize \eqref{eq:space_discr} in time  using the  backwards Euler  method with 
$$
\partial_t \Phi \approx \frac{\Phi^n-\Phi^{n-1}}{\Delta t},
$$
where $\Delta t = T/N$, and  an IMEX time-stepping method, in which the diffusion terms are treated implicitly and the nonlinear reaction terms are treated explicitly \cite{lakkis2013implicit}. The discrete system, with the notation $\Phi^n(x)=\Phi(t_n,x)$, reads
\begin{eqnarray} \label{eq:full_discrete}
    \left<\mathbb{E}(\Phi_a^{n-1})\epsilon(U^n),\epsilon(V) \right>_{Y_h} &=& \left<k_6\mathbb{P}(P_a^{n-1}\hat{\nu}),V \right>_{\Gamma_h},  \nonumber\\
     \Delta t^{-1}\left<\Phi_d^n, \Psi \right>_{Y_h} + \left< D_1 \nabla \Phi_d^n, \nabla \Psi \right>_{Y_h} &+ &\left<C_1\tr(\sigma(U^n))_+\Phi_d^n, \Psi \right>_{Y_h} + \left<n_r\tilde{k}_3\Phi_d^n, \Psi\right>_{\Gamma_h} \nonumber \\
     &=&\Delta t^{-1}\left<\Phi_d^{n-1}, \Psi \right>_{Y_h} + \left<k_1\Phi_a^{n-1},\Psi \right>_{Y_h},  \nonumber\\
    \Delta t^{-1}\left<\Phi_a^n, \Psi \right>_{Y_h} + \left< D_2 \nabla \Phi_a^n, \nabla \Psi \right>_{Y_h} &+& \left<k_1\Phi_a^n,\Psi \right>_{Y_h} \\
    &= & \Delta t^{-1}\left<\Phi_a^{n-1}, \Psi \right>_{Y_h} + \left<C_1\tr(\sigma(U^n))_+\Phi_d^{n}, \Psi \right>_{Y_h} + \left<n_r\tilde{k}_3\Phi_d^{n}, \Psi\right>_{\Gamma_h} , \nonumber\\
   \Delta t^{-1}\left< P_a^n, W \right>_{\Gamma_h} + \left< D_3 \nabla_{\Gamma_h} P_a^n, \nabla_{\Gamma_h} W \right>_{\Gamma_h} &+& \left<\tilde{k}_4(\Phi_a^n)P_a^n,W \right>_{\Gamma_h}\nonumber \\
    &=& \Delta t^{-1}\left< P_a^{n-1}, W \right>_{\Gamma_h} + \left<n_r\tilde{k}_5(\Phi_a^n),W \right>_{\Gamma_h}. \nonumber
\end{eqnarray}
To benchmark the numerical scheme and implementation in FEniCS, we consider $Y$ to be a unit ball and a simplified model  
\begin{equation}\label{eq:benchmark}
\begin{aligned}
    -\nabla \cdot \sigma(u) &= f &&\text{in } \; Y, \\
    \sigma(u)\cdot \nu &= \mathbb{P}(g\rho) &&\text{on }\Gamma,\\
    \partial_t  \phi - \Delta  \phi &= q_1 + \tr(\sigma(u))_{+} \quad &&\text{in } \; Y, \; t>0, \\
    \nabla \phi\cdot \nu &= \rho-\phi &&\text{on }\Gamma, \; t>0, \\
    \partial_t \rho - \Delta_{\Gamma} \rho  &= q_2-\rho+\phi &&\text{on } \Gamma, \; t>0, \\
    u (0, x) = u_{ex}(0,x), &\quad \phi (0, x) = \phi_{ex}(0,x)  \quad &&\mbox{in } Y, \\
    \rho (0, x)& = \rho_{ex}(0,x)  && \mbox{on } \Gamma.
\end{aligned}
\end{equation}
The functions $f$, $g$, $q_1$, and $q_2$ are such that 
$$
\begin{aligned}
    u_{ex}(t,x) &= \left(2x_1^2x_2x_3e^{-4t}, -x_1x_2^2x_3e^{-4t}, -2x_1x_2x_3^2e^{-4t} \right),\\
    \phi_{ex}(t,x) &= \cos(x_1x_2x_3) e^{-4t},\\
    \rho_{ex}(t,x) &= \cos(x_1x_2x_3) e^{-4t} - 3\sin(x_1x_2x_3)e^{-4t}
\end{aligned}
$$
is the exact solution of \eqref{eq:benchmark}. Then for the experimental order of convergence 
$$
{\rm EOC} = \frac{\log(e_n/e_{n-1})}{\log(h_n/h_{n-1})},  
$$
where $h_n $, for $n=1,2,3,4$, are given in Table~\eqref{table:mesh}
and $e_n$ is the error in the $L^2$-norm or the $H^1$-norm, we obtain the second order of convergence in the $L^2$-norm and first order of convergence in $H^1$-norm, see Table~\ref{table:num_error}.  
\begin{table}[!h]
    \centering
    \begin{tabular}{|c|c|c|c|}
        $0.64009064$ & $0.51659533$ & $0.26133991$ & $0.13143219$ 
    \end{tabular}
    \caption{Four mesh sizes used in the calculation of EOC.} \label{table:mesh}
\end{table}
\begin{table}[!h]
    \centering
    \begin{tabular}{|c|c|c|c|}
        EOC for $L^2$ norm for $\phi$ & $1.63892969$ & $1.77954013$ & $1.93991995$  \\
        EOC for $H^1$ norm for $\phi$ & $1.44823967$ & $1.68924575$ & $1.89210920$  \\
        EOC for $L^2$ norm for $\rho$ & $1.58078852$ & $1.78119701$ & $1.96179159$  \\
        EOC for $H^1$ norm for $\rho$ & $1.78659381$ & $1.49338646$ & $1.52956398$  \\
        EOC for $L^2$ norm for $u$ & $0.92059071$ & $1.74758321$ & $1.87989078$  \\
        EOC for $H^1$ norm for $u$ & $1.05358979$ & $1.40005535$ & $1.27787725$ 
    \end{tabular}
     \caption{Experimental order of convergence for numerical scheme~\eqref{eq:full_discrete} considered for model~\eqref{eq:benchmark}.}\label{table:num_error}
\end{table}

\subsection{Conversion from \texorpdfstring{$\mu\mbox{mol}/\mbox{dm}^2$}{micromolar per decimeter squared} to \texorpdfstring{$\# /\mu m^2$}{number of molecules per meter squared}}\label{app:conversion}
Scott et al. \cite{scott_spatial_2021} uses $\mu\mbox{mol}/\mbox{dm}^3$ for concentrations  in the cytoplasm and $\# /\mu m^2$ for concentrations on the plasma membrane, specifically for $\rho_a$. In the models derived and analysed in this work, we use  $\mu\mbox{mol}/\mbox{dm}^3$ for concentrations  in the cytoplasm and $\mu\mbox{mol}/\mbox{dm}^2$ for $\rho_a$. To be able to use the same initial conditions and to compare the results, we need to find a conversion from $\# /\mu m^2$  to $\mu {\rm mol}/dm^2$. We use the fact that the maximum value for $\rho_a$ for large $E$ is $11\cdot 10^{-16} \mu\mbox{mol}/\mu\mbox{m}^2$ as given in \cite[Fig~2B]{scott_spatial_2021}. The maximum value in the numerical results for $\rho_a$ for large $E$ is $593 \# /\mu\mbox{m}^2$ as given in \cite[Fig~3C(ii)]{scott_spatial_2021} and assume these are equivalent. This relation gives the conversion
\begin{equation}\label{conversion}
    10^{-5}\frac{\mu\mbox{mol}}{\mbox{dm}^2} = 10^{-15}\frac{\mu\mbox{mol}}{\mu\mbox{m}^2} = \frac{5930}{11}\# /\mu\mbox{m}^2 = 539.\overline{09} \# /\mu\mbox{m}^2.
\end{equation}

\subsection{Simulations with nucleus}\label{app:sim_nucleus}
To model the inclusion of a nucleus in the cell, we consider the model equations \eqref{eq:coupled_elast} and  \eqref{eq:coupled_reactions} in $Y\setminus \overline Y_{\rm nc}$, where  domain $Y_{\rm nc}$ represents the nucleus. We choose zero flux boundary conditions for $\phi_d$ and $\phi_a$ on $\partial Y_{\rm nc}$. For the mechanics, we choose the interior boundary condition to model the fact that the nucleus is hard to deform
\begin{equation}\label{eq:bc_nucleus}
    \sigma(u) \nu = -\omega u \quad \mbox{ on } \partial Y_{\rm nc},
\end{equation}
where  $\omega$ is a positive constant determining  the rigidity of the nucleus. 

\begin{figure}[!ht]
\begin{minipage}{1\textwidth}
\vspace{-0.04cm}
\begin{minipage}{0.05\textwidth}\centering
\vspace{-0.2cm}
    \textbf{(A)}\\ \vspace{0.3cm}
        $\rho_a$ \\ \vspace{0.7cm}
        $\phi_d$ \\ \vspace{0.7cm}
        $\phi_a$ \\ \vspace{0.7cm}
        $\vert u \vert$ \\ \vspace{0.7cm}
    \textbf{(C)}\\ \vspace{0.3cm}
        $\rho_a$ \\ \vspace{0.7cm}
        $\phi_d$ \\ \vspace{0.7cm}
        $\phi_a$ \\ \vspace{0.7cm}
        $\vert u \vert$ 
\end{minipage}\hfill
\begin{minipage}{0.45\textwidth}
\begin{minipage}{0.33\textwidth} \centering
    {\footnotesize $0.1$kPa } \\ 
        \includegraphics[width=\textwidth]{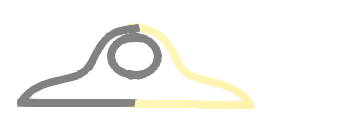}
        \includegraphics[width=\textwidth]{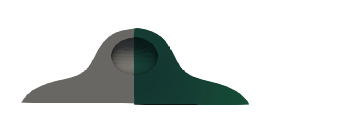}
        \includegraphics[width=\textwidth]{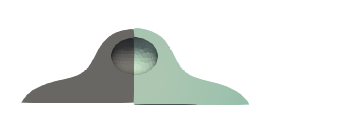}
        \includegraphics[width=\textwidth]{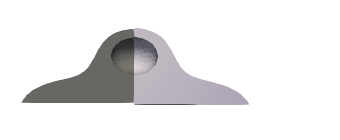}\\
        \vspace{0.5cm}
        \includegraphics[width=\textwidth]{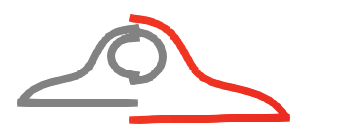}
        \includegraphics[width=\textwidth]{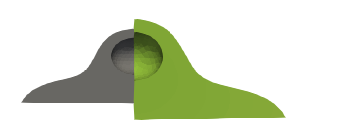}
        \includegraphics[width=\textwidth]{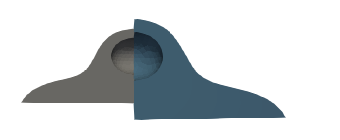}
        \includegraphics[width=\textwidth]{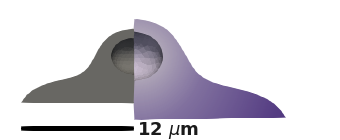}
\end{minipage}\hfill
\begin{minipage}{0.33\textwidth} \centering
    {\footnotesize $5.7$kPa } \\
        \includegraphics[width=\textwidth]{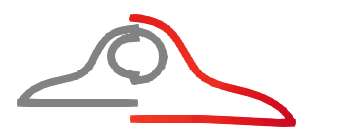}
        \includegraphics[width=\textwidth]{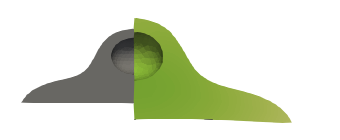}
        \includegraphics[width=\textwidth]{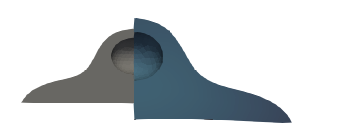}
        \includegraphics[width=\textwidth]{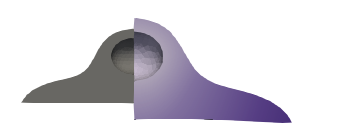}\\
        \vspace{0.5cm}
        \includegraphics[width=\textwidth]{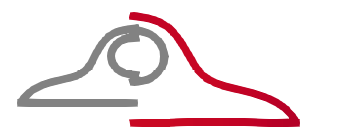}
        \includegraphics[width=\textwidth]{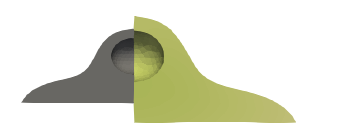}
        \includegraphics[width=\textwidth]{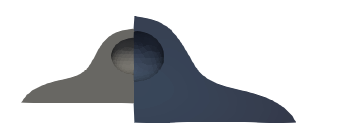}
        \includegraphics[width=\textwidth]{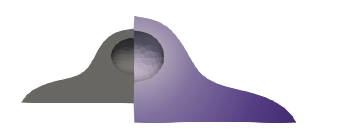}
\end{minipage}\hfill
\begin{minipage}{0.33\textwidth} \centering
    {\footnotesize $7$GPa} \\
        \includegraphics[width=\textwidth]{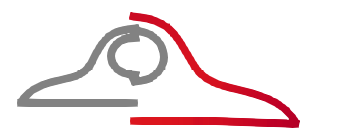}
        \includegraphics[width=\textwidth]{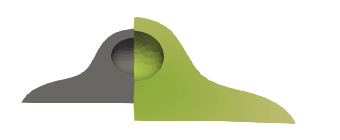}
        \includegraphics[width=\textwidth]{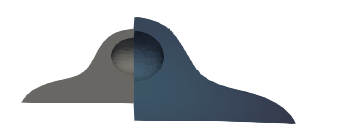}
        \includegraphics[width=\textwidth]{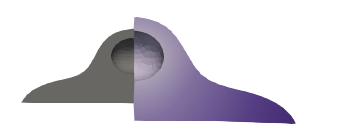}\\
        \vspace{0.5cm}
        \includegraphics[width=\textwidth]{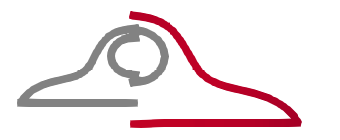}
        \includegraphics[width=\textwidth]{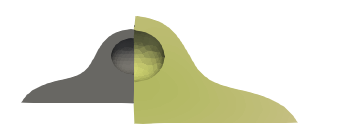}
        \includegraphics[width=\textwidth]{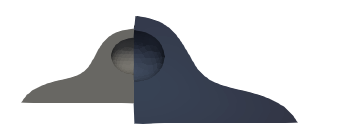}
        \includegraphics[width=\textwidth]{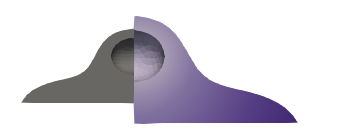}
\end{minipage}\hfill 
\end{minipage}\hfill 
\begin{minipage}{0.05\textwidth}\centering
\vspace{-0.2cm}
    \textbf{(B)}\\ \vspace{4.7cm}
    \textbf{(D)}\\ \vspace{4cm}
\end{minipage}\hfill
\begin{minipage}{0.45\textwidth}
\begin{minipage}{0.33\textwidth} \centering
    {\footnotesize $0.1$kPa } \\ 
        \includegraphics[width=\textwidth]{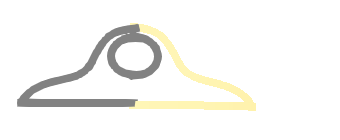}
        \includegraphics[width=\textwidth]{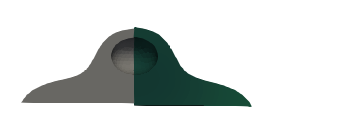}
        \includegraphics[width=\textwidth]{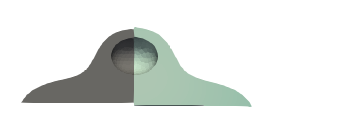}
        \includegraphics[width=\textwidth]{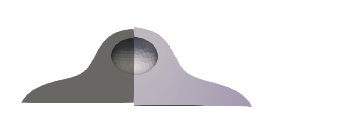}\\
        \vspace{0.5cm}
        \includegraphics[width=\textwidth]{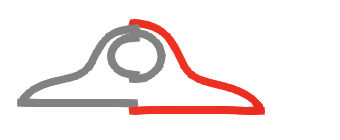}
        \includegraphics[width=\textwidth]{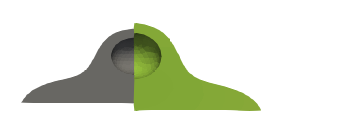}
        \includegraphics[width=\textwidth]{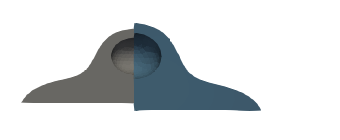}
        \includegraphics[width=\textwidth]{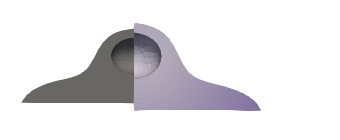}
\end{minipage}\hfill
\begin{minipage}{0.33\textwidth} \centering
    {\footnotesize $5.7$kPa } \\
        \includegraphics[width=\textwidth]{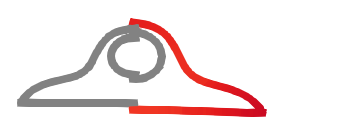}
        \includegraphics[width=\textwidth]{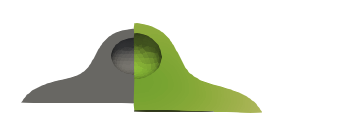}
        \includegraphics[width=\textwidth]{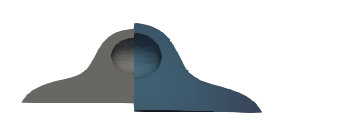}
        \includegraphics[width=\textwidth]{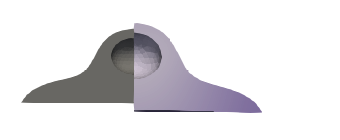}\\
        \vspace{0.5cm}
        \includegraphics[width=\textwidth]{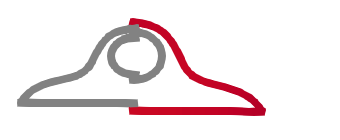}
        \includegraphics[width=\textwidth]{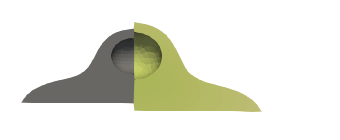}
        \includegraphics[width=\textwidth]{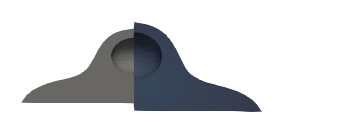}
        \includegraphics[width=\textwidth]{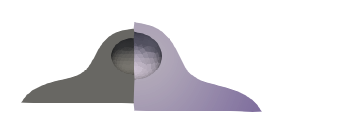}
\end{minipage}\hfill
\begin{minipage}{0.33\textwidth} \centering
    {\footnotesize $7$GPa} \\
        \includegraphics[width=\textwidth]{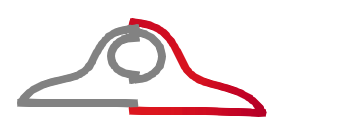}
        \includegraphics[width=\textwidth]{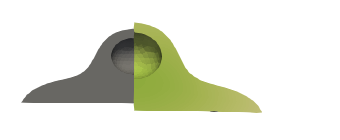}
        \includegraphics[width=\textwidth]{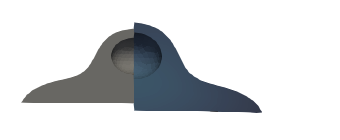}
        \includegraphics[width=\textwidth]{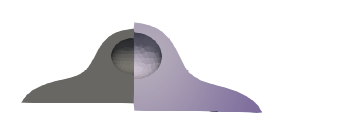}\\
        \vspace{0.5cm}
        \includegraphics[width=\textwidth]{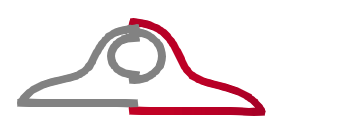}
        \includegraphics[width=\textwidth]{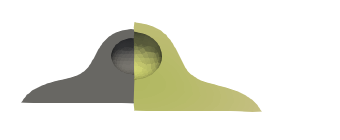}
        \includegraphics[width=\textwidth]{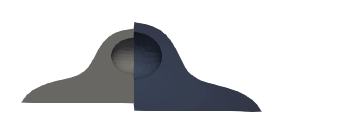}
        \includegraphics[width=\textwidth]{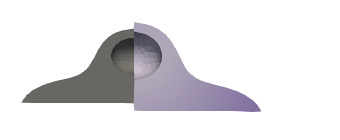}
\end{minipage}\hfill
\end{minipage}
\end{minipage}\hfill
\begin{minipage}{1\textwidth}\raggedright
    \hspace{1cm}
    \includegraphics[width=0.22\textwidth]{figures/coupled_p_bar_grouped.png}
    \includegraphics[width=0.22\textwidth]{figures/coupled_cd_bar_grouped.png}
    \includegraphics[width=0.22\textwidth]{figures/coupled_ca_bar_grouped.png}
    \includegraphics[width=0.22\textwidth]{figures/coupled_u_bar_grouped.png}
\end{minipage}
  \caption{{\bf Numerical simulation results showing $\rho_a$, $\phi_d$, $\phi_a$ and $\vert u\vert$  for model~\eqref{eq:coupled_elast}, \eqref{eq:coupled_force_boundary}, \eqref{eq:coupled_reactions}, and \eqref{eq:bc_nucleus}  for the axisymmetric shape with a nucleus  and in the case of the $3$D stimulus at a steady state at $T=100$~s.}
    Four different scenarios are considered: \textbf{(A)} $C_1=0$~$(\mbox{kPa s})^{-1}$ $(\sigma\not\rightarrow \phi_a)$ and $E_c=0.6$~kPa $(\phi_a\not\rightarrow E_c)$; \textbf{(B)} $C_1=0$~$(\mbox{kPa s})^{-1}$ $(\sigma\not\rightarrow \phi_a)$ and $E_c=f(\phi_a)$ $(\phi_a\rightarrow E_c)$; \textbf{(C)} $C_1=0.1$~$(\mbox{kPa s})^{-1}$ $(\sigma\rightarrow \phi_a)$ and $E_c=0.6$~kPa $(\phi_a\not\rightarrow E_c)$; \textbf{(D)} $C_1=0.1$~$(\mbox{kPa s})^{-1}$ $(\sigma\rightarrow \phi_a)$ and $E_c=f(\phi_a)$ $(\phi_a\rightarrow E_c)$. Within each subfigure, the rows represent $\rho_a$, $\phi_d$, $\phi_a$ and $\vert u\vert$ on a cross-section of the plane $x_1=0$ of the axisymmetric cell, and the columns represent $E=0.1, 5.7, 7\cdot 10^6$~kPa. Parameter values as in Table~\ref{tab:parameters}, and  $\omega=1$. The corresponding results without a nucleus can be found in Fig~\ref{fig:sim_subs}.}
  \label{fig:sim_nucl}
\end{figure}

\newpage
\subsection{Parameter sensitivity analysis} \label{app:param_sens}
As we compare and extend upon the model of \cite{scott_spatial_2021}, we have used the same parameter values where possible. For all new values introduced in this work, see Table~\ref{tab:parameters}, we provide here a simple parameter sensitivity analysis to explore how different choices of these parameters would affect the results. We choose to focus on the case of the 3D stimulus for the model in Eqs~\eqref{eq:coupled_elast}-\eqref{eq:coupled_force_boundary}, \eqref{eq:coupled_reactions} with force boundary conditions on the whole cell membrane for the axisymmetric cell, where $E_c=f(\phi_a)$ and $C_1=0.1$~$(\mbox{kPa s})^{-1}$ when not otherwise specified, and for $E=0.1$~kPa, $5.7$~kPa, $7$~GPa. Parameter values are chosen as in Table~\ref{tab:parameters} and we show summary statistics for a $10\%$ and $20\%$ change in values for each of the following parameters separately: $C_1$, $k_6$, $k_7$, $k_8$, $p$ and $\nu_c$. 

\begin{figure}[!h]
\begin{minipage}{\textwidth} 
    \hspace{4cm} $E=0.1$kPa \hspace{2.5cm} $E=5.7$kPa \hspace{2.5cm} $E=7$GPa \\ 
    \vspace{-0.2cm}
\end{minipage}\hfill
\begin{minipage}{\textwidth} \centering
    \includegraphics[width=0.85\textwidth]{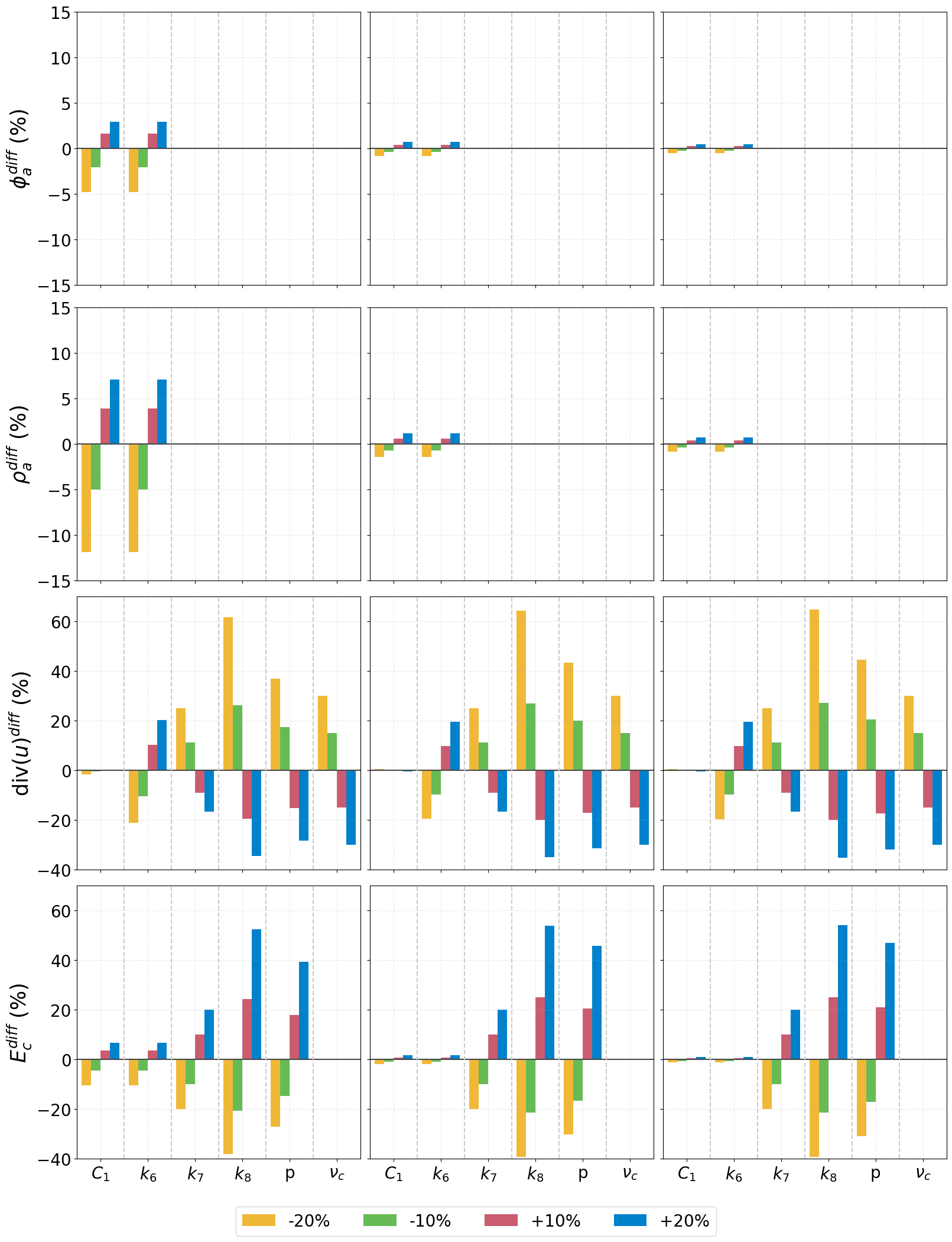}
\end{minipage}
  \caption{{\bf Parameter sensitivity analysis results showing the percentage change in the means of $f(\phi_a)$, ${\rm div}(u)$, $\phi_a$ and $\rho_a$ (\%) as response to a percentage change in the parameters $C_1$, $k_6$, $k_7$, $k_8$, $p$ and $\nu_c$.} Here, the percentage change in the mean is defined as $v^{diff}=\frac{\int_\Omega v \dd{x}-\int_\Omega v^* \dd{x}}{\int_\Omega v^* \dd{x}}\cdot 100$, where $v$ the new result and $v^*$ the original result for parameters in Table~\ref{tab:parameters}. These results are for the model~\eqref{eq:coupled_elast}, \eqref{eq:coupled_force_boundary}, and \eqref{eq:coupled_reactions}, for the axisymmetric shape and in the case of the $3$D stimulus at $T=100$~s by which time the results are at a steady state. The columns represent different values for the substrate stiffness $E$ and the rows represent the percentage change in the mean of the variables $\phi_a$, $\rho_a$, div$(u)$ and $E_c=f(\phi_a)$. If not otherwise specified, parameter values are as in  Table~\ref{tab:parameters}. The corresponding simulation results can be found in Fig~\ref{fig:sim_subs}.}
  \label{fig:param_anal}
\end{figure}

In Fig~\ref{fig:param_anal}, we see all parameters only have a small effect on the concentrations of $\phi_a$ and $\rho_a$ with $C_1$ and $k_6$ being the most influential when $E=0.1$~kPa. The parameters $k_7$, $k_8$, $p$ and $\nu_c$ are part of the equations of linear elasticity and therefore their effect would be through changes in the stress in the reaction term of $\phi_a$, which is also determined by $C_1$. The change in  the stress due to the change in these parameters  barely changes the results of the signalling molecules. Next to that, we see that the effect is larger for smaller substrate stiffness $E$, which could be due to the fact that $\phi_a$ and $\rho_a$ are already close to their maximum values for larger substrate stiffnesses. Fig~\ref{fig:param_anal} shows that the effect of $C_1$ on the change in $\rho_a$ is approximately  tripled compared to the change in $\phi_a$. This could be because the reaction term including $\phi_a$ in $\rho_a$ is to the power of $5$. 

Fig~\ref{fig:param_anal} also shows that $k_7$, $k_8$, and $p$ have an opposite effect on the cell stiffness $E_c$ and the change in volume div$(u)$, where an increase in these parameters increases the stiffness and decreases div$(u)$. The parameter $k_8$ has the largest effect on these variables, where a $-20\%$ change has a tripling effect of a $+60\%$ change in div$(u)$. Lastly, we see that $\nu_c$ only significantly changes  div$(u)$. The effect that $\nu_c$, the Poisson ratio, determining the ratio between the transverse and the axial strain, has on the volume change does not seem to induce changes in other variables.

\clearpage
\subsection{Simulations with linear viscoelasticity} \label{app:viscoelas}
\begin{figure}[!ht]
\begin{minipage}{1\textwidth}
\begin{minipage}{0.05\textwidth}\centering
\vspace{-0.2cm}
    \textbf{(A)}\\ \vspace{0.3cm}
        $\rho_a$ \\ \vspace{0.7cm}
        $\phi_d$ \\ \vspace{0.7cm}
        $\phi_a$ \\ \vspace{0.7cm}
        $\vert u \vert$ \\ \vspace{0.7cm}
    \textbf{(C)}\\ \vspace{0.3cm}
        $\rho_a$ \\ \vspace{0.7cm}
        $\phi_d$ \\ \vspace{0.7cm}
        $\phi_a$ \\ \vspace{0.7cm}
        $\vert u \vert$ 
\end{minipage}\hfill
\begin{minipage}{0.45\textwidth}
\begin{minipage}{0.33\textwidth} \centering
    {\footnotesize $0.1$kPa } \\ 
        \includegraphics[width=\textwidth]{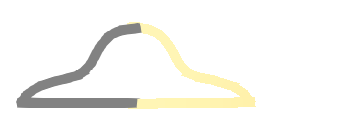}
        \includegraphics[width=\textwidth]{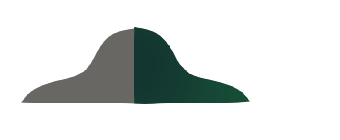}
        \includegraphics[width=\textwidth]{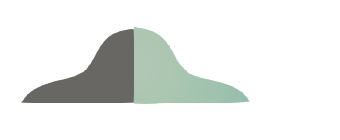}
        \includegraphics[width=\textwidth]{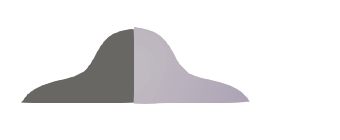}\\
        \vspace{0.5cm}
        \includegraphics[width=\textwidth]{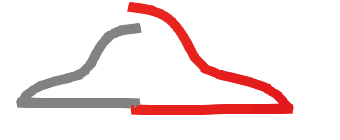}
        \includegraphics[width=\textwidth]{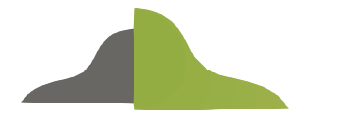}
        \includegraphics[width=\textwidth]{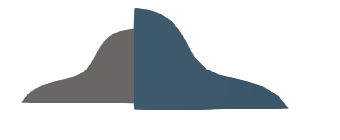}
        \includegraphics[width=\textwidth]{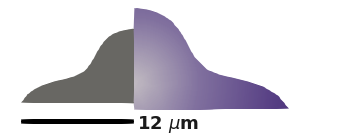}
\end{minipage}\hfill
\begin{minipage}{0.33\textwidth} \centering
    {\footnotesize $5.7$kPa } \\
        \includegraphics[width=\textwidth]{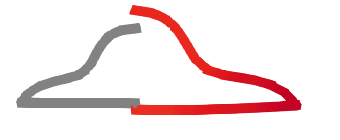}
        \includegraphics[width=\textwidth]{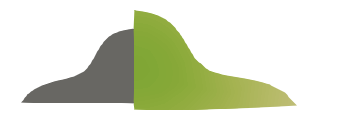}
        \includegraphics[width=\textwidth]{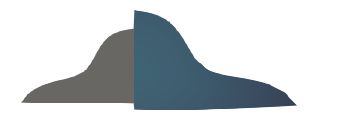}
        \includegraphics[width=\textwidth]{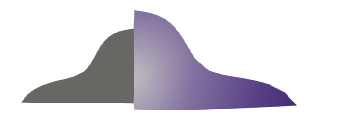}\\
        \vspace{0.5cm}
        \includegraphics[width=\textwidth]{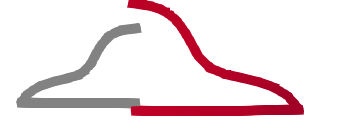}
        \includegraphics[width=\textwidth]{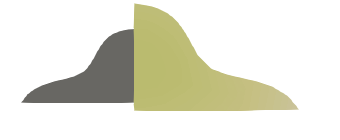}
        \includegraphics[width=\textwidth]{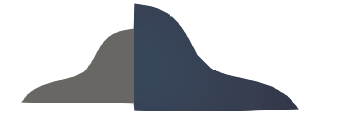}
        \includegraphics[width=\textwidth]{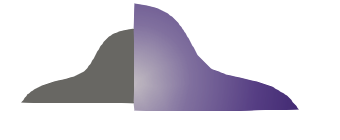}
\end{minipage}\hfill
\begin{minipage}{0.33\textwidth} \centering
    {\footnotesize $7$GPa} \\
        \includegraphics[width=\textwidth]{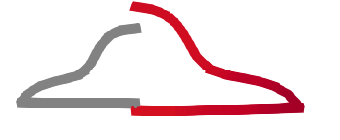}
        \includegraphics[width=\textwidth]{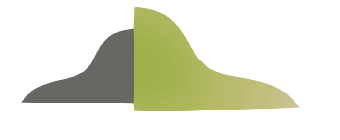}
        \includegraphics[width=\textwidth]{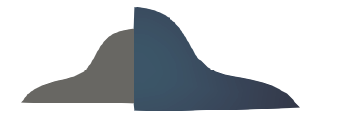}
        \includegraphics[width=\textwidth]{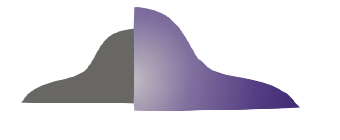}\\
        \vspace{0.5cm}
        \includegraphics[width=\textwidth]{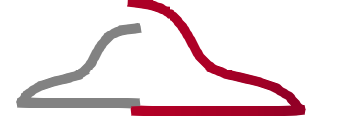}
        \includegraphics[width=\textwidth]{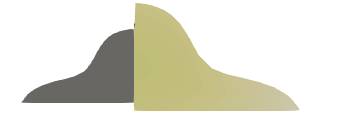}
        \includegraphics[width=\textwidth]{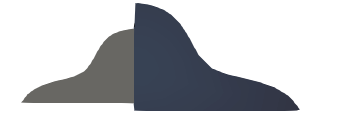}
        \includegraphics[width=\textwidth]{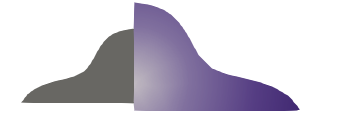}
\end{minipage}\hfill 
\end{minipage}\hfill 
\begin{minipage}{0.05\textwidth}\centering
\vspace{-0.2cm}
    \textbf{(B)}\\ \vspace{4.7cm}
    \textbf{(D)}\\ \vspace{4cm}
\end{minipage}\hfill
\begin{minipage}{0.45\textwidth}
\begin{minipage}{0.33\textwidth} \centering
    {\footnotesize $0.1$kPa } \\ 
        \includegraphics[width=\textwidth]{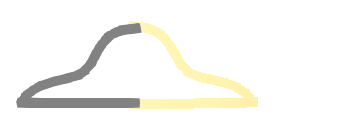}
        \includegraphics[width=\textwidth]{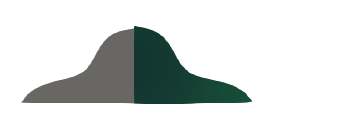}
        \includegraphics[width=\textwidth]{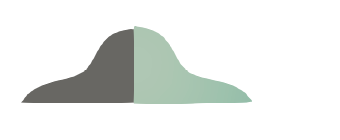}
        \includegraphics[width=\textwidth]{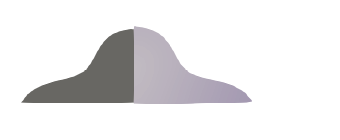}\\
        \vspace{0.5cm}
        \includegraphics[width=\textwidth]{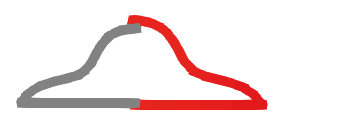}
        \includegraphics[width=\textwidth]{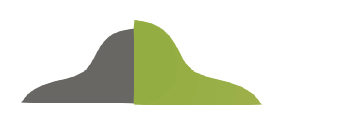}
        \includegraphics[width=\textwidth]{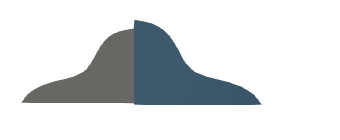}
        \includegraphics[width=\textwidth]{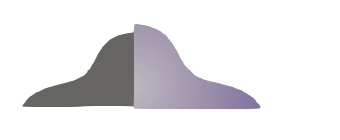}
\end{minipage}\hfill
\begin{minipage}{0.33\textwidth} \centering
    {\footnotesize $5.7$kPa } \\
        \includegraphics[width=\textwidth]{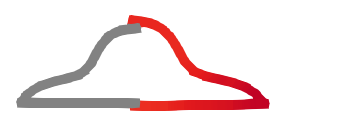}
        \includegraphics[width=\textwidth]{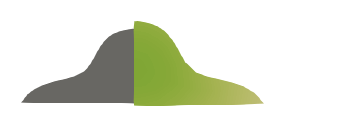}
        \includegraphics[width=\textwidth]{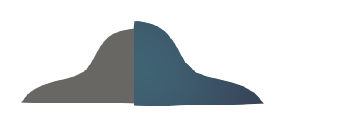}
        \includegraphics[width=\textwidth]{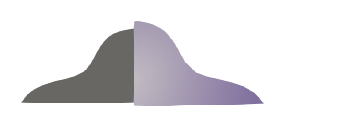}\\
        \vspace{0.5cm}
        \includegraphics[width=\textwidth]{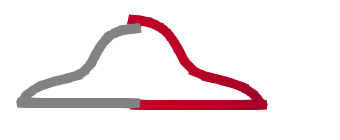}
        \includegraphics[width=\textwidth]{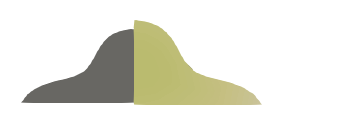}
        \includegraphics[width=\textwidth]{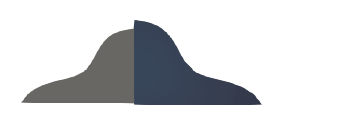}
        \includegraphics[width=\textwidth]{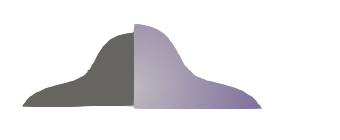}
\end{minipage}\hfill
\begin{minipage}{0.33\textwidth} \centering
    {\footnotesize $7$GPa} \\
        \includegraphics[width=\textwidth]{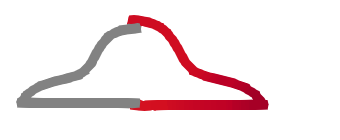}
        \includegraphics[width=\textwidth]{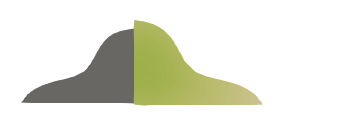}
        \includegraphics[width=\textwidth]{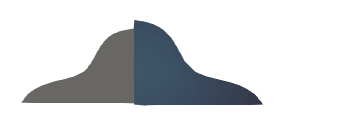}
        \includegraphics[width=\textwidth]{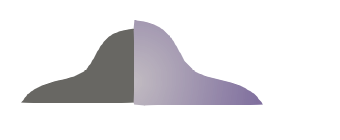}\\
        \vspace{0.5cm}
        \includegraphics[width=\textwidth]{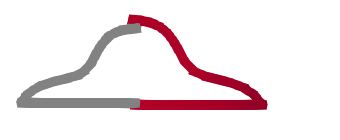}
        \includegraphics[width=\textwidth]{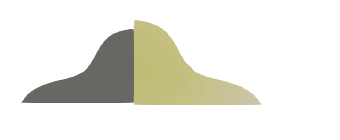}
        \includegraphics[width=\textwidth]{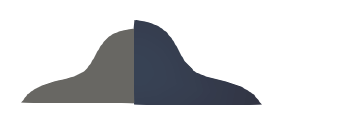}
        \includegraphics[width=\textwidth]{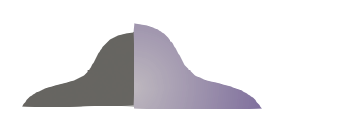}
\end{minipage}\hfill
\end{minipage}
\end{minipage}\hfill
\begin{minipage}{1\textwidth}\raggedright
    \hspace{1cm}
    \includegraphics[width=0.22\textwidth]{figures/coupled_p_bar_grouped.png}
    \includegraphics[width=0.22\textwidth]{figures/coupled_cd_bar_grouped.png}
    \includegraphics[width=0.22\textwidth]{figures/coupled_ca_bar_grouped.png}
    \includegraphics[width=0.22\textwidth]{figures/coupled_u_bar_grouped.png}
\end{minipage}
  \caption{{\bf Numerical simulation results showing $\rho_a$, $\phi_d$, $\phi_a$ and $\vert u\vert$  for model~\eqref{eq:coupled_elast}, \eqref{eq:coupled_force_boundary}, \eqref{eq:coupled_reactions}, and \eqref{eq:sigma_visco}  for the axisymmetric shape and in the case of the $3$D stimulus at a steady state at $T=100$~s.}
    Four different scenarios are considered: \textbf{(A)} $C_1=0$~$(\mbox{kPa s})^{-1}$ $(\sigma\not\rightarrow \phi_a)$ and $E_c=0.6$~kPa $(\phi_a\not\rightarrow E_c)$; \textbf{(B)} $C_1=0$~$(\mbox{kPa s})^{-1}$ $(\sigma\not\rightarrow \phi_a)$ and $E_c=f(\phi_a)$ $(\phi_a\rightarrow E_c)$; \textbf{(C)} $C_1=0.1$~$(\mbox{kPa s})^{-1}$ $(\sigma\rightarrow \phi_a)$ and $E_c=0.6$~kPa $(\phi_a\not\rightarrow E_c)$; \textbf{(D)} $C_1=0.1$~$(\mbox{kPa s})^{-1}$ $(\sigma\rightarrow \phi_a)$ and $E_c=f(\phi_a)$ $(\phi_a\rightarrow E_c)$. Within each subfigure, the rows represent $\rho_a$, $\phi_d$, $\phi_a$ and $\vert u\vert$ on a cross-section of the plane $x_1=0$ of the axisymmetric cell, and the columns represent $E=0.1, 5.7, 7\cdot 10^6$~kPa. Parameter values as in Table~\ref{tab:parameters}, and  $\theta_\lambda=\theta_\mu=1$. The corresponding results with only elastic stress can be found in Fig~\ref{fig:sim_subs}.}
  \label{fig:sim_visc}
\end{figure}

To model a linear viscoelastic material, we consider the model equations \eqref{eq:coupled_elast} and  \eqref{eq:coupled_reactions} with the stress defined as
\begin{equation}\label{eq:sigma_visco} 
    \sigma(u) = \lambda(\phi_a)(\nabla \cdot u + \theta_\lambda \nabla \cdot \partial_t u)I +  2\mu(\phi_a)\left( \big( \nabla u + (\nabla u)^T\big) + \theta_\mu \big( \nabla \partial_t u + (\nabla \partial_t u)^T\big)\right),
\end{equation}
with the initial condition $u(0,x)=0$ for $x \in Y$, where  $\theta_\lambda$ and $\theta_\mu$ are the characteristic retardation times. We use a backward Euler discretization such that
\begin{align*}
    \left<\lambda(\Phi_a^{n-1})\right.&(1+\Delta t^{-1}\theta_\lambda)(\nabla \cdot U^n)I +\left.  2\mu(\Phi_a^{n-1})(1+\Delta t^{-1}\theta_\mu) \epsilon(U^n),\epsilon(V) \right>_{Y_h}\\
    &= \Delta t^{-1} \left<\lambda(\Phi_a^{n-1})\theta_\lambda (\nabla \cdot U^{n-1})I +  2\mu(\Phi_a^{n-1})\theta_\mu \epsilon(U^{n-1}),\epsilon(V) \right>_{Y_h} + \left<k_6\mathbb{P}(P_a^{n-1}\hat{\nu}),V \right>_{\Gamma_h}.
\end{align*}
The numerical simulation results are presented in Fig~\ref{fig:sim_visc}.

\begin{thebibliography}{10}

  \bibitem{Tomar_2009}
  Tomar A, Lim ST, Lim Y, Schlaepfer DD.
  \newblock A {FAK-p120RasGAP-p190RhoGAP} complex regulates polarity in migrating cells.
  \newblock Journal of Cell Science. 2009;122:1852--1862.
  \newblock doi:{10.1242/jcs.046870}.
  
  \bibitem{Valls_2022}
  Valls PO, Esposito A.
  \newblock Signalling dynamics, cell decisions, and homeostatic control in health and disease.
  \newblock J Cell Sci Curr Opin Cell Biol. 2022;75:102066.
  \newblock doi:{10.1016/j.ceb.2022.01.011}.
  
  \bibitem{romani_crosstalk_2021}
  Romani P, Valcarcel-Jimenez L, Frezza C, Dupont S.
  \newblock Crosstalk between mechanotransduction and metabolism.
  \newblock Nature Reviews Molecular Cell Biology. 2021;22(1):22--38.
  \newblock doi:{10.1038/s41580-020-00306-w}.
  
  \bibitem{Cai_2021}
  Cai X, Wang KC, Meng Z.
  \newblock Mechanoregulation of YAP and TAZ in Cellular Homeostasis and Disease Progression.
  \newblock Front Cell Dev Biol. 2021;9:673599.
  \newblock doi:{10.3389/fcell.2021.673599}.
  
  \bibitem{humphrey_mechanotransduction_2014}
  Humphrey JD, Dufresne ER, Schwartz MA.
  \newblock Mechanotransduction and extracellular matrix homeostasis.
  \newblock Nature reviews Molecular cell biology. 2014;15(12):802--812.
  \newblock doi:{10.1038/nrm3896}.
  
  \bibitem{saraswathibhatla_cellextracellular_2023}
  Saraswathibhatla A, Indana D, Chaudhuri O.
  \newblock Cell–extracellular matrix mechanotransduction in {3D}.
  \newblock Nature Reviews Molecular Cell Biology. 2023;24(7):495--516.
  \newblock doi:{10.1038/s41580-023-00583-1}.
  
  \bibitem{burridge_mechanotransduction_2019}
  Burridge K, Monaghan-Benson E, Graham DM.
  \newblock Mechanotransduction: from the cell surface to the nucleus via {RhoA}.
  \newblock Philosophical Transactions of the Royal Society B: Biological Sciences. 2019;374(1779):20180229.
  \newblock doi:{10.1098/rstb.2018.0229}.
  
  \bibitem{xie_cell_2023}
  Xie N, Xiao C, Shu Q, Cheng B, Wang Z, Xue R, et~al.
  \newblock Cell response to mechanical microenvironment cues via {Rho} signaling: {From} mechanobiology to mechanomedicine.
  \newblock Acta Biomaterialia. 2023;159:1--20.
  \newblock doi:{10.1016/j.actbio.2023.01.039}.
  
  \bibitem{martino_cellular_2018}
  Martino F, Perestrelo AR, Vinarský V, Pagliari S, Forte G.
  \newblock Cellular {Mechanotransduction}: {From} {Tension} to {Function}.
  \newblock Frontiers in Physiology. 2018;9.
  
  \bibitem{Sun_2016}
  Sun Z, Guo SS, F\"assler R.
  \newblock Integrin-mediated mechanotransduction.
  \newblock J Cell Biol. 2016;215(4):445--456.
  \newblock doi:{10.1083/jcb.201609037}.
  
  \bibitem{Young_2023}
  Young KM, Reinhart-King CA.
  \newblock Environmental stiffness restores mechanical homeostasis in vimentin‐depleted cells.
  \newblock Current Opinion in Cell Biology. 2023;83:102208.
  \newblock doi:{10.1016/j.ceb.2023.102208}.
  
  \bibitem{gilbert2006computational}
  Gilbert D, Fuss H, Gu X, Orton R, Robinson S, Vyshemirsky V, et~al.
  \newblock Computational methodologies for modelling, analysis and simulation of signalling networks.
  \newblock Briefings in Bioinformatics. 2006;7(4):339--353.
  \newblock doi:{10.1093/bib/bbl043}.
  
  \bibitem{garcia_2014}
  Garc\'ia-Penarrubia P, J~J~G\'alvez JJ, G\'alvez J.
  \newblock Mathematical modelling and computational study of two-dimensional and three-dimensional dynamics of receptor–ligand interactions in signalling response mechanisms.
  \newblock J Math Biol. 2014;69:553--582.
  
  \bibitem{ptashnyk_multiscale_2020}
  Ptashnyk M, Venkataraman C.
  \newblock Multiscale Analysis and Simulation of a Signaling Process With Surface Diffusion.
  \newblock Multiscale Modeling \& Simulation. 2020;18(2):851--886.
  \newblock doi:{10.1137/18M1185661}.
  
  \bibitem{cheng_cellular_2017}
  Cheng B, Lin M, Huang G, Li Y, Ji B, Genin GM, et~al.
  \newblock Cellular mechanosensing of the biophysical microenvironment: {A} review of mathematical models of biophysical regulation of cell responses.
  \newblock Physics of Life Reviews. 2017;22-23:88--119.
  \newblock doi:{10.1016/j.plrev.2017.06.016}.
  
  \bibitem{besser_coupling_2007}
  Besser A, Schwarz US.
  \newblock Coupling biochemistry and mechanics in cell adhesion: a model for inhomogeneous stress fiber contraction.
  \newblock New Journal of Physics. 2007;9(11):425--425.
  \newblock doi:{10.1088/1367-2630/9/11/425}.
  
  \bibitem{novev_spatiotemporal_2021}
  Novev JK, Heltberg ML, Jensen MH, Doostmohammadi A.
  \newblock Spatiotemporal model of cellular mechanotransduction via {Rho} and {YAP}.
  \newblock Integrative Biology. 2021;13(8):197--209.
  \newblock doi:{10.1093/intbio/zyab012}.
  
  \bibitem{scott_spatial_2021}
  Scott KE, Fraley SI, Rangamani P.
  \newblock A spatial model of {YAP}/{TAZ} signaling reveals how stiffness, dimensionality, and shape contribute to emergent outcomes.
  \newblock Proceedings of the National Academy of Sciences. 2021;118(20):e2021571118.
  \newblock doi:{10.1073/pnas.2021571118}.
  
  \bibitem{kang_structurally_2015}
  Kang J, Puskar KM, Ehrlicher AJ, LeDuc PR, Schwartz RS.
  \newblock Structurally {Governed} {Cell} {Mechanotransduction} through {Multiscale} {Modeling}.
  \newblock Scientific Reports. 2015;5(1):8622.
  \newblock doi:{10.1038/srep08622}.
  
  \bibitem{bar-ziv_pearling_1999}
  Bar-Ziv R, Tlusty T, Moses E, Safran SA, Bershadsky A.
  \newblock Pearling in cells: {A} clue to understanding cell shape.
  \newblock Proceedings of the National Academy of Sciences. 1999;96(18):10140--10145.
  \newblock doi:{10.1073/pnas.96.18.10140}.
  
  \bibitem{vianay_single_2010}
  Vianay B, Käfer J, Planus E, Block M, Graner F, Guillou H.
  \newblock Single {Cells} {Spreading} on a {Protein} {Lattice} {Adopt} an {Energy} {Minimizing} {Shape}.
  \newblock Physical Review Letters. 2010;105(12):128101.
  \newblock doi:{10.1103/PhysRevLett.105.128101}.
  
  \bibitem{albert_dynamics_2014}
  Albert PJ, Schwarz US.
  \newblock Dynamics of cell shape and forces on micropatterned substrates predicted by a cellular {Potts} model.
  \newblock Biophysical Journal. 2014;106(11):2340--2352.
  \newblock doi:{10.1016/j.bpj.2014.04.036}.
  
  \bibitem{sun_computational_2016}
  Sun M, Spill F, Zaman MH.
  \newblock A Computational Model of {YAP}/{TAZ} Mechanosensing.
  \newblock Biophysical Journal. 2016;110(11):2540.
  \newblock doi:{10.1016/j.bpj.2016.04.040}.
  
  \bibitem{eroume_exploring_2021}
  Eroumé KS, Cavill R, Staňková K, de~Boer J, Carlier A.
  \newblock Exploring the influence of cytosolic and membrane {FAK} activation on {YAP}/{TAZ} nuclear translocation.
  \newblock Biophysical Journal. 2021;120(20):4360--4377.
  \newblock doi:{10.1016/j.bpj.2021.09.009}.
  
  \bibitem{gould_introduction_2013}
  Gould PL.
  \newblock Introduction to Linear Elasticity.
  \newblock Springer New York; 2013.
  \newblock Available from: \url{http://link.springer.com/10.1007/978-1-4614-4833-4}.
  
  \bibitem{dziuk2013finite}
  Dziuk G, Elliott CM.
  \newblock Finite element methods for surface {PDEs}.
  \newblock Acta Numerica. 2013;22:289--396.
  \newblock doi:{10.1017/S0962492913000056}.
  
  \bibitem{Gilbert_2017}
  Gilbert PM, Weaver VM.
  \newblock Cellular adaptation to biomechanical stress across length scales in tissue homeostasis and disease.
  \newblock Seminars Cell Develop Biol. 2017;67:141--152.
  \newblock doi:{10.1016/j.semcdb.2016.09.004}.
  
  \bibitem{Grolleman_2023}
  Grolleman J, van Engeland NCA, Raza M, Azimi S, Conte V, Sahlgren CM, et~al.
  \newblock Environmental stiffness restores mechanical homeostasis in vimentin‐depleted cells.
  \newblock Nature, Scientific Reports. 2023;13:18374.
  \newblock doi:{10.1038/s41598-023-44835-8}.
  
  \bibitem{mcnicol_theoretical_2025}
  McNicol GR, Dalby MJ, Stewart PS.
  \newblock A theoretical model for focal adhesion and cytoskeleton formation in non-motile cells.
  \newblock Journal of Theoretical Biology. 2025;596:111965.
  \newblock doi:{10.1016/j.jtbi.2024.111965}.
  
  \bibitem{Di_2023}
  Di X, Gao X, Peng L, Ai J, Jin X, Qi S, et~al.
  \newblock Cellular mechanotransduction in health and diseases: from molecular mechanism to therapeutic targets.
  \newblock Sig Transduct Target Ther. 2023;8:282.
  \newblock doi:{10.1038/s41392-023-01501-9}.
  
  \bibitem{kasza_cell_2007}
  Kasza KE, Rowat AC, Liu J, Angelini TE, Brangwynne CP, Koenderink GH, et~al.
  \newblock The cell as a material.
  \newblock Current Opinion in Cell Biology. 2007;19(1):101--107.
  \newblock doi:{10.1016/j.ceb.2006.12.002}.
  
  \bibitem{moeendarbary_cytoplasm_2013}
  Moeendarbary E, Valon L, Fritzsche M, Harris AR, Moulding DA, Thrasher AJ, et~al.
  \newblock The cytoplasm of living cells behaves as a poroelastic material.
  \newblock Nature Materials. 2013;12(3):253--261.
  \newblock doi:{10.1038/nmat3517}.
  
  \bibitem{banerjee_controlling_2013}
  Banerjee S, Marchetti MC.
  \newblock Controlling cell–matrix traction forces by extracellular geometry.
  \newblock New Journal of Physics. 2013;15(3):035015.
  \newblock doi:{10.1088/1367-2630/15/3/035015}.
  
  \bibitem{oakes_geometry_2014}
  Oakes PW, Banerjee S, Marchetti MC, Gardel ML.
  \newblock Geometry {Regulates} {Traction} {Stresses} in {Adherent} {Cells}.
  \newblock Biophysical Journal. 2014;107(4):825--833.
  \newblock doi:{10.1016/j.bpj.2014.06.045}.
  
  \bibitem{chojowski_reversible_2020}
  Chojowski R, Schwarz US, Ziebert F.
  \newblock Reversible elastic phase field approach and application to cell monolayers.
  \newblock The European Physical Journal E. 2020;43(10):63.
  \newblock doi:{10.1140/epje/i2020-11988-1}.
  
  \bibitem{chengPredicting2023}
  Cheng B, Li M, Wan W, Guo H, Genin GM, Lin M, et~al.
  \newblock Predicting {{YAP}}/{{TAZ}} Nuclear Translocation in Response to {{ECM}} Mechanosensing.
  \newblock Biophysical Journal. 2023;122(1):43--53.
  \newblock doi:{10.1016/j.bpj.2022.11.2943}.
  
  \bibitem{Graham_2016}
  Graham DM, Burridge K.
  \newblock Mechanotransduction and nuclear function.
  \newblock Current Opinion in Cell Biology. 2016;40:98--105.
  \newblock doi:{10.1016/j.ceb.2016.03.006}.

  \bibitem{gardel_elastic_2004}
  Gardel ML, Shin JH, {MacKintosh} FC, Mahadevan L, Matsudaira P, Weitz DA.
  \newblock Elastic Behavior of Cross-Linked and Bundled Actin Networks.
  \newblock Science. 2004;304(5675):1301--1305.
  \newblock doi:{10.1126/science.1095087}.
  
  \bibitem{Burridge_2016}
  Burridge K, Guilluy C.
  \newblock Focal adhesions, stress fibers and mechanical tension.
  \newblock Experimental Cell Research. 2016;343:14--20.
  \newblock doi:{10.1016/j.yexcr.2015.10.029}.
  
  \bibitem{Doyle_2015}
  Doyle AD, Carvajal N, Jin A, Matsumoto K, Yamada KM.
  \newblock Local {3D} matrix microenvironment regulates cell migration through spatiotemporal dynamics of contractility-dependent adhesions.
  \newblock J Cell Sci. 2015;6:8720.
  \newblock doi:{10.1038/ncomms9720}.
  
  \bibitem{Zhao_2007}
  Zhao XH, Laschinger C, Arora P, Sz{\'a}szi K, Kapus A, McCulloch CA.
  \newblock Force activates smooth muscle $\alpha$-actin promoter activity through the Rho signaling pathway.
  \newblock J Cell Sci. 2007;120:1801--1809.
  \newblock doi:{10.1242/jcs.001586}.
  
  \bibitem{mokbelPoisson2020}
  Mokbel M, Hosseini K, Aland S, Fischer-Friedrich E.
  \newblock The Poisson Ratio of the Cellular Actin Cortex Is Frequency Dependent.
  \newblock Biophysical Journal. 2020;118(8):1968--1976.
  \newblock doi:{10.1016/j.bpj.2020.03.002}.

  \bibitem{chrzanowska-wodnicka_rho-stimulated_1996}
  Chrzanowska-Wodnicka M, Burridge K.
  \newblock Rho-stimulated contractility drives the formation of stress fibers and focal adhesions.
  \newblock Journal of Cell Biology. 1996;133(6):1403--1415.
  \newblock doi:{10.1083/jcb.133.6.1403}.
  
  \bibitem{beamish_engineered_2017}
  Beamish JA, Chen E, Putnam AJ.
  \newblock Engineered extracellular matrices with controlled mechanics modulate renal proximal tubular cell epithelialization.
  \newblock PLOS ONE. 2017;12(7):e0181085.
  \newblock doi:{10.1371/journal.pone.0181085}.
  
  \bibitem{chen_cell_2003}
  Chen CS, Alonso JL, Ostuni E, Whitesides GM, Ingber DE.
  \newblock Cell shape provides global control of focal adhesion assembly.
  \newblock Biochemical and Biophysical Research Communications. 2003;307(2):355--361.
  \newblock doi:{10.1016/S0006-291X(03)01165-3}.
  
  \bibitem{mcbeath_cell_2004}
  McBeath R, Pirone DM, Nelson CM, Bhadriraju K, Chen CS.
  \newblock Cell {Shape}, {Cytoskeletal} {Tension}, and {RhoA} {Regulate} {Stem} {Cell} {Lineage} {Commitment}.
  \newblock Developmental Cell. 2004;6(4):483--495.
  \newblock doi:{10.1016/S1534-5807(04)00075-9}.
  
  \bibitem{kassab_biomechanical_2024}
  Kassab GS.
  \newblock Biomechanical {Homeostasis}.
  \newblock In: Kassab GS, editor. Coronary {Circulation}: {Mechanobiology}, {Growth}, {Remodeling}, and {Clinical} {Implications}. Cham: Springer International Publishing; 2024. p. 1--43.
  \newblock Available from: \url{https://doi.org/10.1007/978-3-031-62652-4_1}.
  
  \bibitem{andersen_cell_2023}
  Andersen T, Wörthmüller D, Probst D, Wang I, Moreau P, Fitzpatrick V, et~al.
  \newblock Cell size and actin architecture determine force generation in optogenetically activated cells.
  \newblock Biophysical Journal. 2023;122(4):684--696.
  \newblock doi:{10.1016/j.bpj.2023.01.011}.
  
  \bibitem{weaver_cellular_2016}
  Weaver VM, Gilbert PM.
  \newblock Cellular adaptation to biomechanical stress across length scales in tissue homeostasis and disease.
  \newblock Seminars in cell \& developmental biology. 2016;67:141.
  \newblock doi:{10.1016/j.semcdb.2016.09.004}.
  
  \bibitem{chanduri_cellular_2024}
  Chanduri M, Kumar A, Weiss D, Emuna N, Barsukov I, Shi M, et~al.
  \newblock Cellular stiffness sensing through talin 1 in tissue mechanical homeostasis.
  \newblock Science Advances. 2024;10(34):eadi6286.
  \newblock doi:{10.1126/sciadv.adi6286}.
  
  \bibitem{cacopardo_characterizing_2022}
  Cacopardo L, Guazzelli N, Ahluwalia A.
  \newblock Characterizing and {Engineering} {Biomimetic} {Materials} for {Viscoelastic} {Mechanotransduction} {Studies}.
  \newblock Tissue Engineering Part B, Reviews. 2022;28(4):912--925.
  \newblock doi:{10.1089/ten.TEB.2021.0151}.
  
  \bibitem{byfield_endothelial_2009}
  Byfield FJ, Reen RK, Shentu TP, Levitan I, Gooch KJ.
  \newblock Endothelial actin and cell stiffness is modulated by substrate stiffness in {2D} and {3D}.
  \newblock Journal of Biomechanics. 2009;42(8):1114--1119.
  \newblock doi:{10.1016/j.jbiomech.2009.02.012}.
  
  \bibitem{chaudhuri_substrate_2015}
  Chaudhuri O, Gu L, Darnell M, Klumpers D, Bencherif SA, Weaver JC, et~al.
  \newblock Substrate stress relaxation regulates cell spreading.
  \newblock Nature Communications. 2015;6(1):6365.
  \newblock doi:{10.1038/ncomms7365}.

  \bibitem{isomursuDirected2022}
  Isomursu A, Park K-Y, Hou J, Cheng B, Mathieu M, Shamsan GA, et~al.
  \newblock Directed Cell Migration towards Softer Environments.
  \newblock Nature Materials. 2022;21(9):1081--1090.
  \newblock doi:{10.1038/s41563-022-01294-2}.
  
  \bibitem{ross_physical_2012}
  Ross AM, Jiang Z, Bastmeyer M, Lahann J.
  \newblock Physical {Aspects} of {Cell} {Culture} {Substrates}: {Topography}, {Roughness}, and {Elasticity}.
  \newblock Small. 2012;8(3):336--355.
  \newblock doi:{10.1002/smll.201100934}.
  
  \bibitem{jafarinia_insights_2024}
  Jafarinia H, Khalilimeybodi A, Barrasa-Fano J, Fraley SI, Rangamani P, Carlier A.
  \newblock Insights gained from computational modeling of {YAP}/{TAZ} signaling for cellular mechanotransduction.
  \newblock npj Systems Biology and Applications. 2024;10(1):1--14.
  \newblock doi:{10.1038/s41540-024-00414-9}.
  
  \bibitem{le_devedec_residence_2012}
  Le~Dévédec SE, Geverts B, de~Bont H, Yan K, Verbeek FJ, Houtsmuller AB, et~al.
  \newblock The residence time of focal adhesion kinase ({FAK}) and paxillin at focal adhesions in renal epithelial cells is determined by adhesion size, strength and life cycle status.
  \newblock Journal of Cell Science. 2012;125:4498--4506.
  \newblock doi:{10.1242/jcs.104273}.
  
  \bibitem{logg_automated_2012}
  Logg A, Mardal KA, Wells G, editors.
  \newblock Automated {Solution} of {Differential} {Equations} by the {Finite} {Element} {Method}. vol.~84 of Lecture {Notes} in {Computational} {Science} and {Engineering}.
  \newblock Berlin, Heidelberg: Springer Berlin Heidelberg; 2012.
  \newblock Available from: \url{http://link.springer.com/10.1007/978-3-642-23099-8}.
  
  \bibitem{geuzaine_gmsh_2009}
  Geuzaine C, Remacle JF.
  \newblock Gmsh: {A} 3-{D} finite element mesh generator with built-in pre- and post-processing facilities.
  \newblock International Journal for Numerical Methods in Engineering. 2009;79(11):1309--1331.
  \newblock doi:{10.1002/nme.2579}.
  
  \bibitem{lakkis2013implicit}
  Lakkis O, Madzvamuse A, Venkataraman C.
  \newblock Implicit-explicit timestepping with finite element approximation of reaction-diffusion systems on evolving domains.
  \newblock SIAM Journal on Numerical Analysis. 2013;51(4):2309--2330.
  \newblock doi:{10.1137/120880112}.
  
  \end{thebibliography}
\end{document}